\documentclass[11pt]{amsart}
\usepackage{amsmath}
\usepackage{amssymb}
\usepackage{amscd}
\usepackage{mathrsfs}
\usepackage{graphicx}
\usepackage{psfrag}
\usepackage{times}
\usepackage[all]{xy}
\def\ol#1{\overline{#1}}     
\def\wh#1{\widehat{#1}}   
\def\wt#1{\widetilde{#1}}   

\def\pair#1{\langle #1 \rangle}
\def\makeop#1{\expandafter\def\csname #1\endcsname{\mathop{\mathrm{#1}}\nolimits}} 
\def\makefrak#1{\expandafter\def\csname f#1\endcsname{\mathfrak{#1}}}
\def\makebb#1{\expandafter\def\csname bb#1\endcsname{\mathbb{#1}}}
\def\makebf#1{\expandafter\def\csname b#1\endcsname{\mathbf{#1}}}
\def\makecal#1{\expandafter\def\csname c#1\endcsname{\mathcal{#1}}}
\def\Exp#1{\exp \left[#1\right]}
\makeop{End}\makeop{Ker}\makeop{div}%
\makeop{Re}\makeop{Im}\makeop{Pic}\makeop{hol}
\makeop{Hom}\makeop{Ext}\makeop{Sym}%
\makeop{Spm}\makeop{Spec}\makeop{Spf}%
\makeop{Gal}\makeop{Cont}\makeop{tor}%
\makeop{Im}\makeop{id}\makeop{ab}%
\makeop{pol}\makeop{rec}\makeop{Res}%
\makebb{C}\makebb{Z}\makebb{Q}\makebb{F}\makebb{G}%
\makebb{T}%
\makecal{O}\makecal{L}\makecal{P}\makecal{E}\makecal{M}%
\makecal{R}%
\makefrak{p}\makefrak{a}\makefrak{f}\makefrak{m}\makefrak{c}\makefrak{g}%
\makefrak{n}%
\makefrak{P}%
\makebf{w}\makebf{z}\makebf{c}%
\makebf{N}%

\def\bundlexi{\xi}
\def\elementalpha{\alpha}

\def\smallalpha{\epsilon}

\def\Nfp{N\kern-0.5mm\fp}
\theoremstyle{plain}
	\newtheorem{theorem}{Theorem}[section]
    	\newtheorem{proposition}[theorem]{Proposition}
	\newtheorem{theoremint}{Theorem} 
    \newtheorem{lemma}[theorem]{Lemma}
    \newtheorem{corollary}[theorem]{Corollary}
\theoremstyle{definition}
    \newtheorem{definition}[theorem]{Definition}
    
    \newtheorem{example}[theorem]{Example}
    \newtheorem{remark}[theorem]{Remark}

\newcommand{\C}{\mathbb{C}}

\newcommand{\G}{\mathbb{G}}
\newcommand{\Z}{\mathbb{Z}}
\newcommand{\Q}{\mathbb{Q}}
\newcommand{\R}{\mathbb{R}}

\renewcommand{\labelenumi}{\roman{enumi})}
\begin{document}
\title[Algebraic Theta Functions and Eisenstein-Kronecker Numbers]{Algebraic theta functions and the
$p$-adic interpolation of Eisenstein-Kronecker numbers}
\author{Kenichi Bannai and Shinichi Kobayashi}
\email{bannai@math.nagoya-u.ac.jp} 
\urladdr{http://www.math.nagoya-u.ac.jp/\textasciitilde bannai/ }
\email{shinichi@math.nagoya-u.ac.jp}
\urladdr{http://www.math.nagoya-u.ac.jp/\textasciitilde shinichi/ }
\thanks{Both authors were supported in part by the JSPS postdoctoral fellowship for research abroad.}
\address{Graduate School of Mathematics, Nagoya University}
\address{Furo-cho Chikusa-ku, Nagoya, Japan 464-8602}
\date{December 11, 2007}

\begin{abstract}
	We study the properties of Eisenstein-Kronecker numbers,
	which are related to special values of Hecke $L$-function of imaginary quadratic  fields.
	We prove that the generating function of these numbers
	is a reduced  (normalized or canonical in some literature)
	theta function associated to the Poincar\'e bundle of an elliptic curve.
	We introduce general methods to study the algebraic and $p$-adic properties
	of reduced theta functions for CM abelian varieties.
	As a corollary, when the prime $p$ is ordinary, we give a new construction of the two-variable $p$-adic measure
	interpolating special values of Hecke $L$-functions of imaginary quadratic fields, 
	originally constructed by   Manin-Vishik and Katz.
	Our method via theta functions also gives insight for the case when $p$ is supersingular.
	The method of this paper will be used in subsequent papers to study the precise $p$-divisibility of critical values of 
	Hecke $L$-functions associated to Hecke characters of quadratic imaginary fields
	for supersingular $p$, as well as explicit calculation in two-variables of the $p$-adic elliptic 
	polylogarithm for CM elliptic curves.
\end{abstract}

\setcounter{tocdepth}{1}
\maketitle

\setcounter{section}{-1}
\section{Introduction}
%

\subsection{Introduction}

The Eisenstein-Kronecker-Lerch series is a generalization of the classical real analytic Eisenstein series,
and was reintroduced by Andr\'e Weil in his inspiring book \cite{We1}.  In this paper, we investigate the algebraic and 
$p$-adic properties of the \textit{Eisenstein-Kronecker numbers},  which we define to be the special values of 
the Eisenstein-Kronecker-Lerch series.  These numbers may  be regarded as elliptic 
analogues of the classical generalized Bernoulli numbers.  

Let $\Gamma$ be a lattice in $\bbC$, and let $A$ be the area of the fundamental domain of $\Gamma$ 
divided by $\pi$.  Then Eisenstein-Kronecker numbers for integers $a$, $b$ such that
$b > a+2$ are defined as the sum
$$
	e_{a,b}^*(z_0, w_0) := \sum_{\gamma \in \Gamma \setminus \{ -z_0\}} 
	\frac{(\ol z_0 + \ol \gamma)^a}{(z_0 + \gamma)^b}
	\pair{\gamma, w_0}_\Gamma,
$$
where $z_0$, $w_0 \in \bbC$ and $\pair{z,w}_\Gamma$ is the pairing defined by
$\pair{z,w}_\Gamma := \exp[(z\ol w - w \ol z)/A]$.
The sum converges only when $b > a+2$, but one can give it
meaning  for any $a \geq 0$, $b >0$ by analytic continuation.
Similarly to the fact that generalized Bernoulli numbers may be used to express
special values of Dirichlet $L$-functions, when the lattice $\Gamma$ has complex multiplication, 
Eisenstein-Kronecker numbers may be used to express special values of Hecke $L$-functions
pertaining to the field of complex multiplication.
Being of considerable arithmetic interest, these numbers have been studied by many authors.

We investigate the problem using a new approach, through the systematic use of algebraic theta functions.
The crucial observation for our approach is the fact that the two-variable generating function for
Eisenstein-Kronecker numbers is a certain theta function $\Theta(z,w)$, associated to the Poincar\'e 
bundle of the elliptic curve $E(\bbC) = \bbC/\Gamma$. 
The Poincar\'e bundle is a line bundle on $E \times E^\vee$, 
where $E^\vee$ is the dual of $E$.
Our observation allows us to view the two variables of the generating function as
coming from the elliptic curve and the dual elliptic curve.   This differs from the standard viewpoint, 
originally established by Katz \cite{Ka2}, where the second variable
is regarded as a parameter on the moduli space.
The advantage of our view point is that the theta function $\Theta(z,w)$ is a reduced theta function 
(referred in literature also to as a normalized or canonical theta function) with an algebraic divisor, and
the properties of our generating function may be studied 
through Mumford's theory  of
algebraic theta functions on abelian varieties \cite{Mum5}, applied
to the case when the abelian variety is $E \times E^\vee$.

We study in detail the algebraic and $p$-adic properties of
reduced theta functions, on a general abelian variety with complex multiplication.
Using this theory, we prove the algebraicity of Eisenstein-Kronecker numbers  when the 
corresponding lattice has complex multiplication by the ring of integers of an imaginary  quadratic 
field $K$.  As a corollary, this gives the classical theorem of Damerell concerning
the algebraicity of the critical values of Hecke $L$-functions of imaginary quadratic fields.
Our proof is conceptual, in a sense that the values are algebraic 
because the generating function is an algebraic theta function.
We further apply our method to construct a $p$-adic measure on $\bbZ_p \times \bbZ_p$ which $p$-adically 
interpolates the Eisenstein-Kronecker numbers when $p \geq 5$ is an ordinary prime, i.e. 
when $p$  splits as $(p)= \fp \ol\fp$ in $K$.
Not surprisingly, our $p$-adic measure is related to the measure constructed by Manin-Vishik \cite{MV}
and Katz \cite{Ka2} which were used in the construction of the $p$-adic $L$-function for algebraic 
Hecke characters associated to $K$.  We relate our measure to various $p$-adic measures 
which appear in literature.

\vskip2mm

One application of our approach is the following.
The fact that we have a generating function at our disposal allows us to understand in detail
the $p$-adic properties of Eisenstein-Kronecker numbers,
even for the case when $p$ is supersingular, i.e. when $p$ remains prime in $K$.
Calculations given at the end of this paper reveal that the $p$-adic radius of convergence of the two-variable generating function in the supersingular case has radius of convergence strictly less than one. 
In a subsequent paper \cite{BK}, we build upon our technique to give a construction of  
$p$-adic distributions for inert $p$ previously studied by Boxall \cite{Box1}, \cite{Box2}, 
Schneider-Teitelbaum \cite{ST}, Fourquaux \cite{Fou} and Yamamoto \cite{Yam}, which 
interpolates Eisenstein-Kronecker numbers in one variable.  We then study in two-variables 
the $p$-divisibility of critical values of Hecke $L$-functions associated to Hecke characters
of imaginary quadratic fields for inert $p$, extending previous works of 
Katz \cite{Ka6}, Chellali \cite{Ch} and Fujiwara \cite{Fuj}.
%

As further application, in \cite{BKT}, the construction of $p$-adic distributions in this paper and 
\cite{BK} will be used to study the relation between such distributions and
$p$-adic elliptic polylogarithms for CM elliptic curves, for any prime $p \geq 5$ of good reduction.  
This result gives a $p$-adic analogue of the result of Beilinson and Levin \cite{BL}, which expressed the
Hodge realization of the elliptic polylogarithm in terms of complex Eisenstein-Kronecker series.
This extends previous results of the first author \cite{Ban}, which dealt only with one-variable 
measures in the ordinary case.  Since the $p$-adic elliptic polylogarithms are expected to be of
motivic in origin, and since the $p$-adic distributions interpolate special values of Hecke $L$-functions, 
this result may be interpreted as a $p$-adic analogue of Beilinson's conjecture.

The method of this paper can also be used to investigate any set of invariants which appear as Taylor coefficients of a
reduced theta function on an abelian variety with complex multiplication.  Theta functions for elliptic curves 
were used in the construction of the Euler system of elliptic units, which played an important role in 
Iwasawa theory.  There has been attempts to generalize this result to higher genus
(see for example \cite{dSG}),
but with preliminary success.  Considering that the Poincar\'e bundle exists for abelian varieties,
in light of the conjectured relation between $p$-adic $L$-functions and Euler systems,
we hope our approach via algebraic theta functions will give future insight to this problem.\\

\tableofcontents

\subsection{Overview}
The paper consists of four sections.  In the first section, we introduce Eisenstein-Kronecker numbers
and the Kronecker theta function, which is a reduced theta function
on the Poincar\'e bundle of an elliptic curve.  
Our main theorem is that the Kronecker theta function is the generating function of the Eisenstein-Kronecker numbers.  In the second section, we prove general theorems concerning
the algebraic and $p$-adic properties of reduced theta functions on abelian varieties with complex 
multiplication.  In the third section, we apply this theory to the Kronecker theta function to
construct our $p$-adic measure.  In the final section, we discuss explicit calculations, especially 
in the case of supersingular primes $p$.  

The precise content of each section is as follows.\vskip2mm

Let $\Gamma \subset \bbC$ be a lattice, and let $E$ be an elliptic curve such that $E(\bbC)=\bbC/\Gamma$.
Let  $\Delta$ be the kernel of the multiplication $E \times E \rightarrow E$, and 
consider the divisor
$$
	D: = \Delta - (E \times \{ 0 \}) - (\{ 0 \} \times E)
$$
in $E \times E$.   We define the Kronecker theta function $\Theta(z,w)$ to be the unique reduced 
theta function on $E \times E$ characterized by the property that its divisor is $D$ and 
its residue along $z=0$ is equal to one.   If we identify $E$ with its dual, then
this function is a meromorphic section of the Poincar\'e bundle. 
In general, a theta function is determined by a divisor up to multiplication by a quadratic exponential.  
The reducedness condition removes the ambiguity up to a constant multiple,
and the condition on residue normalizes the constant.

Let $z_0$, $w_0 \in \Gamma \otimes \bbQ$.  We define the algebraic translation by $(z_0, w_0)$
of $\Theta(z,w)$ to be the function
$$
	\Theta_{z_0, w_0}(z,w) := \Exp{ - \frac{z_0 \ol w_0}{A}}%
	 \Exp{ - \frac{z \ol w_0 + w \ol z_0}{A} } \Theta(z_0 + z, w_0 + w),
$$
which is again a reduced theta function.
Such translation is defined for general reduced theta functions on complex tori, and
the extra exponential factors play a role in preserving algebraicity.

The main theorem of the first section of this paper is our key observation, namely,
we prove that the Kronecker theta function $\Theta_{z_0, w_0}(z,w)$ is a two-variable generating function
for the Eisenstein-Kronecker numbers.

\begin{theoremint}[= Theorem \ref{theorem: generating function theorem}]\label{theoremint: one}
	We have the Laurent expansion
	\begin{multline*}
		\Theta_{z_0, w_0}(z,w) = \pair{w_0, z_0} \frac{\delta(z_0)}{z}  
		+ \frac{\delta(w_0)}{w} \\+ \sum_{a \geq 0, b> 0} 
		(-1)^{a+b-1} \frac{e^*_{a, b}(z_0, w_0)}{a! A^a}  z^{b-1} w^a,
	\end{multline*}
	where $\delta(u) = 1$ if $u \in \Gamma$ and $\delta(u)=0$ otherwise.
\end{theoremint}
\noindent
The key result which will be used in the proof of the above theorem
is Kronecker's theorem (Theorem \ref{theorem: Kronecker}).
The normalization of the function $\Theta(z,w)$ differs by an exponential factor
from that of the two-variable Jacobi theta function previously studied 
by Zagier \cite{Zag}.  See \eqref{eq: compare zagier} for the precise relation. The generating function property of the two-variable Jacobi theta function at the origin is given in \cite{Zag} \S 3 Theorem.  

\vskip2mm

In the second section, we work with a general abelian variety $A$, defined over a number field, with complex multiplication (CM).   Fix an algebraic uniformization $\bbC^g/\Lambda \cong A(\bbC)$  compatible with the 
complex multiplication (see \S \ref{subsection:3-1} for the precise statement).
For any divisor $D \subset A$, we may associate a reduced theta
function $\vartheta_s(z)$ with divisor $D$, determined up to a constant multiple, on $A(\bbC) \cong\bbC^g/\Lambda$.
We first prove that, assuming a slight admissibility condition, 
if $D$ is defined over a number field and the  leading term of the Laurent expansion at the origin
of  $\vartheta_s(z)$ is algebraic, then the $g$-variable Laurent expansion of $\vartheta_s(z)$ 
at the origin has algebraic coefficients (Proposition \ref{proposition; algebraicity at the origin}).
We then define an algebraic translation  $U_{v_0} \vartheta_s(z)$ of $\vartheta_s(z)$ by a torsion point $v_0 \in A(\ol\bbQ)$.  This translation may be interpreted by Mumford's theory, and using this theory,
we prove the algebraicity of the Laurent expansion of  $U_{v_0} \vartheta_s(z)$ again at the origin 
by reducing to the case of $\vartheta_s(z)$ (Theorem \ref{theorem: algebraicity}).

In the case of the Kronecker theta function,  we have
$$
	U_{(z_0, w_0)} \Theta(z,w) = \Theta_{z_0, w_0}(z,w).
$$ 
We apply the above theory to Eisenstein-Kronecker numbers for the case when $\Gamma \subset \bbC$ is 
the period lattice of an elliptic curve $E$ with complex multiplication defined over a number field $F$.
Consider an imaginary quadratic field $K$, and suppose
$E$ has complex multiplication by the ring of integers $\cO_K$
of $K$.  If we fix an invariant differential $\omega$ of $E$ defined over $F$,
then we have a uniformization $\bbC/\Gamma\cong E(\bbC)$
such that the pull-back of $\omega$ is $dz$.
The above result and Theorem \ref{theoremint: one} 
gives the algebraicity of the Eisenstein-Kronecker numbers in this case
(Corollary \ref{corollary: AT2}), when $z_0$, $w_0 \in \Gamma \otimes \bbQ$.  
This implies the classical theorem of Damerell (Corollary \ref{corollary: Damerell})
on the algebraicity of the special values of Hecke $L$-functions.

We then prove the $p$-integrality of the Laurent expansion of $\vartheta_s(z)$ (Propositions 
\ref{proposition: integrality of theta}  and \ref{proposition: integrality of translation}), when
$p$ satisfies a certain ordinarity condition for the abelian variety.  This applied to the
case of the Kronecker theta function is as follows.
Let $p \geq 5$ be an ordinary prime, i.e. splits as $(p) = \fp \ol\fp$ in $K$.  We suppose that
the elliptic curve $E$ defined over $F$ has good reduction at a prime $\fP$ over $\fp$,
and  we choose a smooth model $\cE$ of $E$ over $\cO_{F_{\fP}}$, where $F_{\fP}$
is the completion of $F$ at $\fP$.
We denote by $\wh\cE$  the formal group of $\cE$ at the origin, and we let
$\lambda(t)$ be the formal logarithm of $\wh\cE$.

\begin{theoremint}[=Corollary \ref{corollary: for measure}]
	Let $z_0$ and $w_0$ be torsion points of $E(\ol\bbQ)$ of order prime to $p$, and let
	$$
		\wh\Theta_{z_0, w_0}(s,t) := \left.\Theta_{z_0, w_0}(z,w) \right|_{z = \lambda(s), w = \lambda(t)}
	$$
	be the formal composition of the two-variable
	Laurent expansion of $\Theta_{z_0, w_0}(z,w)$ at the origin
	with the power series $z = \lambda(s)$  in $s$ and $w= \lambda(t)$ in $t$.
	Then we have
	$$
		\wh\Theta_{z_0, w_0}(s,t)  - \pair{w_0, z_0} \delta(z_0) s^{-1} -  \delta(w_0) t^{-1}
		\in \cO[[s,t]],
	$$
	where $\cO$ is the ring of integers of a finite extension of $F_\fP$ unramified over $\fP$.
\end{theoremint}

In the third section, we use the above theorem and well-established relation between power series and 
$p$-adic measures 
to construct the
$p$-adic measure $\mu_{z_0, w_0}$ on $\bbZ_p \times \bbZ_p$, interpolating the Eisenstein-Kronecker numbers.
We relate this to various $p$-adic measures used in the construction of $p$-adic $L$-functions interpolating
special values of Hecke $L$-function of imaginary quadratic  fields.  Namely, we compare our measure 
with that of Yager (Theorem \ref{theorem: Yager}) and that of Katz (Theorem \ref{theorem: Katz}).
We also give an explicit construction of the $p$-adic measure defined by Mazur and Swinnerton-Dyer
interpolating twists of Hasse-Weil $L$-function of elliptic curves, in the CM case.\vskip2mm

In the last section, we calculate the Laurent expansions of the Kronecker theta function
for explicit examples in the supersingular case.  
One of the goals of our research was to apply  the methods of
Boxall \cite{Box1} \cite{Box2} and Schneider-Teitelbaum \cite{ST}
to the two-variable power series $\wh\Theta(s,t)$ to construct a two-variable 
$p$-adic $L$-function in this case.  Our calculations show, however,
that this naive strategy is not directly applicable, since the two-variable power series
$\wh \Theta(s,t)$ in the supersingular case is not only non-integral, but has radius of 
convergence strictly less than one. In a subsequent paper \cite{BK}, we refine the theory
of Schneider-Teitelbaum to study the $p$-adic properties of $\wh\Theta(s,t)$ when $p$ is inert.
We then use this theory to study the $p$-divisibility of critical values of Hecke $L$-function associated
to Hecke characters of imaginary quadratic fields.

\subsection{Acknowledgment}

Part of this research was conducted while the first author was visiting the \'Ecole Normale Sup\'erieure at Paris,
and the second author Institut de Math\'ematiques de Jussieu.  The authors would like to thank their
hosts Yves Andr\'e and Pierre Colmez for hospitality.  The authors would also like to thank 
Guido Kings for discussion, and Ahmed Abbes for pointing out some references.
We greatly appreciate Seidai Yasuda for discussion, especially for informing the authors
important formulas  useful for explicit calculations of elliptic curves.  Finally, we would like to thank 
our colleagues Yasushi Komori  and Hyohe Miyachi at Nagoya University for introducing the authors
to computational software.

\section{Kronecker Theta Function}

In this section, we will define and study the basic properties of the Kronecker theta function. 
In \S \ref{subsection:2-1}, we review the definition and properties of the Eisenstein-Kronecker-Lerch series, and we define the Eisenstein-Kronecker numbers as the special values of the Eisenstein-Kronecker-Lerch series. 
Special values of Hecke $L$-functions are expressed in terms of Eisenstein-Kronecker numbers. 
In \S \ref{subsection:2-2}, we review the definition of the Poincar\'e bundle for elliptic curves, and 
we define the Kronecker theta function  $\Theta(z,w; \Gamma)$ as the reduced (or normalized) theta function 
associated to a certain canonical meromorphic section of the Poincar\'e bundle. 
Then we give the relation between Eisenstein-Kronecker-Lerch series and  Kronecker theta functions. 
In \S \ref{subsection:2-3},  we define a theta-theoretic  translation operator $U_{(z_0,w_0)}$ of reduced theta functions. 
We then prove in \S \ref{subsection:2-4} that the Laurent expansion of $U_{(z_0, w_0)}\Theta(z,w; \Gamma)$ at the origin 
 is a generating function of the Eisenstein-Kronecker numbers.

\subsection{Eisenstein-Kronecker Numbers}\label{subsection:2-1}
%
%
We introduce the Eisenstein-Kronecker-Lerch series following Weil \cite{We1}. 
The results of this subsection are contained in \cite{We1}, and we refer the reader there for details.

Let $\Gamma = \bbZ \omega_1 \oplus \bbZ \omega_2$ be a lattice in $\mathbb{C}$ generated
by $\omega_1$ and $\omega_2$ with $\Im(\omega_2/\omega_1) > 0$, and let 
$$
  A(\Gamma)  = \frac{1}{\pi} \Im( \omega_2 \overline{\omega_1})
=  \frac{1}{2 \pi i}(\omega_2  \overline{\omega}_1 -   \omega_1 \overline{\omega_2}).
$$
We define a pairing for $z$, $w \in \mathbb{C}$ by 
$\pair{z,w}_{\Gamma} :=  \Exp{(z \ol{w}  - w \ol{z})/A(\Gamma)}$.  We will freely use the following 
properties of this pairing.
\begin{enumerate}
	\item $\pair{z,w}_\Gamma = \pair{-w,z}_\Gamma= \pair{w,z}_\Gamma^{-1}$,
	\item $\pair{az, w}_\Gamma = \pair{z, \ol a w}_\Gamma $ for any $a \in \C$,
	\item We have $z \in \Gamma$ if and only if $\pair{z, \gamma}_\Gamma = 1$ for 
	all $\gamma \in \Gamma$.
\end{enumerate}

\begin{definition}(\cite{We1} VIII \S 12)\label{definition: Eisenstein-Kronecker-Lerch series} 
 Let $a$ be an integer $\geq 0$.  
For $z$, $w \in \bbC\setminus \Gamma$,  we define the 
  \textit{Eisenstein-Kronecker-Lerch series} $K_a(z,w,s; \Gamma)$ by
  \begin{equation} \label{equation: definition of Eisenstein-Kronecker}
     K_a(z,w,s; \Gamma) = {\sum}_{\gamma \in \Gamma} 
    \frac{(\ol{z} + \ol{\gamma})^a}{|z + \gamma|^{2s}} \pair{\gamma, w}_{\Gamma}, \quad
    (\Re s > a/2 + 1).
  \end{equation}
  For a fixed $z_0$, $w_0 \in \bbC$,  we define $K^*_a(z_0,w_0,s; \Gamma)$ by
  \begin{equation} \label{equation: definition of Eisenstein-Kronecker*}
     K_a^*(z_0,w_0,s; \Gamma) = {\sum}^*_{\gamma \in \Gamma} 
    \frac{(\ol{z}_0 + \ol{\gamma})^a}{|z_0 + \gamma|^{2s}} \pair{\gamma, w_0}_{\Gamma}, \quad
    (\Re s > a/2 + 1),
  \end{equation}
   where $\sum^*$  means the sum taken over all $\gamma \in \Gamma$ other than $-z_0$
   if $z_0$ is in $\Gamma$.
  The series on the right  converges absolutely for $\Re s > a/2+1$.
 \end{definition}

\begin{remark}
	The series $K_a(z,w,s; \Gamma)$ defines a $\mathscr{C}^\infty$ function for $z$, $w \in \bbC\setminus \Gamma$.
	In this case, we have $K_a(z,w,s; \Gamma)=K^*_a(z,w,s; \Gamma)$.  If we were to let 
	$K_a(z,w,s; \Gamma) = K^*_a(z,w,s;\Gamma)$ for $z$, $w \in \C$, then this function 
	would not be continuous when $z \in \Gamma$ or $w \in \Gamma$.
	In order to avoid using non-continuous functions, we take the convention that $z$, $w \notin \Gamma$ for $K_a(z,w,s; \Gamma)$,
	and we fix $z_0$, $w_0 \in \C$ when considering the series $K_a^*(z_0,w_0,s; \Gamma)$.
\end{remark}

When there is no fear of confusion, we will often omit the $\Gamma$ from the notation and
simply denote $K_a(z,w,s)$ for $K_a(z,w,s; \Gamma)$ and $A$ for $A(\Gamma)$, etc.  
The function $K_a^*(z_0,w_0,s)$ is known to satisfy the following properties.

\begin{proposition}[\cite{We1} \rm{VIII} \S 13] \label{proposition: properties of K}
  Let $a$ be an integer $\geq 0$. 
  \begin{enumerate}
     \renewcommand{\theenumi}{\roman{enumi}}
     \renewcommand{\labelenumi}{(\theenumi)}
     \item The function $K_a^*(z_0,w_0,s)$  for $s$ continues meromorphically to a function on the whole $s$-plane, 
     with possible poles only at 
     $s=0$  $($if $a=0$, $z_0 \in \Gamma${}$)$ and at 
     $s=1$ $($if $a=0$, $w_0 \in \Gamma${}$)$.
     \item $K_a^*(z_0,w_0,s)$ satisfies the functional equation
      \begin{equation} \label{equation: functional equation for K}
         \Gamma(s) K_a^*(z_0,w_0,s) = A^{a+1-2s} \Gamma(a+1-s) K_a^*(w_0,z_0,a+1-s) \pair{w_0,z_0}.
      \end{equation}
 \end{enumerate}
\end{proposition}	

 We note that in the case $a>0$, the analytic
continuation and the functional equation of $K_a^*(z_0,w_0,s)$ follows from the representation
\begin{equation} \label{equation: representation as integral} 
   \Gamma(s) K_a^*(z_0,w_0,s) =  I_{a}(z_0, w_0, s) 
              + A^{a+1-2s} I_{a}(w_0,z_0, a+1-s)
   \pair{w_0,z_0},
\end{equation}
where $I_{a}(z_0, w_0, s)$ is the integral
\begin{equation*}
 I_{a}(z_0, w_0, s) = \int_{A^{-1}}^{\infty} \theta_a(t,z_0,w_0) t^{s-1}  dt
\end{equation*}
for
\begin{equation} \label{equation: theta for analytic continuation}
   \theta_a(t,z_0,w_0) = \sum_{\gamma \in \Gamma} \Exp{-t|z_0+ \gamma|^2} 
 (\overline{z}_0 + \overline{\gamma})^a  \pair{\gamma, w_0},  
\end{equation}
which follows from the functional equation of $ \theta_a(t,z_0,w_0)$; namely, 
$$\theta_a(t,z_0,w_0)=(At)^{-a-1}\theta_a(A^{-2}t^{-1},w_0,z_0)\pair{w_0,z_0}.$$

The Eisenstein-Kronecker-Lerch series for various integers $a$ and $s$  are related 
by the following differential equations.

\begin{lemma} \label{lemma: calculation of differentials}
   Let $a$ be an integer $>0$.
   The function $K_a(z,w,s)$ satisfies differential equations
   \begin{equation*}
       \begin{aligned}  
           \partial_{z} K_a(z,w,s) &= -s K_{a+1}(z,w,s+1) \\
           \partial_{\ol{z}}K_a(z,w,s) &= (a-s)K_{a-1}(z,w,s) \\
           \partial_{w} K_a(z,w,s) &= - A^{-1}( K_{a+1}(z,w,s) - \ol{z} K_{a}(z,w,s) ) \\
           \partial_{\ol{w}}K_a(z,w,s) &= A^{-1} (K_{a-1}(z,w,s-1) - z K_{a}(z,w,s) ).
        \end{aligned}
   \end{equation*}
\end{lemma}

\begin{proof}
	From the definition,  the statement is true for $\Re s > a/2+1$.
	The statement for general $s$ is obtained by analytic continuation.
\end{proof}

%

We now give the definition of the Eisenstein-Kronecker numbers.

\begin{definition}  \label{definition: EK number}
	Let $z_0$, $w_0 \in \C$.
	For any integer $a \geq 0$, $b >0$,  we define the \textit{Eisenstein-Kronecker number}  
	$e^*_{a,b}(z_0, w_0; \Gamma)$ by
	$$
		e_{a,b}^*(z_0, w_0 ; \Gamma) := K_{a+b}^*(z_0, w_0, b; \Gamma).
	$$
	We will  omit $\Gamma$ from the notation if there is no fear of confusion.
\end{definition}

The numbers $e_{a,b}^*(\Gamma) := e_{a,b}^*(0,0; \Gamma)$ 
are those that appear in \cite{We1} (VI \S 5, VIII \S 15).
For $b > a \geq 0$, we have by definition 
$$e_{a,b}^*(\Gamma) = \lim_{s \rightarrow 0} \sum_{\gamma \in \Gamma \setminus \{ 0 \}} 
\frac{
\overline{\gamma}^a}{ \gamma^b|\gamma|^{2s}}.$$  

Let $\ol\bbQ$ be the algebraic closure of $\bbQ$, and fix once and for all an embedding  $i_{\infty}: \ol\bbQ \hookrightarrow \bbC$.    Let $K$ be an imaginary  quadratic  field with ring of integers $\cO_K$.  We consider $\cO_K$ and any ideal 
$\ff$ of $\cO_K$ to be lattices in $\bbC$ through $i_\infty$.
One may express the special values of Hecke $L$-functions of $K$ by using the Eisenstein-Kronecker numbers.   
Let $\varphi$ be an algebraic Hecke character of ideals  in $K$, with conductor $\ff$ and
infinity type $(m, n)$.  In other words, $\varphi$ is a homomorphism from 
$I_K(\ff)$, the group of the fractional ideals of $\cO_K$ 
prime to $\ff$, to $\ol{\bbQ}^{\times}$ of the form 
$\varphi((\elementalpha)) = \varepsilon(\elementalpha) \elementalpha^{m}
 \ol{\elementalpha}^n$ for 
some finite character $$\varepsilon : (\cO_K/\ff)^{\times} \rightarrow \ol{\mathbb{Q}}^{\times}$$ 
if $\elementalpha$  in $\cO_K$ is prime to $\ff$. 
Then the \textit{Hecke $L$-function} associated to $\varphi$ is defined by 
	$$
		L_{\ff}(\varphi,s) = \sum_{(\fa, \ff) = 1} \frac{\varphi(\fa)}{N \fa^s}	$$
	where the sum is over integral ideals of $\cO_K$ prime to $\ff$, 
	and this series is absolutely convergent if  
	$\Re s > (m+n)/2+1$. 
The Hecke $L$-function is known to have a meromorphic continuation to the whole 
complex $s$-plane and a functional equation.  
The function $L_{\ff}(\varphi,s)$ has a pole at $s=s_0$ if and only if 
the conductor of $\varphi$ is trivial, $m=n$ and $s_0=(m+n)/2+1$. 
The Hecke $L$-function may be expressed in terms of Eisenstein-Kronecker series as follows. 

\begin{proposition}\label{proposition: EK and L}
Let $\ff$  be an integral ideal of $K$.  
Let $I_K(\ff)$ be the subgroup of the ideal group of $K$ consisting of 
ideals prime to $\ff$, and let $P_K(\ff)$ be the subgroup of $I(\ff)$ consisting of 
 principal ideals $(\alpha)$ such that $\alpha \equiv 1 \mod \ff$.  
Let $a$ be a non-negative integer.  \\
	i) 
	Let $\varphi$ be an algebraic Hecke character of $K$ of conductor $\ff$ 
	and infinity type $(1,0)$. 
Let $\alpha_\fa$ be an element of $\fa^{-1} \cap P_K(\ff)$. 
	Then  we have
		\begin{align*}
		L_{\ff}(\ol{\varphi}^{a}, s) &= 
		\sum_{\fa \in I_K(\ff)/P_K(\ff)} 
		K^*_{a}(\varphi(\elementalpha_\fa \fa), 0, s; \varphi(\fa) \fa^{-1}\ff ). 
	\end{align*}
		ii) Let $w_\ff$ be 
	 the number of roots of unity in $K$ congruent to $1$ modulo $\ff$. 
	For a Hecke character $\varphi$ of $K$ of conductor $\ff$ and infinity type $(a,b)$ ($a\not= b$),  we have 
$$L_{\ff}(\varphi, s)=
  \frac{1}{w_\ff}
		\sum_{\fa \in I_K(\ff)/P_K(\ff)}  \frac{{\varphi}(\fa)}{N \fa^{s}}
	K^*_{|b-a|}(\alpha,  0,s-\mathrm{min}\{a,b\} ;(\fa^{-1}\ff)^{\delta} ),
	$$	
	where $\delta \in \mathrm{Gal}(\C/\R)$ and 
	$\delta$ is trivial if and only if $b-a>0$. 
\end{proposition}

\begin{proof} 
	For i), we first consider the case $s > a/2+1$.  Then we have
	\begin{equation*}\begin{split}
		L_{\ff}(\ol{\varphi}^{a}, s)& = 
		\sum_{(\fa, \ff)=1} \frac{\ol{\varphi}(\fa)^{a}}{N \fa^s} = \frac{1}{w_\ff}
		\sum_{\fa \in I_K(\ff)/P_K(\ff)} \frac{\ol{\varphi}(\fa)^{a}}{N \fa^{s}}
		 \sum_{\elementalpha \in \fa^{-1}\cap P_K(\ff)} \frac{\ol{\varphi}(\elementalpha)^{a}}{|\elementalpha|^{2s}}\\
				&= \frac{1}{w_\ff} \sum_{\fa \in I_K(\ff)/P_K(\ff)} 
		\sum_{\gamma \in \fa^{-1}\ff} 
		\frac{(\ol{\varphi(\elementalpha_\fa\fa)+\varphi(\fa)\gamma})^{a}}
		{|\varphi(\elementalpha_\fa\fa) + \varphi(\fa)\gamma|^{2s}}. 
		\end{split}\end{equation*}
		The existence of $\varphi$ of type $(1,0)$ shows $w_\ff=1$, hence our formula.  The statement for
		general $s$ follows by analytic continuation. For ii), suppose that $b>a$.  Then we have 
	\begin{equation*}\begin{split}
			L_{\ff}({\varphi}, s)=
		\sum_{(\fa, \ff)=1} \frac{{\varphi}(\fa)}{N \fa^s} =  \frac{1}{w_\ff}
		 \sum_{\fa \in I_K(\ff)/P_K(\ff)} \frac{{\varphi}(\fa)}{N \fa^{s}}
		 \sum_{\elementalpha \in \fa^{-1}\cap P_K(\ff)} \frac{
		 {\ol{\elementalpha}}^{b-a}}{|\elementalpha|^{2(s-a)}}
			\end{split}\end{equation*}
		as desired.
\end{proof}

From the above proposition, the study of special values of Hecke $L$-functions of $K$ is reduced 
to the study of Eisenstein-Kronecker numbers.  

\subsection{Poincar\'e bundle and reduced theta Function}\label{subsection:2-2}
%
%

Let $V$ be a finite dimensional  complex vector space and $\Lambda$ a lattice in $V$.  
It is known that every holomorphic 
line bundle $\mathscr{L}$ on $\bbT =V/\Lambda$ is  a quotient of the trivial bundle 
$V \times \mathbb{C}$
over $V$ by the action of $\Lambda$ of the 
form 
\begin{equation}
	\gamma \cdot ( \bundlexi, v) = (e_\gamma (v) \bundlexi, v + \gamma), 
	\qquad
	(\gamma \in \Lambda, (\bundlexi, v) \in \mathbb{C} \times V)
\end{equation} 
for some $1$-cocycle $\gamma \mapsto e_\gamma$ of  $\Lambda$ 
with values in the group $\cO^* (V)$ of invertible holomorphic 
functions on $V$. 
This correspondence extends to an isomorphism of groups 
$\mathrm{Pic}(\bbT) \cong H^1(\Lambda, \cO^* (V))$. 
A meromorphic section $s: \bbT \rightarrow \mathscr{L}$  is given by $v \mapsto (\vartheta(v), v)$ for some meromorphic function
$\vartheta(v)$ on $V$ satisfying $\vartheta(v+\gamma) = e_\gamma(v) \vartheta(v)$. 
We call $\vartheta$ the  theta function corresponding to $s$. 


The isomorphism classes of  holomorphic line bundles on  $\bbT$ are classified by the following theorem.

\begin{theorem}[Appell and Humbert]
	The group $\Pic(\bbT)$ of isomorphism classes of holomorphic line bundles on $\bbT$ is isomorphic to the 
	group of pairs $(H, \alpha)$, where
	\begin{enumerate}
		\renewcommand{\theenumi}{\roman{enumi}}
		\renewcommand{\labelenumi}{(\theenumi)}
		\item $H$ is an Hermitian form on $V$.
		\item $E = \Im H$ takes integral values on $\Lambda \times \Lambda$.
		\item $\alpha : \Lambda \rightarrow U(1) : = \{ z \in \mathbb{C} \mid |z|=1 \}$ is a map
		such that $\alpha(\gamma+ \gamma') = \exp(\pi i E(\gamma, \gamma')) \alpha(\gamma) \alpha(\gamma')$.
	\end{enumerate}	
\end{theorem}
The correspondence is given by associating to a pair $(H, \alpha)$ the line bundle 
$\mathscr{L}(H, \alpha)$ whose cocycle is given by
$$e_\gamma(v):=\alpha(\gamma) \exp\left(\pi H(v, \gamma) +\frac{\pi}{2}%
	 H(\gamma,\gamma) \right) .$$
	We put 
	$\pair{v_1,v_2}_\mathscr{L}:=
	\Exp{2 \pi i E(v_1, v_2)}.$

\begin{definition}
	We call any theta function associated to a 
	(meromorphic) section of a line bundle of the form $\mathscr{L}(H, \alpha)$ a 
	\textit{reduced theta function}.  Namely, a reduced theta function is a meromorphic function $\vartheta$ on $V$ 
	which satisfies 
	$$\vartheta(v+\gamma)=\alpha(\gamma) \exp\left(\pi H(v, \gamma) +\frac{\pi}{2}%
	 H(\gamma,\gamma) \right) \vartheta(v)$$	
	for all $v \in V$ and $\gamma \in \Lambda$. 
	The term \textit{normalized} or \textit{canonical} theta function is used in some literature.
\end{definition} 

For a divisor $D$ of  $\bbT$, there exists,  up to constant, a unique meromorphic 
section $s$ with divisor $D$ of the line bundle $\mathscr{L}(D)$ associated to the invertible sheaf 
$\cO_{\bbT}(D)$. 
We take $(H, \alpha)$ such that $\mathscr{L}(D)\cong \mathscr{L}(H,\alpha)$ and 
we have a reduced theta function associated to $s$.   Hence we can associate to $D$ a reduced theta 
function determined up to a constant multiple whose divisor $D$.
Without the condition of reducedness, a theta function with divisor  $D$ is determined
up to multiplication by a trivial theta function, which is an exponential of a quadratic form. 

Now we give two important examples of reduced theta functions. 

\begin{example}\label{robert}
	Let $E$ be an elliptic curve such that $E(\bbC) = \bbC/\Gamma$.
	We consider the line bundle $\mathscr{L}=\mathscr{L}([0])$ associated to the divisor 
	$[0]$ of an $E$. Then the pair $(H, \alpha)$ corresponding to 
	$\mathscr{L}$ is 
	$$
		H (z_1, z_2) = \frac{z_1 \ol z_2}{\pi A}, 
	$$
	and  $\alpha : \Gamma \rightarrow \{ \pm 1 \}$ is such that
	$\alpha(u) = -1$ if $u \not\in 2 \Gamma$ and $\alpha(u)=1$ otherwise. 
	This is checked, for example, by constructing a reduced theta function 
	associated to $[0]$ explicitly as follows. 

	 Let $\sigma(z)$ be the Weierstrass $\sigma$-function
	\begin{equation}\label{product expansion of sigma}
		\sigma(z) = z \prod_{\gamma \in \Gamma \setminus \{ 0 \}}
		\left( 1 - \frac{z}{\gamma} \right) \Exp{ \frac{z}{\gamma} + \frac{z^2}{2 \gamma^2}}.
	\end{equation}
	This function is  known to satisfy the transformation formula 
	$$
		\sigma(z+ u) =\alpha(u) \Exp{ \eta(u)  \left( z + \frac{u}{2} \right)} \sigma(z)
	$$
	for any $u \in \Gamma$, %
	where $\eta(z) := A^{-1} \,\overline{z} +e^*_2\, z$ where  $e^*_2:=e_{0,2}^*(\Gamma)$.  

	The function  $\sigma(z)$ is not reduced and we define the function $\theta(z)$ to be
	$$
		\theta(z) = \Exp{-\frac{ e_{2}^*  }{2}z^2}  \sigma(z).
	$$
	Then $\theta(z)$ is a holomorphic function on $\bbC$, having simple zeros at the points in $\Gamma$ and satisfies the 
	transformation formula 
	\begin{equation}\label{equation: transformation formula for theta}
				\theta(z + \gamma) = \alpha(\gamma)
				\exp\left(\pi H(z, \gamma) +\frac{\pi}{2}%
		 H(\gamma,\gamma) \right) 
				\theta(z)
			\end{equation}
			for $u \in \Gamma$. 
			Namely, $\theta(z)$ is a reduced theta function  
			associated to the divisor $[0]$. 
	The function $\theta(z)$ was used by Robert to construct the Euler system of elliptic units.  
	As Robert pointed out  in \cite{Rob2} \S 1, the function $\theta(z)$ is characterized as 
	the reduced theta function  associated to the divisor $[0]$ normalized by the condition $\theta'(0)=1$. 
\end{example}

\begin{example} [The Kronecker theta function]\label{Theta}
Given a complex torus $\bbT$ and its dual torus $\bbT^\vee$, the 
Poincar\'e bundle is a line bundle $\mathscr{P}$ on $\bbT \times \bbT^\vee$ giving the isomorphism
$$
	\bbT^\vee \xrightarrow\cong \Pic^0(\bbT),
$$
defined by mapping $w \in \bbT^\vee$ to the line bundle on $\bbT$ obtained as the
restriction of $\mathscr{P}$ to $\bbT = \bbT \times \{ w \} \hookrightarrow \bbT \times \bbT^\vee$.

On $E\times E^\vee$,  the Poincar\'e bundle $\mathscr{P}$  is a line bundle characterized by the properties 
i) the restriction $\mathscr{P} \vert_{\{0\}\times E^\vee}$ is trivial and ii)
for all $w \in  E^\vee$,
the restriction  to $E\times \{w\}$ represents the element of 
$\mathrm{Pic}^0(E)$ given by $w$ under the isomorphism $\mathrm{Pic}^0(E) \cong E^\vee$. 
We identify $E$ with $E^\vee$ by the canonical polarization given by the divisor $[0]$ of $E$.  
Then  $\mathscr{P}$ is what is called a Mumford bundle on $E\times E$, namely, is of the form
$$
	\mathscr{P}=m^*\mathscr{L} \otimes p_1^*\mathscr{L}^{-1} \otimes p_2^*\mathscr{L}^{-1},
$$ 
where $\mathscr{L}=\mathscr{L}([0])$ and   
the morphisms $m$, $p_1$ and $p_2$ are the multiplication,  the first and the second projections from 
$E \times E$ to $E$.   In other words, $\mathscr{P}$ is the line bundle associated to the invertible sheaf 
$\cO_{E\times E} (D)$ where $D$ is the divisor $\Delta-(\{0\}\times E)-(E\times\{0\})$ and 
$\Delta$ is the kernel of the multiplication morphism $E \times E \rightarrow E$. 
We have $\mathscr{P} \cong \mathscr{L}(H_\mathscr{P}, \alpha_\mathscr{P})$ where 
$$H_\mathscr{P}=m^*H_\mathscr{L}- p_1^*H_\mathscr{L}- p_2^*H_\mathscr{L}, 
\qquad \alpha_\mathscr{P}=m^*\alpha_\mathscr{L}\cdot p_1^*\alpha_\mathscr{L}^{-1}\cdot p_2^*\alpha_\mathscr{L}^{-1}$$
for $\mathscr{L}=\mathscr{L}(H_\mathscr{L}, \alpha_\mathscr{L})$. 
These are explicitly given as 
\begin{equation} \begin{split}\label{Pairing and character of P}
      H_\mathscr{P}((z_1, w_1), (z_2, w_2)) 
&= \frac{z_1 \overline{w}_2 + w_1 \overline{z}_2}{\pi A}, \\
      \alpha_\mathscr{P}(u,v) &= \exp \left(  \frac{u \overline{v}- \overline{u} v}{2A} \right). 
\end{split} \end{equation}
	Hence, the cocycle $\Gamma \times \Gamma \rightarrow \cO^* (V), 
	(u,v) \mapsto e_{ (u,v)}(z,w)$ associated to $\mathscr{P}$ is given by 
	$$
		e_{(u,v)}(z,w)= \exp \left[ \frac{u \overline{v}}{A} \right]
			\exp \left[ \frac{z \overline{v} + w \overline{u}}{A}\right] .
	$$
	Hence for $u, v \in \Gamma$, any reduced theta function $\vartheta(z,w)$ for $\mathscr{P}$ satisfies the 
 	transformation formula 
	\begin{equation}\label{transformation formula for poincare bundle}
			\vartheta(z+u, w+v) = \exp \left[ \frac{u \overline{v}}{A} \right]
			\exp \left[ \frac{z \overline{v} + w \overline{u}}{A}\right] \vartheta(z,w). 
	\end{equation}

Since $H_\mathscr{P}$ is not definite, the Poincar\'e bundle $\mathscr{P}$ 
does not have any non-zero holomorphic section.  
This fact may be proved easily as follows. 
\begin{lemma}\label{global section of P}
	Any holomorphic function $f(z,w)$ on $\bbC \times \bbC$ satisfying the transformation formula
	\begin{equation}\label{equation: transformation formula}
		f(z+u, w+v) = \Exp{\frac{u \ol v}{A}} \Exp{\frac{z \ol v + w \ol u}{A}} f(z,w)
	\end{equation}
	for any $u$, $v \in \Gamma$ is identically equal to \textit{zero}.
\end{lemma}

\begin{proof}
			For any point $(z_0, w_0) \in \mathbb{C}
			\times \mathbb{C}$ and $u \in \Gamma$, we have
			$$
				|f(z_0+u, w_0-u)| = |f(z_0, w_0)| \Exp{ - \frac{|u|^2 + \Re((z_0-w_0)\ol{u})}{A} }.
			$$
			The right hand side goes to $0$ as $|u| \rightarrow \infty$.  Applying the
			maximum principal to the holomorphic function $f(z_0 + s, w_0 -s)$ with
			respect to the variable $s$, we obtain $f(z_0, w_0) = 0$ as desired.
\end{proof}

Therefore there are no holomorphic sections of $\mathscr{P}$.  However,  there exists a canonical meromorphic section 
$s_D$, namely, a section corresponding to the divisor $D$.   Since $s_D$ is of the form $m^*s\otimes p_1^*s^{-1} \otimes p_2^*s^{-1}$ 
for a section $s$ of $\mathscr{L}_E([0])$ corresponding the divisor $[0]$, the reduced theta function 
 corresponding to $s_D$ is, up to constant, 
$$\Theta(z,w):=\frac{\theta(z+w)}{\theta(z)\theta(w)}.$$ 
We write $\Theta(z,w)$ as $\Theta(z,w; \Gamma)$ if we want to specify the lattice, 
and we call it {\it the Kronecker theta function for $\Gamma$}. 
By definition, the Kronecker theta function is characterized as 
the reduced theta function associated to $D$ whose residue 
along $z=0$ is equal to $1$. 
The relation between the theta function $\Theta(z,w)$ for the lattice $\Gamma = \bbZ\tau \oplus  \bbZ$
and the two-variable Jacobi theta function $F_\tau(z,w)$ defined in \cite{Zag} \S 3  is given by
\begin{equation}\label{eq: compare zagier}
	\Theta(z,w) = \exp(zw/A)F_\tau(z,w).
\end{equation}

\begin{proposition}\label{proposition: sections of the  poincare bundle}
i) For a positive real number $t$, the real analytic function
 $$\Exp{{z\overline{w}}/{A}}\theta_a(t,z,w)$$ \
is a $\mathscr{C}^\infty$-section of $\mathscr{P}$. 
$($see  \eqref{equation: theta for analytic continuation} for the definition of $\theta_a(t,z,w))$. 
In particular, the function $\Exp{{z\overline{w}}/{A}}K_a(z,w,s)$ 
is also a $\mathscr{C}^\infty$-section of $\mathscr{P}$  on an open set of $\C \times \C$.  \\
ii) The function $\Exp{{z\overline{w}}/{A}}K_1(z,w,1)$ is holomorphic on $\C \times \C$ except for simple
poles at the divisor corresponding to $\{0\} \times E$ and $E \times \{0\}$ with residue 
$1$ along $z=0$ and $w=0$. 
In particular, the function $\Exp{{z\overline{w}}/{A}}K_1(z,w,1)$ 
is a meromorphic section of $\mathscr{P}$. 
\end{proposition}
\begin{proof}
	i) It suffices to show that $\vartheta_a(t,z,w)=\Exp{{z\overline{w}}/{A}}\theta_a(t,z,w)$
	satisfies 
	$$
  	  	\vartheta_a(t, z+u, w+v) = e_{(u,v)}(z,w) \vartheta_a(t, z,w)
  	$$ 
	for any $u$, $v \in \Gamma$. 
	This is checked by direct calculation.  \\
	ii) We put $f(z,w)= \Exp{{z\overline{w}}/{A}}K_1(z,w,1)$.  
		By Lemma \ref{lemma: calculation of differentials}, we have for any integer $a>0$ 
	$$
		\partial_{\ol{z}} K_a(z,w,s) = (a-s) K_{a-1}(z,w,s).
  	$$
  	Substituting $a=s=1$, we see that $\partial_{\ol{z} } f(z,w) = 0$, hence $f(z,w)$ is holomorphic in $z$. 
	By the functional equation (\ref{equation: functional equation for K}), we have $f(w,z)=f(z,w)$. 
	Hence we also have  $\partial_{\ol{w} } f(z,w) = \partial_{\ol{w} } f(w,z)= 0$. 	
	We show that $f(z,w)$ has simple poles. 
	By the integral expression (\ref{equation: representation as integral}), 
	it has poles possibly at $z \in \Gamma$ or $w \in \Gamma$, and 
	by the transformation formula of $f(z,w)$ we may assume that 
	$z=0$ or $w=0$. 
	Then by the integral expression (\ref{equation: representation as integral}), the function 
	\begin{equation}\label{eq: calc integral}
 		K_1(z,w,1) - \int_{A^{-1}}^{\infty}  \Exp{-t|z|^2} \overline{z} dt 
		-  \pair{w,z}\int_{A^{-1}}^{\infty}  \Exp{-t|w|^2} 
		 \overline{w} dt 
	\end{equation}
	is analytic at the origin, hence  the function $zwf(z,w)$ is bounded around the origin. 
	The calculation of the residue follows by calculating the integrals in \eqref{eq: calc integral}.
\end{proof}

We now compare the Kronecker theta function $\Theta(z,w)$ with the Eisenstein-Kronecker-Lerch
series $K_1(z,w,1)$.  We obtain  the following theorem, essentially known to 
Kronecker  (for example, see Weil (\cite{We1} VIII, \S 4,  p.71, (7))). 

\begin{theorem}\label{theorem: Kronecker}
The Kronecker theta function $\Theta(z,w)$ associated to 
the divisor  $D=\Delta-(\{0\}\times E)-(E\times\{0\})$ on $E\times E$ is 
related to the Eisenstein-Kronecker-Lerch series $K_1(z,w,1)$ as follows.
$$\Theta(z,w)=\frac{\theta(z+w)}{\theta(z)\theta(w)}
=\Exp{\frac{z\overline{w}}{A}}K_1(z,w,1).$$
\end{theorem} 
\begin{proof}
By Proposition \ref{proposition: sections of the  poincare bundle}, 
the functions $\Theta(z,w)$ and $\Exp{\frac{z\overline{w}}{A}}K_1(z,w,1)$ are both meromorphic sections of
$\mathscr{P}$ and have the same simple poles with the same residues. 
	Hence the function 
	$$\frac{\theta(z+w)}{\theta(z)\theta(w)}
-\Exp{\frac{z\overline{w}}{A}}K_1(z,w,1)$$
	defines a holomorphic section of $\mathscr{P}$, which is zero by Lemma \ref{global section of P}. 
	\end{proof}
		
\end{example}

\subsection{Translations of reduced theta functions}\label{subsection:2-3}
%

Let $\mathbb{T}=V/\Lambda$ be a complex torus.  
We consider the translation $\tau_{v_0} : \mathbb{T} \rightarrow \mathbb{T}, 
v \mapsto v+v_0$ by  a point $v_0 \in V$.  
Since the naive translation of a reduced function 
$\vartheta(v) \mapsto \vartheta(v+v_0)$ does not preserve 
the reducedness of theta functions,
we would like to define a notion of translation for reduced theta functions. 

Let $\mathscr{L}(H, \alpha)$ be the line bundle on $\mathbb{T}$ associated to $(H, \alpha)$. 
Then it is know that 
\begin{equation}\label{transisom}
	\tau_{v_0}^*\mathscr{L}(H, \alpha)  \xrightarrow{\cong} \mathscr{L}(H, \alpha \cdot \nu_{v_0}) 
\end{equation}
where $\nu_{v_0} : \Lambda \rightarrow U(1)$ is the 
character $\nu_{v_0}(\gamma) = \exp(2\pi i E(v_0, \gamma))$. 
Hence for a meromorphic section $s$ of  $\mathscr{L}(H, \alpha )$, we have 
a reduced theta function for $\mathscr{L}(H, \alpha \cdot \nu_{v_0})$ 
corresponding to the section $\tau^*_{v_0}s$. 
However, the choice of the isomorphism (\ref{transisom}) gives ambiguity by a constant multiple
in determining the reduced theta function corresponding to $\tau^*_{v_0}s$.

{\it We now assume that $\alpha : \Lambda \rightarrow U(1)$ is the restriction of 
a certain map $V \rightarrow \C^\times$ to $\Lambda$.} We also denote this map by $\alpha$ 
but the relation
$$\alpha(v_0+v_1)=\alpha(v_0)\alpha(v_1)\pair{v_0, v_1/2}_\mathscr{L}$$
might not be valid for non-lattice points $v_0, v_1$.
Since the reduced theta functions of  $\mathscr{L}(H, \alpha )$ satisfy 
$\vartheta(v+u)=e_{u}(v)  \vartheta(v)$ for $u \in \Lambda$, 
it would be natural to define the reduced theta function for $\mathscr{L}(H, \alpha \cdot \nu_{v_0})$ 
corresponding to the section $\tau^*_{v_0}s$ as  
\begin{align*}
U_{v_0}\vartheta_s(v)&:=e_{v_0}(v)^{-1} \vartheta_s(v+v_0) \\
&=\alpha(v_0)^{-1} \Exp{- \pi H(v, v_0) - \frac{\pi}{2} H(v_0, v_0) } 
	  \vartheta_s(v+v_0)
	\end{align*}
for the reduced theta function $\vartheta_s$  corresponding to $s$. 
One can check that $U_{v_0}\vartheta_s(v)$ is in fact 
a reduced theta function for $\mathscr{L}(H, \alpha \cdot \nu_{v_0})$.
We call $U_{v_0}$ the translation by $v_0$ of reduced theta functions
of $\mathscr{L}(H, \alpha )$ (however, it depends on the 
choice of $\alpha: V \rightarrow \C^\times$.)  
We can also consider 
another translation of reduced theta function by 
$$U_{v_0}^{\cM}\vartheta_s(v):=\Exp{- \pi H(v, v_0) - \frac{\pi}{2} H(v_0, v_0) } 
	\vartheta_s(v+v_0)$$
(it does not need  on a map $\alpha: V \rightarrow \C^\times$.)  

We will see that these translations preserve certain algebraic properties. 
In fact, the translation $U_{v_0}^{\cM}$ is obtained by the theory 
of algebraic theta functions of Mumford. 
For a symmetric algebraic line bundle $\mathscr{L}$ with a fixed isomorphism 
$\mathscr{L} \cong  \mathscr{L}(H, \alpha )$ over $\C$, Mumford's theory gives an algebraic way to 
fix an appropriate  isomorphism $\tau_{v_0}^*\mathscr{L} \cong \mathscr{L}(H, \alpha \cdot \nu_{v_0})$. 
Then the reduced theta function corresponding to $\tau^*_{v_0}s$ by this isomorphism is 
$U_{v_0}^{\cM}\vartheta_s(v)$. 
See \S \ref{subsection:3-2} for details. 
For a convenience, we denote some basic properties of these translations.  

\begin{proposition}\label{proposition; transformations }
Let $v_0, v_1$ be elements of $V$, and let 
$$
	e_{v_0}^{\cM}(v)=\Exp{\pi H(v, v_0) + \frac{\pi}{2} H(v_0, v_0) } .
$$ 
Then we have
\begin{enumerate}
\item \quad $e_{v_0+v_1}^{\cM}(v)=\pair{v_0,v_1/2}_\mathscr{L}\,e_{v_0}^{\cM}(v+v_1)e_{v_1}^{\cM}(v).$ \\
\item  \quad $ U_{v_0}^{\cM}\vartheta_s(v+v_1)
	=\pair{v_0,v_1/2}_\mathscr{L}\,e_{v_1}^{\cM}(v)
	U_{v_0+v_1}^{\cM}\vartheta_s(v). $\\
\item  \quad $U_{v_1}^{\cM} \circ U_{v_0}^{\cM}=\pair{v_0,v_1/2}_\mathscr{L} \,U_{v_0+v_1}^{\cM}
=\pair{v_0,v_1}_\mathscr{L}\,U_{v_0}^{\cM} \circ U_{v_1}^{\cM}.$ \\

$\! \! \! \!\! \! \! \!\! \! \!\! \! \! \!\! \! \!\! \! \!\! \! \!\! \! \!$ Suppose that $\alpha : \Lambda \rightarrow U(1)$ 
is the restriction of 
a map $\alpha : V \rightarrow \C^\times$.  We put 
$$\chi(v_0, v_1):=\alpha(v_0+v_1)\alpha(v_0)^{-1}\alpha(v_1)^{-1}\pair{v_0, v_1/2}_\mathscr{L}.$$ 
$\! \! \! \!\! \! \! \!\! \! \!\! \! \! \!\! \! \!\! \!$ Then we have 

\item \quad $e_{v_0+v_1}(v)=\chi(v_0, v_1)\,e_{v_0}(v+v_1)e_{v_1}(v).$ \\
\item  \quad $ U_{v_0}\vartheta_s(v+v_1)
	=\chi(v_0, v_1)\,e_{v_1}(v)
	U_{v_0+v_1}\vartheta_s(v). $\\
\item  \quad $U_{v_1}\circ U_{v_0}=\chi(v_0, v_1) \,U_{v_0+v_1}
=\pair{v_0, v_1}_\mathscr{L}\,U_{v_0}\circ U_{v_1}.$
\end{enumerate}
\end{proposition}
\begin{proof}
i) is proven by direct calculation. 
These are proved by direct calculation. ii) follows from i) and iii) follows from ii). 
iv)-vi) follows from i)-iii) by comparing  $e_{v_0}^{\cM}(v)$ and $e_{v_0}(v)$. 
\end{proof}

Now consider
the Mumford bundle 
$$\mathscr{M}:=m^*\mathscr{L}\otimes p_1^*\mathscr{L}^{-1}\otimes p_2^*\mathscr{L}^{-1}$$ on 
$V/\Lambda \times V/\Lambda$. If we write 
$\mathscr{M}=\mathscr{L}(H_\mathscr{M}, \alpha_\mathscr{M})$, 
then we have $\alpha_\mathscr{M}(v,w)=\pair{v, w/2}_\mathscr{L}$, which also have 
meaning on $V \times V$. Hence we consider a map 
$\alpha: V \times V \rightarrow \C^\times, (v,w) \mapsto \pair{v, w/2}_\mathscr{L}$ to 
define the translation $U_{v_0}$ for $v_0 \in V \times V$ of reduced theta functions on $V\times V$.  
We also have 
$$
\pair{(v_0,w_0), (v_1,w_1)}_\mathscr{M}=\pair{v_0, w_1}_\mathscr{L}\pair{w_0,v_1}_\mathscr{L},$$
$$
\chi((v_0,w_0), (v_1,w_1) )=\pair{v_0, w_1}_\mathscr{L}.$$ 

If there is no fear of confusion, we put $\vartheta(v):=\vartheta_s(v)$, 
$\vartheta_{v_0}(v):=U_{v_0}\vartheta(v)$ and 
$\vartheta^{\cM}_{v_0}(v):=U^{\cM}_{v_0}\vartheta(v)$ for simplicity. 
We also put 
$$\varTheta(v,w):=\frac{\vartheta(v+w)}{\vartheta(v)\vartheta(w)}$$
and for $v_0, w_0 \in V/\Lambda$ we put 
$$\varTheta_{v_0, w_0}(v,w) :=U_{(v_0,w_0)} \varTheta(v,w)
.$$  
By Proposition \ref{proposition; transformations } vi), we have 
$$\varTheta_{v_0+\gamma, w_0+\gamma'}(v,w)=\pair{w_0, \gamma}_\mathscr{L} 
\varTheta_{v_0, w_0}(v,w) $$
for $\gamma, \gamma' \in \Lambda$. We define $\varTheta^{\cM}_{v_0, w_0}(v,w)$ similarly.

In particular, the case of the Poincar\'e bundle of elliptic curves, 
namely the case of our  Kronecker theta function $\Theta(z,w)$, 
the translation is given explicitly as 
\begin{equation}\label{equation: definition of modified theta}
		\Theta_{z_0, w_0}(z,w) =	\Exp{- \frac{z_0 \ol{w}_0}{A}} \Exp{ - \frac{z \ol{w}_0 + w \ol{z}_0}{A} }%
		 \Theta(z+ z_0, w + w_0).
	\end{equation}

We next prove a distribution property for the Kronecker theta function, which will be 
important for later calculation.  We first begin with a lemma.

\begin{lemma} \label{lemma, u-op} 
For an element $\gamma \in \Gamma$, we have 
$$\lim_{z \rightarrow -u+\gamma}(z+u-\gamma)\Theta_{u,v}(z,w)=
\pair{v, \gamma}\Exp{\frac{(\ol {\gamma}-\ol{u}) w}{A}}$$
and 
$$\lim_{w \rightarrow -v+\gamma}(w+v-\gamma)\Theta_{u,v}(z,w)=
\pair{u, \gamma-v}\Exp{\frac{z (\ol \gamma-\ol v)}{A}}.$$ 
\end{lemma}
\begin{proof}
This follows from direct calculations. 
\end{proof}

\begin{proposition}[Distribution relation]\label{proposition; distribution}
Let $\mathfrak{a},\mathfrak{b}$ be integral ideals of  $\cO_K$ such that 
$(\mathfrak{a}\mathfrak{b}, \overline{\mathfrak{b}})=1$. 
Let $\epsilon \in \cO_K$ be such that 
$\epsilon \equiv 1 \mod \mathfrak{a}\mathfrak{b}$ and 
$\epsilon \equiv 0 \mod \overline{\mathfrak{b}}$. 
Then 
\begin{align*}
\sum_{\alpha \in {\mathfrak{a}}^{-1}\Gamma/\Gamma, \;
\beta \in {\mathfrak{b}}^{-1}\Gamma/\Gamma}
\pair{\epsilon \alpha,  w_0}_\Gamma\,
&\Theta_{z_0+\epsilon\alpha, w_0+\epsilon \beta}(z,w ; 
\Gamma)\\
&=N(\mathfrak{a}\mathfrak{b}) \Theta_{N\!\mathfrak{a}\,z_0, N\!\mathfrak{b}\, w_0}(N\!\mathfrak{a}\,z, N\!\mathfrak{b}\, w; \overline{\mathfrak{a}\mathfrak{b}}\Gamma). 
\end{align*}
\end{proposition}
\begin{proof}
First we note that $\pair{\epsilon \alpha,  w_0}_\Gamma\Theta_{z_0+\epsilon\alpha, w_0+\epsilon\beta}(z,w; \Gamma)$ 
does not depend on a 
choice of the representative
 $\alpha \in {\mathfrak{a}}^{-1}\Gamma/\Gamma$ and $\beta \in {\mathfrak{b}}^{-1}\Gamma/\Gamma$. 
 By considering the action of $U_{ (z_0,  w_0)}$, 
we may assume that $z_0=w_0=0$. 

We show that the both sides have the same transformation formula 
with respect to $\overline{\mathfrak{a}\mathfrak{b}}\Gamma$. 
For  $u, v \in \overline{\mathfrak{a}\mathfrak{b}}\Gamma$, we have 
$$\Theta(N\!\mathfrak{a}(z+u), N\!\mathfrak{b} (w+v); \overline{\mathfrak{a}\mathfrak{b}}\Gamma)
=e_{u,v}(z,w; \Gamma)\Theta(N\!\mathfrak{a}\,z, N\!\mathfrak{b}\, w; \overline{\mathfrak{a}\mathfrak{b}}\Gamma). 
 $$ 
On the other hand, 
$$U_{(u,v)}\Theta_{\epsilon\alpha, \epsilon \beta}(z,w;\Gamma)
=\pair{\epsilon\alpha, v}_\Gamma\pair{\epsilon\beta, u}_\Gamma
\Theta_{\epsilon\alpha, \epsilon \beta}(z,w;\Gamma)=\Theta_{\epsilon\alpha, \epsilon \beta}(z,w;\Gamma).$$
Hence we have 
$\Theta_{\epsilon\alpha, \epsilon \beta}(z+u,w+v;\Gamma)
=e_{u,v}(z,w; \Gamma)\Theta_{\epsilon\alpha, \epsilon \beta}(z,w;\Gamma)$.

Next we show that both sides have the same poles with the same residues. 
The function on the left hand side  has simple poles at most on $(z,w)$ where 
$z=-\epsilon \alpha_0 +\gamma$ or $w=-\epsilon \beta_0 +\gamma$ for some 
$\alpha_0 \in \fa^{-1}$, $\beta _0 \in \mathfrak{b}^{-1}$ and 
$\gamma \in \Gamma$. By the functional equation, it suffices to calculate the residues 
for $z=-\epsilon \alpha_0 +\gamma$. By Lemma \ref{lemma, u-op}
the  value 
$$\lim_{z \rightarrow -\epsilon \alpha_0+\gamma} (z+\epsilon \alpha_0-\gamma)\sum_{\alpha \in {\mathfrak{a}}^{-1}\Gamma/\Gamma, \;
\beta \in {\mathfrak{b}}^{-1}\Gamma/\Gamma}
\Theta_{\epsilon\alpha, \epsilon \beta}(z,w ; 
\Gamma)$$ is equal to 
\begin{equation*}
\sum_{
\beta \in{\mathfrak{b}^{-1}}\Gamma/\Gamma} \pair{\epsilon \beta, \gamma}_\Gamma 
\exp\left[\frac{(\ol \gamma-\ol\epsilon \ol\alpha_0) w}{A(\Gamma)}\right]=
\begin{cases}
N(\mathfrak{b})\exp\left[\frac{(\ol \gamma-\ol\epsilon \ol\alpha_0)w}{A(\Gamma)}\right] \; &\text{if} \; \gamma \in \ol{\mathfrak{b}}\,\Gamma, \\
\quad 0 \quad &\text{otherwise.} 
\end{cases}
\end{equation*}
Hence the left hand side has a pole at $z \in (\epsilon \mathfrak{a}^{-1}+\ol{\mathfrak{b}})\Gamma =\overline{\mathfrak{b}}\mathfrak{a}^{-1}\Gamma$. It is straightforward to see that 
the function on the right hand side  has the same poles with the same residues. 
Hence the difference of these functions defines a holomorphic 
section of a non-ample line bundle, which should be zero (See Lemma \ref{global section of P}). 
\end{proof}

\subsection{Generating Function for Eisenstein-Kronecker numbers}\label{subsection:2-4}
%

	The main result of this subsection is as follows.

\begin{theorem}\label{theorem: generating function theorem}%
	Let $z_0$, $w_0 \in \bbC$.
	Then the Laurent expansion of $\Theta_{z_0, w_0}(z,w)$  at $(z,w) = (0,0)$ is given by
	\begin{multline*}
		\Theta_{z_0, w_0} (z,w)%
		=   \pair{w_0, z_0}  \delta(z_0) z^{-1} + %
		\delta(w_0) w^{-1} \\
		 + \sum_{a \geq 0, b > 0} %
		 (-1)^{a+b-1} \frac{e^*_{a,b}(z_0, w_0)}{a! A^a}%
		z^{b-1} w^a,
		\end{multline*}
	where $\delta(z) = 1$ if $z \in \Gamma$ and $\delta(z) = 0$ otherwise. 
	 In other words,  $\Theta_{z_0, w_0}(z,w)$ is the generating function for the 
	Eisenstein-Kronecker numbers $e^*_{a, b}(z_0,w_0)$.
\end{theorem}%

The fact that $K_1(z,w,1)$ is a (non-holomorphic) generating function of the Eisenstein-Kronecker numbers 
was already observed by Colmez and Schneps \cite{CS}. \vskip2mm

The proof of Theorem \ref{theorem: generating function theorem} is simple 
if $z,w \in \C\setminus \Gamma$ and essentially it  follows from Lemma 
\ref{lemma: calculation of differentials} and Theorem \ref{theorem: Kronecker}. 
To reduce the general case to 
this case, we introduce auxiliary functions. 

For any integer $a \geq 1$ and subset $\Gamma' \subset \Gamma$, we let
\begin{equation}\label{equation: expression sum}
	\theta_a(t,z,w; \Gamma') = \sum_{\gamma \in \Gamma'} \Exp{-t|z + \gamma|^2} (\ol z + \ol \gamma)^a
	 \pair{\gamma, w}_\Gamma 
\end{equation}
and put 
\begin{equation}\label{equation: expression integral}
	I_a(z,w,s ; \Gamma') = \int_{A^{-1}}^\infty \theta_a(t,z,w,\Gamma') t^{s-1} dt
\end{equation}
and  $\wt I_a(z,w,s; \Gamma') = \Exp{-  \ol z w/A} I_a(z,w,s; \Gamma')$ where 
$A=A(\Gamma)$.  

\begin{lemma}
For an integer $a >0$, the function
$\wt I_{a}(z, w,s ; \Gamma' )$ is analytic at $(z,w)=(z_0,w_0)$ if 
$-z_0  \notin \Gamma'$, and we have 
\begin{equation}\begin{split}
	\partial_z \wt I_{a}(z, w,s ; \Gamma' ) &
	= - \wt I_{a+1}(z,w, s+1 ; \Gamma') \\
	\partial_w \wt I_{a}(z,w, s ;  \Gamma') 
	&= - A^{-1} \wt I_{a+1}(z,w, s; \Gamma'). 
\end{split}\end{equation}
\end{lemma}
\begin{proof}
By the assumption for $z_0$, derivations commute with the 
integral symbol in (\ref{equation: expression integral}). 
Then the assertion is straightforward. 
\end{proof}

For subsets $\Gamma_1, \Gamma_2$ of $\Gamma$, we let  
\begin{multline}\label{equation: integral expression A}
	(b-1)! \wt K_{a, b}(z,w; \Gamma_1, \Gamma_2)\\
	:= \wt I_{a+b}(z,w, b; \Gamma_1) 
	+ A^{a-b+1}\wt I_{a+b}(w,z, a+1; \Gamma_2). 
\end{multline}
(Compare with formula (\ref{equation: representation as integral})).
 
 \begin{lemma}\label{lemma: K**} 
The function $\wt K_{a, b}(z,w; \Gamma_1, \Gamma_2)$ 
is analytic  at $(z,w) = (z_0, w_0)$ if $-z_0 \notin \Gamma_1$ and $-w_0 \notin \Gamma_2$ 
and satisfies 
\begin{equation}\begin{split}\label{equation: expression8}
	\partial_z \wt K_{a, b}(z, w; \Gamma_1, \Gamma_2 ) &
	= -b \wt K_{a, b+1}(z,w; \Gamma_1, \Gamma_2) \\
	\partial_w \wt K_{a,b}(z,w; \Gamma_1, \Gamma_2) 
	&= - A^{-1} \wt K_{a+1, b}(z,w; \Gamma_1, \Gamma_2). 
\end{split}\end{equation}
Moreover, for $\Gamma_1=\Gamma \setminus \{  - z_0 \}$ and 
$\Gamma_2=\Gamma \setminus \{  - w_0 \}$, 
we have 
\begin{equation}\label{equation: expression9}
	\wt K_{a,b}(z_0, w_0; \Gamma_1, \Gamma_2) 
	= \Exp{-  \ol z_0 w_0/A} e_{a, b}^*(z_0, w_0).
\end{equation}
\end{lemma}
\begin{proof}
The derivative formulae follows directly from the previous lemma. 
We consider the case $\Gamma_1=\Gamma \setminus \{  - z_0 \}$ and 
$\Gamma_2=\Gamma \setminus \{  - w_0 \}$. 
If  $a > 0$,  we have 
$
	\theta_a(t, z_0, w_0; \Gamma) = \theta_a(t, z_0, w_0; \Gamma_1)
$
and
$
	\theta_a(t, w_0, z_0; \Gamma) = \theta_a(t, w_0, z_0; \Gamma_2). 
$
Hence by the integral formula
(\ref{equation: representation as integral})  for $K_{a+b}^*(z,w,s;\Gamma)$,
we have 
\begin{align*}
	(b-1)!\wt K_{a, b}(z_0,w_0; & \Gamma_1, \Gamma_2)\\
	&=\wt I_{a+b}(z_0,w_0, b; \Gamma) 
	+ A^{a-b+1}\wt I_{a+b}(w_0,z_0, a+1; \Gamma)\\
	& =(b-1)! \Exp{- \ol z_0  w_0/A} K_{a+b}^*(z_0,w_0,b;\Gamma).
\end{align*}
as desired.
\end{proof}

\begin{lemma}\label{lemma: expression5} 
We let $\Gamma_1=\Gamma \setminus \{  - z_0 \}$ and 
$\Gamma_2=\Gamma \setminus \{  - w_0 \}$.
Then we have 
	\begin{equation}\label{equation: expression1}
	\Theta_{z_0, w_0}(z, w)= u_1
	\Exp{\frac {\ol z_0 w_0}{A}} \wt K_{0,1}(z+z_0, w+w_0; \Gamma, \Gamma)
	\end{equation} 
	and 
	\begin{multline}
		\wt K_{0,1}(z+z_0, w+w_0; \Gamma, \Gamma) =  u_2  \Exp{- \frac{z_0 \ol w_0}{A}}
		 \frac{\delta(z_0)}{z} \\+ u_3  \Exp{- \frac{\ol z_0 w_0 }{A}}
		\frac{\delta(w_0)}{w} 
		+\wt K_{0,1}(z+z_0, w+w_0; \Gamma_1, \Gamma_2)
			\end{multline}
where $u_i=u_i(z,w, \ol z, \ol w)$ ($i=1,2,3$) are real analytic functions for $z,w, \ol z, \ol w$ 
such that  $u_i(z,w, 0, 0)=1$. 
\end{lemma}

\begin{proof}
Since  
$$\wt K_{0,1}(z, w; \Gamma, \Gamma)
= \Exp{- {\ol{z} w}/{A}}K_{0,1}(z, w, 1; \Gamma),$$  we may take  
$$u_1(z,w, \ol z, \ol w)=\Exp{ \frac{(z + z_0) \ol w + (w + w_0) \ol z}{A} }.$$
If  $z_0 \in \Gamma$, we have 
	 $$
		 \wt I_1(z, w, 1; \Gamma) \, =\,   \wt I_1(z, w, 1; \{ - z_0 \}) \,+\,  \wt I_1( z, w, 1; \Gamma_1)
$$
	and 
	\begin{align*}
		\wt I_1  (z + z_0,  & w+ w_0, 1; \{ - z_0 \} )  \\
		&= \frac{1}{z}\Exp{-\frac{(\ol z + \ol z_0)(w + w_0)}{A}} \Exp{- \frac{z \ol z}{A}} 
		\pair{w, z_0}.
	\end{align*}
Hence we may take 
$u_2(z,w, \ol z, \ol w)=\Exp{\frac{z_0 \ol w_0}{A}} \wt I_1  (z + z_0,   w+ w_0, 1; \{ - z_0 \} )$. 
Similarly, we may take 
$u_3(z,w, \ol z, \ol w)=\Exp{\frac{w_0 \ol z_0}{A}} \wt I_1  (w + w_0,   z+z_0, 1; \{ - w_0 \} )$.
	\end{proof}

\begin{proof}[Proof of Theorem \ref{theorem: generating function theorem}]
	By  the previous  Lemma, we have
	\begin{multline}\label{equation: expression6}
		\Theta_{z_0, w_0}(z,w) = v_1 \pair{w_0, z_0} \frac{\delta(z_0)}{z} + v_2 \frac{\delta(w_0)}{w} \\
		 + v_3 \Exp{\frac{ \ol z_0 w_0}{A}} \wt K_{0,1}(z+z_0, w+w_0; \Gamma_1, \Gamma_2),
	\end{multline}
	where $v_i=v_i(z,w, \ol z, \ol w)$ ($i=1,2,3$) are real analytic functions for $z,w, \ol z, \ol w$ 
at the origin such that  $v_i(z,w, 0, 0)=1$. 
We take the differential
	$\partial^{b-1}_z \partial^a_w$ of 
	both sides of \eqref{equation: expression6} and substitute 
	$z=w=\ol z=\ol w=0$. Since the operators 
	$\partial_z, \partial_w$ commute with the evaluation  $\ol z = \ol w = 0$ 
	for real analytic functions with variable $z, w, \ol z, \ol w$ at the origin, 
	we may substitute $\ol z=\ol w=0$ first, and 
	apply \eqref{equation: expression8}.  
	Then we have
	\begin{multline*}
		\partial_z^{b-1} \partial_w^a \left( \Theta_{z_0, w_0}(z,w) -
		 \pair{w_0, z_0} \delta(z_0) z^{-1} - \delta(w_0) w^{-1}\right)
		  \\ = \frac{(-1)^{a+b-1}(b-1)!}{A^a}\Exp{ \frac{ w_0 \ol z_0}{A}} 
		 \wt K_{a,b}(z+z_0, w+w_0; \Gamma_1, \Gamma_2). 
	\end{multline*}
	Our assertion now follows from  
	 \eqref{equation: expression9}.
\end{proof}

\section{Algebraicity and $p$-Integrality on CM abelian varieties}\label{section:3}

The purpose of this section is to study algebraic and $p$-adic properties of 
the Taylor coefficients of reduced theta functions for CM abelian varieties at torsion points. 
In particular, the algebraicity and the $p$-adic integrality of 
Eisenstein-Kronecker numbers are proved.  In \S \ref{subsection:3-1}, we prove that 
the Taylor coefficients of reduced theta functions for CM abelian varieties 
at the origin are algebraic.
Then, in \S \ref{subsection:3-2}, we will prove the corresponding statement at  
torsion points using the theory of algebraic 
theta functions constructed by Mumford to reduce to the case at the origin.  
As a corollary, we obtain Damerell's 
theorem concerning the algebraicity of the special values of the Hecke $L$-function.    
In  \S \ref{subsection:3-3} and \S \ref{subsection:3-3},  we prove
at ordinary primes  the $p$-adic integrality of 
the Taylor coefficients of reduced theta functions for CM abelian varieties at torsion points.
In particular, for $p \geq 5$ that splits as $(p) = \fp \ol\fp$ in $\cO_K$,  
the Laurent expansion of  $\Theta_{z_0, w_0}(z,w)$ with respect to the formal parameter of the elliptic curve has $\fp$-integral coefficients. 
Finally, in  \S \ref{subsection:3-5}, we give the relation between our result and 
the theory of $p$-adic theta functions by Norman. 

\subsection{Algebraicity  of reduced theta functions 
 at the origin}\label{subsection:3-1}
%

In this section, we let $K$ be a CM field of degree $2g$ 
and let $K_0$ be its  totally real subquadratic extension $K_0$ of $K$. 
We fix a CM type $\Phi$ of $K$, namely, 
a set of embeddings $\phi_i: K \hookrightarrow \C$ ($i=1, \dots, g$) such that  
$(\phi_i)$ and its complex conjugates $(\ol \phi_i)$ form all the embeddings of $K$ into $\C$. 
We embed $K$ into $\C^g$ by $\Phi$, that is, $\Phi(a)=(\phi_1(a), \dots, \phi_g(a))$. 
For simplicity, we assume that $K$ is normal over $\Q$. 

Let $(A, \iota)$ be a pair consisting of an abelian variety $A$ defined over a number field $F$ and an embedding 
$\iota: K \hookrightarrow \mathrm{End}(A) \otimes \Q$. 
We assume that $F$ contains $K$ and   
for every $x \in K$, the morphism $\iota(x)$ is also 
defined over $F$. Moreover, for simplicity, we assume that 
$\iota(K) \cap \mathrm{End}(A)=\iota(\cO_K)$ and $\mathrm{dim} A=g.$ 
We also assume that $(A, \iota)$ is of type $(K, \Phi)$ over $F$, namely, 
the representation of $K$ on the space of invariant differential $1$-forms on $A$ over $F$, 
which is obtained through $\iota$, is equivalent to $\Phi$. 
By definition, there exists a non-zero invariant differential $1$-form 
$\omega_i$ such that $\iota(a)^*\omega_i=\phi_i(a)\omega_i$ for all $a \in K$.
Let $\pi$ be a uniformization 
$\C^g/\Lambda \rightarrow A(\C)$ such that the pull back of $\omega_i$ is 
equal to $dz_i$, where $z_i$ is the $i$-coordinate of the canonical basis of $\C^g$. 
If we identify $A(\C)$ with $\C^g/\Lambda$ through $\pi$, then the endomorphism $\iota(a)$ is 
the multiplication by $\phi_i(a)$ on each component.
Then  there exists a fractional ideal $\fa$ and 
an element $\Omega_\Phi={}^t(\Omega_1, \dots, \Omega_g) \in \C^g$ whose component-wise multiplication 
gives $\Lambda=\Phi(\fa) \Omega_\Phi$.
Let $\mathscr{L}$ be a line bundle on $A$. Suppose that  
$\pi^*\mathscr{L} \cong \mathscr{L}(H, \alpha)$.  
We assume that $E=\mathrm{Im}\, H$ is $\Phi$-admissible, namely, 
it satisfies $E(\iota(a)z,w)=E(z, \iota(a^c)w)$ 
for a non-trivial element of 
$c\in \mathrm{Gal}(K/K_0)$. 
It is known that if $A$ is simple, then any non-zero $E$ is automatically 
$\Phi$-admissible. (cf. Shimura \cite{Sh3}, Theorem 4 or Lang \cite{Lan3}, Theorem 4.5). 
We say simply that $\mathscr{L}$ is $\Phi$-admissible if $E$ is 
$\Phi$-admissible. 

\begin{proposition}\label{proposition; algebraicity at the origin}
Let  $\mathscr{L}=\mathscr{L}(H, \alpha)$ be a $\Phi$-admissible line bundle on $A$, and 
assume that $\alpha^N=1$ for some integer $N$. 
$($For example, if $\mathscr{L}$ is symmetric, then we have $\alpha^2=1$. $)$
Let $s$ be a meromorphic section of $\mathscr{L}$ whose 
divisor $D$ is defined over $F$. 
Let $f$ be a rational function over $F$ which defines the 
Cartier divisor $D$ in a neighborhood of the origin.  
Let $\vartheta_s(z)$ be a reduced theta function on 
$A(\C)=\C^g/\Phi(\fa)\Omega_\Phi$ corresponding to $s$. 
We assume that $(\vartheta_s/f)(0) \in F$. 
Then the Taylor  coefficients of $f^{-1}(z) 
\vartheta_s(z)$ are in $F$. 
\end{proposition}
\begin{proof}
First suppose that $\vartheta_s(z)$ is periodic and holomorphic at the origin.
Then there exists a non-zero constant $a \in \C$ such that $a \vartheta_s(z)$ is a rational function 
defined over $F$.  Since the derivation $\partial/\partial z_i$ with respect to $\omega_i$ is defined over $F$,
the function $a \vartheta_s(z)$ has Taylor  coefficients in $F$. 
However, since $(\vartheta_s/f)(0)\in F$, the element $a$ is in $F$.  Hence the proposition is proved in this case.
We reduce the general case to this case. 
In particular, we may assume that $f=1$, namely, that 
$\vartheta$ is holomorphic at the origin and $\vartheta(0)=1$. 

Let $p$ be a prime number prime to $N$ that split completely over $K$ and 
let $\fp$ be a prime ideal of $K$ over $p$. 
Then the Chinese remainder theorem shows that there exists 
an element  $a$  in $N\cO_K$ such that 
$\phi_i(a) \in \fp$ but  $\ol \phi_i(a) \notin \fp$ for all $i$.  
We consider the function 
	$$
		F_{a}(z) =\vartheta_s(\iota(a) z)/ \vartheta_s(\iota({a^c}) z). 
	$$
	Then from the transformation formula and $\Phi$-admissibility, 
	the function $F_{a}(z)$ is periodic with respect to the lattice $\Lambda$,
	and its  divisor is defined over $F$.  
	Since $F_a(0)=1$, as we have seen, $F_a(z)$ has the Taylor coefficients 
	in $F$. We write the Taylor expansion $\vartheta_s(z)=\sum c_n z^n$ using multi-index $n$, and we consider 
	the equality
	$$
	      \vartheta_s(\iota(a) z) = \vartheta_s(\iota({a}^c) z)F_{a}(z).
	$$
	Comparing the degree $n=(n_1, \dots, n_g)$ term of the above equality, we have
	\begin{equation*}
		\left(\prod_{i=1}^g \phi_i(a)^{n_i}\right) c_{n} =
		\left( \prod_{i=1}^g \ol\phi_i(a)^{n_i}\right) c_{n} + u(\{c_k\}),
	\end{equation*}
	where $u(\{z_k\})$ is a polynomial over $F$ of several variables $\{z_k\}$ 
         where $k$ runs through 
	multi-indices whose total degree is less than $n$. 
	Then by the assumption on $a$, the element 
	$\prod_{i=1}^g \phi_i(a)^{n_i}$ is different from  
	$\prod_{i=1}^g \ol\phi_i(a)^{n_i}$ if $n\not=0$. 
	Since $\vartheta(0)=1$, we have $c_0 = 1$, hence by induction,
 	we have $c_n \in F$. 
\end{proof}

Now suppose that $K$ is an imaginary  quadratic  field and 
let $E$ be an elliptic curve over  $F$ with complex multiplication by 
$\cO_K \cong\End_{F} (E)$. 
We fix a Weierstrass model $ y^2 = 4 x^3 - g_2x - g_3$ of 
$E$ over $F$ and let $\Gamma$ be the period lattice  
associated to the invariant differential $\omega = dx/y$.  
Let $\pi$ be the complex uniformization
\begin{equation}\label{equation: complex uniformization}
	\pi\colon\mathbb{C} / \Gamma \xrightarrow{\cong} E(\mathbb{C}),
\end{equation}
defined by $\pi(z) = \left(\wp(z),\wp'(z) \right)$, where 
$\wp(z)$  is the Weierstrass $\wp$ function associated to $\Gamma$.  
We fix an isomorphism
$
	\iota :  \cO_K \cong \End_{F} (E)
$
such that $\iota(\alpha)^* \omega = \alpha \omega$.

\begin{corollary}\label{corollary: AT1}
Let $\theta(z)$ be the theta function of Example \ref{robert}. 
Then the Taylor coefficients of $\theta(z)$ are elements in $F$. 
Moreover, the Laurent expansion of $\Theta(z,w)$ at $z=w=0$ has coefficients in $F$.
\end{corollary}

\begin{proof}
We put  $t(z):=-2x/y=-2\wp(z)/\wp'(z)=z+\cdots$, and we apply 
Proposition \ref{proposition; algebraicity at the origin} 
by taking $t$ to be $f$.  The last statement  follows from Theorem 
\ref{theorem: Kronecker}. 
	\end{proof}

\subsection{Algebraicity of the translation $U_{v_0}$ }\label{subsection:3-2}
%

Mumford's theory of algebraic theta functions is a systematic method to reduce the investigation
of the properties of theta functions at torsion points to investigation at the origin.
In this subsection, we  relate the translation $U_{v_0}$ in \S \ref{subsection:2-3} 
to the translation which appears in Mumford's theory.
Using this theory and the fact that the Laurent expansion of $\varTheta(z,w)$ 
at the origin has algebraic coefficients,  we prove that the Laurent expansion of $\varTheta_{z_0,w_0}(z,w)$ at 
the origin also has algebraic coefficients. 
First we recall Mumford's  theory of algebraic theta functions. 
Let $A$ be an abelian variety  defined over a field $F$ (with or without CM).
Let $\mathscr{L}$ be a symmetric  line bundle on $A$, and suppose $P$ is a torsion point in $A(\ol F)$.
In order to compare properties of sections of $\mathscr{L}$ at the origin and at $P$, we would like
to construct an isomorphism $\tau_P^* \mathscr{L} \cong \mathscr{L}$, where $\tau_P$ is the
translation on $A$ by $P$.  However, such an isomorphism does not exists in general and 
we also have to consider not the translation of  a point of $A$ but of a 
point of the universal cover of $A$.
In order to circumvent this problem, Mumford proceeds as follows.

Let $V(A)$ be the set of the systems $(P_n)_{n \in \mathbb{N} }$ such that 
$P_n \in A(\overline{F})$ and $mP_{mn}=P_n$ and$NP_1=0$ for some non-zero 
integer $N$.  In other words, the group 
$V(A)$ is the adelic Tate module of $A$, namely,  
$$V(A)=\{ (\alpha_p) \in 
{\prod}_p V_p(A) \;|\; \text{all but finitely many $\alpha_p$ belong to $T_p(A)$}\}. 
$$
We may think of $V(A)$ as an algebraic version of  
the universal cover of $A$.   For $P=(P_n)_{n \in \mathbb{N} } \in V(A)$ and $NP_1=0$, Mumford
constructed the following.
\begin{proposition}[Mumford]
	Let $\mathscr{L}$ be a symmetric  line bundle on $A$. 
	For any integer $n$, suppose $n : A \rightarrow A$ is the multiplication by $n$ map.
	Then there exists a system of canonical isomorphisms
	$$
		n^* \tau_P^\mathcal{M} \colon  \; n^* \tau_{P_1}^* \mathscr{L}\cong n^*\mathscr{L}
	$$ 
	defined over $F(A[4N^2])$  for any $n$ such that $2N | n$.
\end{proposition}

\begin{proof}
	Consider $\mathscr{L}$ to be an invertible sheaf.
	Let $Q=(Q_n)$ be the unique element in $V(A)$ 
	such that $2Q=P$, that is,  $Q_n=P_{2n}$.  
	We put $\mathscr{L}_{Q_n}=\tau_{Q_n}^*(n^*\mathscr{L})\otimes (n^*\mathscr{L})^{-1}$. 
	Then if $2N \vert n$, we have
	$$
		\mathscr{L}_{Q_n}=\tau_{Q_n}^*(n^*\mathscr{L})\otimes (n^*\mathscr{L})^{-1}
		=n^*(\tau_{Q_1}^*\mathscr{L}\otimes \mathscr{L}^{-1})
		\cong \tau_{nQ_1}^*\mathscr{L}\otimes \mathscr{L}^{-1} =  \mathcal{O}_A
	$$ 
	where the isomorphism is given by the theorem of the square.
	We choose  isomorphisms $ \rho_{Q_n} : \mathscr{L}_{Q_n} \cong   \mathbb{A}^1_A$ 
	and $\rho_{-1}:  (-1)^*\mathscr{L}\cong  \mathscr{L}$. 
	We consider the isomorphism  
	$$
		(-1)^*\mathscr{L}_{Q_n}=\tau_{-Q_n}^*(-n)^*\mathscr{L}\otimes (-n)^*\mathscr{L}^{-1} 
		\cong \tau_{-Q_n}^*(n^*\mathscr{L})\otimes (n^*\mathscr{L})^{-1}
	$$ 
	where the last isomorphism is given by $\rho_{-1} \otimes \rho_{-1}^{\otimes-1}$. Since
	$\rho_{-1}$ is unique up to a constant multiple, this 
	isomorphism does not depend on the choice of $\rho_{-1}$. 
	Then using this isomorphism, we have 
	\begin{align*}
		n^*(\tau_{P_1}^* \mathscr{L})&\otimes(n^*\mathscr{L})^{-1} \\
		&=\tau_{P_n}^*(n^*\mathscr{L})\otimes(n^*\mathscr{L})^{-1} 
		=\tau_{Q_n}^*\left(\mathscr{L}_{Q_n} \otimes (-1)^*\mathscr{L}_{Q_n}^{-1} \right)
		\cong \mathcal{O}_A
	\end{align*}
	where the last isomorphism is given by $\rho_{Q_n} \otimes \rho_{Q_n}^{-1}$ and is 
	independent of the choice of $\rho_{Q_n}$. 
	Since $Q_{2N}$ is defined over $F(A[4N^2])$, we have a canonical 
	isomorphism 
	$n^*\tau_{P}^\mathcal{M}: n^*(\tau_{P_1}^*\mathscr{L}) \cong n^*\mathscr{L}$
	defined over $F(A[4N^2])$ as desired.
\end{proof}
 
The above proposition may be regarded in geometric terms as follows.
We consider the line bundle $\mathscr{L}=\mathscr{L}(D)$ associated to 
an invertible sheaf $\cO_A(D)$ corresponding to the Cartier divisor $D=\{(U_i, f_i)\}_i$. 
Then by definition, the geometric line bundle $\mathscr{L}(D):=\mathbb{V}(\cO_A(-D))$ is 
given as $\mathcal{S}\!\operatorname{pec}\left(\mathrm{Sym} (\cO_A(-D))\right)$, namely, 
it is the scheme obtained by patching 
$\mathbb{A}^1_{U_i}=\Spec\, \cO_{U_i} [x_i] $ on  $U_{ij}=U_i \cap U_i$ by 
$$
	\phi_{ij}\colon \mathrm{Spec}\, \cO_{U_{ij}} [x_i]  \rightarrow 
	\mathrm{Spec}\, \cO_{U_{ji}} [x_j] , \quad x_j \mapsto  {f_j}{ f_i}^{-1} x_i. 
$$ 
Using this notation, the above proposition implies the following.
\begin{corollary} 
 For $P=(P_n)_{n \in \mathbb{N} } \in V(A)$ and $NP_1=0$, 
	there exists canonical rational functions 
	$(f_{D, P_n}^{\mathcal{M}})_{n \in 2N \mathbb{N}}$ on $A$ 
	over $F(A[2N^2])$   
	without any ambiguity of constant multiple such that  
	$$
		\mathrm{div} (f_{D, P_n}^{\mathcal{M}})=n^*\tau^*_{P_1}(D)-n^*D
	$$
  	and  the isomorphism   
	$n^*\tau_{P}^\mathcal{M}:  n^*\tau^*_{P_1}\mathscr{L}(D) \cong n^*\mathscr{L}(D) $ is given by patching 
	$$
 		\mathrm{Spec}\, \cO_{ n^{-1}(U_{i})} [y_i] \rightarrow
		\mathrm{Spec}\, \cO_{n^{-1}(U_{i})} [x_i]   , \quad  x_i \mapsto  
		(f^{\mathcal{M}}_{D, P_n}) \cdot 
		(n^*f_i)\cdot(n^*\tau^*_{P_1} f_i)^{-1}   y_i. 
	$$
\end{corollary}
The above construction also works scheme-theoretically 
for abelian schemes
(see Mumford \cite{Mum5}, \S 5, Appendix I).
Mumford uses the isomorphisms $\tau_{P}^{\cM}$ to construct his algebraic theta functions.
Theta function is a section of  a line bundle if we fixed a trivialization of the line bundle 
on the universal cover, which may be  algebraically regarded  as $V(A)$. 
Hence it would be natural to consider not only a single trivialization of a line bundle 
but also a system of trivializations. 
Let $(X_n)_{n\in \mathbb{N}}$ be a projective system of schemes over $A$ 
such that the diagram
\begin{equation*}
	\begin{CD}
		X_{mn} @>>> X_n \\
		@V\pi_{mn}VV @VV\pi_nV \\
		A @>m>>A 
	\end{CD}
\end{equation*}
is commutative for all positive integer $m$ and $n$. 

\begin{definition}
	We say that $((X_n, \pi_n, \varphi_n))_{n \geq 1}$ is a 
	\textit{system of trivializations of $\mathscr{L}$}, if 
	$$
		\varphi_n: \; \pi_n^*n^*\mathscr{L} \cong \mathbb{A}_{X_n}^1
	$$
	 is an isomorphism compatible with the natural projections.  
\end{definition}
For example, if we are over $\bbC$, then we may take $X_n $ for any $n$ as the complex universal cover of $A$ and 
consider the natural system of complex analytic trivializations.  As another example, we may take 
a canonical system of  trivialization coming from a point 
$(P_n)_n \in V(A)$ defined as follows. 
Let $X_n$ be $\operatorname{Spec} {\ol F}$ and consider the morphism  
$\pi_n : X_n \rightarrow A$ coming from the inclusion $P_n \hookrightarrow A$.  
If we  fix an isomorphism $[0]^*\mathscr{L} \cong \mathbb{A}_{{F}}^1$, 
we have a trivialization 
$n^*\mathscr{L}\times_A A[n]= [0]^*\mathscr{L} \times_{\ol{F}} A[n] \cong \mathbb{A}_{A[n]}^1$. 
From this we have a system of trivializations 
$ \varphi_n: \pi_n^*(n^*\mathscr{L}) \cong \mathbb{A}_{X_n}^1$.

For a system of trivializations $((X_n, \pi_n, \varphi_n))_{n \in \mathbb{N}}$, 
the translation isomorphism $\tau_{P}^{\cM}$ gives a system of trivializations
$((X_n, \pi_n, \varphi_{P,n}))_{n \in 2 N \mathbb{N}}$
of $\tau^*_{P_1}\mathscr{L}$, where  $\varphi_{P,n}$
is the composition
\begin{equation*}
\begin{CD}
  \varphi_{P,n}\;:\; \pi^*_n n^* \tau^*_{P_1} \mathscr{L} @>\pi^*_n n^*\tau_{P}^\mathcal{M}>>  \pi^*_n n^*\mathscr{L}
@>\varphi_n>> \mathbb{A}_{X_n}^1.
\end{CD}
\end{equation*}

Let $s=s_D$ be a meromorphic (rational) section of $\mathscr{L}(D)$ defined by 
$s_D \vert_{U_i} = f_i \vert_{U_i}$ for a Cartier divisor $D=((U_i, f_i))_i$. 
We define $n^*\vartheta_{s}$ to be the rational morphism on $X_n$ given by 
\begin{equation*}
\begin{CD}
	n^*\vartheta_{s} \;:\; X_n @>\pi^*_n n^*s>>   \pi^*_n n^*\mathscr{L} 
	@>\varphi_n>> \mathbb{A}_{X_n}^1.
\end{CD}
\end{equation*}
\begin{definition}
	We define the rational morphism 
	$ n^*U_{P}^{\mathcal{M}}\vartheta_{s} $ to be the composition
	\begin{equation*}
		\begin{CD}
 			n^*U_{P}^{\mathcal{M}}\vartheta_{s} \colon \; 
			X_n @>\pi^*_n n^*\tau^*_{P_1}s>>   \pi^*_n n^* \tau^*_{P_1}\mathscr{L}  
			@>\varphi_{P,n}>> \mathbb{A}_{X_n}^1.
		\end{CD}
	\end{equation*}
	Then $n^*U_{P}^{\mathcal{M}}\vartheta_{s}$
	is  explicitly given as $n^*\vartheta_{s}\cdot  \pi^*_n f^{\mathcal{M}}_{D, P_n}$. 
	We will see in Proposition \ref{proposition; explicit U} below 
	that notations are compatible with that in 
	\S \ref{subsection:2-3}. 
\end{definition}

We remark that  for a holomorphic section $s$, 
taking the canonical trivialization coming from a point $(P_n)_n \in V(A)$ 
for a fixed  isomorphism $[0]^*\mathscr{L} \cong \mathbb{A}_{\ol{F}}^1$, 
we have a morphism 
$$
	\vartheta_s \colon V(A) \rightarrow \ol{F}, \quad  
	(P_n)_n \mapsto \varinjlim_{n \in 2N\mathbb{N}}\; n^*U^{\mathcal{M}}_{P}\vartheta_{s}. 
$$
We note that since $2N | n$, the value $n^*U^{\mathcal{M}}_{P}\vartheta_{s}$ is constant.
The above function $\vartheta_s$ on $V(A)$ is Mumford's adelic theta function. 

We now relate the function $n^* U_{P}^\mathcal{M} \vartheta_s$ to the translation 
of theta reduced functions defined in \S \ref{subsection:2-3}.  Assume that $F$ is a subfield of $\C$. 
We first give an explicit description of the function $f_{D,P_n}^{\mathcal{M}}$ 
in terms of reduced theta functions. The idea is, since 
this function is independent of the choices of $\rho_{-1}$ and $\rho_{Q_n}$, we 
pick convenient  $\rho_{-1}$ and $\rho_{Q_n}$ over $\C$ using reduced theta functions. 

Fix a complex uniformization  $\pi: V/\Lambda \cong A(\bbC)$ and consider the morphism
\begin{equation}
	\iota: \Lambda \otimes \Q \rightarrow (\Lambda \otimes \Q) /\Lambda 
	\hookrightarrow  A(\C)_{\mathrm{tor}}=A(\overline{F})_{\mathrm{tor}}. 
\end{equation}
Then we obtain a canonical injection $\widetilde{\iota}: \Lambda \otimes \Q \hookrightarrow V(A)$ by 
mapping $v \in \Lambda \otimes \Q$ to 
$P=(P_n)_n$,  where $P_n$ is the image by $\iota$ of $v/n$. 
Let $v_0$ be an element $\Lambda \otimes \Q$ and let 
$w_0$ be $v_0/2$. 
We denote the image by $\widetilde{\iota}$  of $v_0$ (resp. $w_0$) as $P=(P_n)$ (resp. $Q=(Q_n)$). 
Since $\mathscr{L}=\mathscr{L}(D)$ is symmetric, 
the function $\rho_{-1}(v):=\vartheta_D(-v)/\vartheta_D(v)$ is periodic with 
divisor $(-1)^*D-D$. We take $\rho_{-1}: (-1)^*\mathscr{L} \cong \mathscr{L}$ as an isomorphism 
defined by  the function $\rho_{-1}(v)$.
Let $\vartheta_D(v)$ be a reduced theta function associated to $D$. 
By the transformation formula of  reduced theta functions, if 
$nw_0 \in  \Lambda$, we see that 
the function 
	$$
		\rho_{Q_n} (v) := \Exp{ - \pi H(n v, w_0)} \vartheta_D(nv+w_0) \vartheta_D(nv)^{-1}
	$$
	is meromorphic in $v$ and is periodic with respect to $\Lambda$ with 
	divisor  $\tau^*_{Q_n}\left(n^*D\right)-\left(n^*D\right)$. 
We take $\rho_{Q_n}: \mathscr{L}_{Q_n} \cong \mathbb{A}^1_A$ as an isomorphism 
defined by  $\rho_{Q_n}(v)$.

\begin{lemma}\label{lemma; explicit U}
	The rational function  $\pi^* f^{\mathcal{M}}_{D, P_n}$ is given by 
	$$
		f^{\mathcal{M}}_{D, P_n}(v)=\frac{n^*U_{v_0}^{\mathcal{M}}\vartheta_D(v) }{n^*\vartheta_D(v)}
	$$
	where $n^*U_{v_0}^{\mathcal{M}}$ is the translation of reduced theta functions 
	defined in \S \ref{subsection:2-3}. 
\end{lemma}

\begin{proof}
	We have
	\begin{align*}\label{equation: rho}
		&f^{\mathcal{M}}_{D, P_n}(v)= \tau^*_{Q_n} \left( \rho_{Q_n} (v) \cdot \rho_{Q_n} (-v)^{-1} 
		\cdot  \tau^*_{-Q_n}\rho_{-1}(nv) \cdot  \rho_{-1}(nv)^{-1}\right)
		  \\
		&= \Exp{- \pi H(n v, v_0)- \frac{\pi}{2} H(v_0, v_0)} 	
		\frac{\vartheta_D( nv+v_0 ) \vartheta_D(-nv-w_0) 
		\rho_{-1}(nv)}{\vartheta_D(nv+w_0) \vartheta_D(-nv)\tau^*_{Q_n}(\rho_{-1}(nv))} \\
		&= \Exp{- \pi H(n v, v_0)- \frac{\pi}{2} H( v_0, v_0)} 
		\frac{\vartheta_D(nv+v_0)}{\vartheta_D(nv)}
		=\frac{n^*U_{v_0}^{\mathcal{M}}\vartheta_D(v) }{n^*\vartheta_D(v)}. \nonumber
	\end{align*}
	as desired.
\end{proof}

Now we consider a system of complex analytic trivialization of $\mathscr{L}(D)$ given  by  
the reduced theta function  $\vartheta_D(v)$. Namely, the trivialization 
$((X_n, \pi_n, \varphi_n))_{n \geq 1}$ is  such that  
$X_n=V$, $\pi_n$ is a projection induced by $\pi$ and $\varphi_n=n^*\varphi_{\vartheta_D}$ where 
\begin{equation}\label{equation; trivialization by theta_D}
\varphi_{\vartheta_D}:  \pi^*\mathscr{L} \cong \pi^*\mathscr{L}(H, \alpha)=\mathbb{A}^1_V
\end{equation} is 
given by  
$$
	\mathrm{Spec}\, \cO_{U_{i}} [x_i]   \;\rightarrow  \;
	\mathbb{A}^1_{\pi^{-1}(U_i)}=\mathrm{Spec}\, \cO_{\pi^{-1}(U_i)} [X]
 	, \quad X \; \mapsto\; { \vartheta_D(v)} {f_i(v)}^{-1} x_i
$$
where $D=\{(U_i, f_i)\}_i$.    Then we have the following.

\begin{proposition}\label{proposition; explicit U}
Let $D$ be a Cartier divisor  $\{(U_i, f_i)\}_i$ and $s_D$ 
a meromorphic section defined by $s_D \vert_{U_i} = f_i \vert_{U_i}$. 
Let  $\vartheta_D(v)$ a reduced theta function associated to $D$. 
Let $((X_n, \pi_n, \varphi_{n}))_{n \in \mathbb{N}}$ be 
the system of complex analytic trivializations of $\mathscr{L}(D)$ given  by  
the reduced theta function  $\vartheta_D(v)$ as above. 
Then the rational morphism $n^*\vartheta_{s}: V \rightarrow 
	\mathbb{A}^1_V $ corresponding to the section $n^*s_D$ 
	by this trivialization is 
	the reduced theta function $n^*\vartheta_D(v)$, and 
the rational morphism  $n^*U_{v_0}^{\mathcal{M}} \vartheta_{s}: V \rightarrow 
	\mathbb{A}^1_V $ corresponding to the section $n^*\tau_{P_1}^*s_D$ 
	by the system of trivializations $((X_n, \pi_n, \varphi_{P,n}))_{n \in 2N \mathbb{N}}$   is given by 
	$$
		n^*U_{v_0}^{\mathcal{M}}\vartheta_{s}\;:\; V  \;\longrightarrow \; \mathbb{A}^1_V, \qquad 
		v \;\mapsto\; n^*U_{v_0}^{\mathcal{M}}\vartheta_D(v).
	$$ 
In other words, the trivialization $\varphi_{P,n}$ is the pull-back by the multiplication $n$ of 
\begin{equation*}
\varphi_{U_{v_0}^{\mathcal{M}}\vartheta_D}:\;  \pi^*\tau^*_{v_0}\mathscr{L} \cong \pi^*\mathscr{L}(H, \alpha \cdot \nu_{v_0})=\mathbb{A}^1_V
\end{equation*}
for $\tau^*_{v_0}D=\{(\tau^*_{v_0}U_i, \tau^*_{v_0} f_i)\}_i$ defined as in (\ref{equation; trivialization by theta_D}).  
\end{proposition}

\begin{proof}
The first part directly follows the definition of  $((X_n, \pi_n, \varphi_{n}))_{n \in \mathbb{N}}$. 
The last part follows from the first part of this proposition and Lemma \ref{lemma; explicit U} since 
$n^*U_{P}^{\mathcal{M}}\vartheta_{s}=n^*\vartheta_{s}\cdot  \pi^*_n f^{\mathcal{M}}_{D, P_n}$. 
	\end{proof}

Since $U_{v_0}^{\mathcal{M}}\vartheta_D(v)$ is a reduced theta function for 
$\mathscr{L}(H, \alpha \cdot \nu_{v_0})$, the trivialization $\varphi_{P,n}$ gives 
an isomorphism 
$n^*\tau_{v_0}^*\mathscr{L} \cong n^*\mathscr{L}(H, \alpha \cdot \nu_{v_0}).$
Hence Mumford's theory gives an algebraic way to 
choose an isomorphism $\tau_{v_0}^*\mathscr{L} \cong \mathscr{L}(H, \alpha \cdot \nu_{v_0})$ 
from $\mathscr{L} \cong \mathscr{L}(H, \alpha)$ (up to a $n$-th root of unity).


Now we show that the translation $U_{v_0}$ preserves the algebraicity of the coefficients of the Laurent 
expansion of reduced theta functions.

\begin{theorem}\label{theorem: algebraicity}
	Let $A$ be an abelian variety over a number field $F$ of dimension $g$ with or 
	without complex multiplication, and let $\omega_1, \dots, \omega_g$ be invariant differential $1$-forms on $A$ 
	which form a global basis of the K\"ahler differential $\Omega_{A/F}^1$. 
	We let $\pi$ be the complex uniformization $\C^g/\Lambda \rightarrow 
	A(\C)$ such that the pull back of $\omega_i$ is the 
	differential form $dz_i$ where $z_i$ is the canonical $i$-th coordinate of $\C^g$.  
	Let $\mathscr{L}$ be a symmetric line bundle on $A$ and 
	$s$ a meromorphic section of $\mathscr{L}$ with divisor $D$ defined over $F$.  
	Let $\vartheta_D(v)$ be a reduced theta function associated to $s$. 
	We denote by $v_0$ an element of $\Lambda \otimes \mathbb{Q}$ such that $N v_0 \in \Lambda$
	for an integer $N > 0$.
	Let $f$ be a rational function over $F$ that defines the Cartier divisor $D$  in a neighborhood of the origin, and 
	let $g_{v_0}$ be  a rational function over $F(A[N])$ that defines 
	the Cartier divisor $\tau_{v_0}^*D$  in a neighborhood of the origin. 
	We  assume that the coefficients of the Taylor expansion of 
	$(\vartheta_D/f)(v)$ at the origin 
	are  contained in $F$. Then 
	the coefficients of the Taylor expansion of 
	$g_{v_0}(v)^{-1}U_{v_0}^{\mathcal{M}} \vartheta_D(v)$ 
	and  $g_{v_0}(v)^{-1}U_{v_0} \vartheta_D(v)$ are contained in $F(A[4N^2])$. 
\end{theorem} 

\begin{proof}
Let $n$ be an integer such that $2N | n$.
Since $U_{v_0} \vartheta(v)$ differs from $U^{\mathcal{M}}_{v_0} \vartheta(v)$ only by a $n$-th root of unity, 
it suffices to show the theorem for $U^{\mathcal{M}}_{v_0} \vartheta_D(v)$. 
Then by the construction, the rational function $f^{\mathcal{M}}_{D, P_n}$ is defined over $F(A([4N^2])$. 
Since the derivation $\partial/\partial z_i$ with respect to $\omega_i$ is defined over $F$, 
the Taylor coefficients of the holomorphic function 
$g_{v_0}(v)^{-1}  f^{\mathcal{M}}_{D, P_n}(v) f(v)$
at origin are in $F(A([4N^2]) $. 
The theorem follows from this since 
$n^*U^{\mathcal{M}}_{v_0} \vartheta_D(v)=n^*\vartheta_D(v)\cdot f^{\mathcal{M}}_{D, P_n}(v)$. 
\end{proof}

\begin{corollary}\label{corollary: AT2} 
Assume that the complex torus $\C/\Gamma$ has complex multiplication 
in  $\cO_K$ 
and 
has an Weierstrass model $E: y^2 = 4 x^3 - g_2x - g_3$ 
such that $g_2, g_3 \in F$. 
	For $z_0$, $w_0 \in \Gamma \otimes \frac{1}{n} \Z$, the Laurent expansion of 
	$\Theta_{z_0, w_0}(z,w; \Gamma)$ at $z=w=0$ has coefficients in 
	$F(E([4n^2]))$. 
	In particular, the Eisenstein-Kronecker numbers
	$$
		e^*_{a,b}(z_0, w_0; \Gamma)/A^a
	$$
	are algebraic. 
	\end{corollary}
\begin{proof}
We apply the theorem to $A=E\times E$ and divisor $\Delta-(E\times \{0\})- (\{0\}\times E)$. Then 
this follows from Theorem \ref{theorem: generating function theorem}, 
Corollary \ref{corollary: AT1} and 
 Theorem \ref{theorem: algebraicity}. 
	\end{proof}

\begin{corollary}[Damerell's Theorem]\label{corollary: Damerell}
Let $\varphi$ be a Hecke character of an imaginary  quadratic  field $K$ with conductor $\ff$ 
and infinity type $(a,-b)$, where $a, b$ are distinct non-negative integers. 
Let $\Omega$ be a complex number such that 
$\ff\Omega$ is the period lattice of a pair $(E, \omega_E)$ consisting of 
an elliptic curve $E$ and an invariant differential $\omega_E$ defined over an algebraic number field. Then the numbers
	$$
		\left(\frac{2\pi}{\sqrt{d_K}}\right)^a \frac{L_{\ff}(\varphi, 0)}{\Omega^{a+b}}
	$$
are algebraic, where $-d_K$ is the discriminant of $K$.
\end{corollary}

\begin{proof}
Suppose that $b>a \geq 0$. 
By proposition \ref{proposition: EK and L} ii), we have 
$$L_{\ff}(\varphi, 0)=
  \frac{1}{w_\ff}
		\sum_{\fa \in I_K(\ff)/P_K(\ff)} \varphi(\fa)
	e^*_{a, b}(\alpha_\fa,  0; \fa^{-1}\ff ). 
	$$	
	It is known (for example, see \S \ref{subsection:5-2} below) that 
the lattice $\fa^{-1}\ff \Omega$ also comes from
 a Weierstrass model defined over an algebraic number field. 
Therefore Corollary \ref{corollary: AT2}
shows that the number 
$$
A(\fa^{-1}\ff\Omega)^{-a}e^*_{a, b}(\alpha\Omega,  0; \fa^{-1}\ff\Omega )
=A(\fa^{-1}\ff\Omega)^{-a}\Omega^{-b}\ol{\Omega}^a 
e^*_{a, b}(\alpha\Omega,  0; \fa^{-1}\ff)
$$ 
is algebraic. 
Since 
$$A(\fa^{-1}\ff\Omega)=N(\fa^{-1}\ff)\Omega\ol\Omega A(\cO_K)
=N(\fa^{-1}\ff)\Omega\ol\Omega \sqrt{d_K}/2\pi, 
$$
the number 
$$
		\left(\frac{2\pi}{\sqrt{d_K}}\right)^a \frac{e^*_{a,b}(\alpha_\fa\Omega, 0; 
		\fa^{-1}\ff\Omega)}{\Omega^{a+b}}
	$$
	is  algebraic.
\end{proof}

We digress here to give the algebraic property of the special values of 
the Weierstrass $\sigma$-function, which may be of independent interest.

\begin{theorem}
	Let $E$ be an elliptic curve defined over an algebraic number field $F$ 
	(with or without complex multiplication). 
	Let $\Gamma$ be the period lattice of $E$ for an invariant differential defined over $F$. 
	Let $\sigma(z)$ be the Weierstrass $\sigma$-function for the lattice $\Gamma$. 
	Let $z_n$ be an element of $\C$ such that $nz_n \in \Gamma$. 
	Then the Taylor coefficients of   
	$$
		\exp\left[-\eta(z_n)\left(z+\frac{z_n}{2}\right)\right]\sigma(z+z_n)
	$$
	at the origin is in  $F(E[4N^2])$.   In particular,  we have $\exp[-\eta(z_n)z_n/2]\sigma(z_n) \in F(E[4N^2])$. 
	This last value is the value of Mumford's adelic theta function 
	for the section corresponding to $1 \in \Gamma(E, \mathscr{L}([0]))$ at 
	$z_n \in \Lambda \otimes \Q \hookrightarrow V(E)$. 
\end{theorem}

\begin{proof}
	We let $\mathscr{L}([0]) \cong \mathscr{L}(H, \alpha)$ and 
	let  $\theta(z)$ be the theta function of Example \ref{robert}. 
	Then 
	\begin{multline*}
		\exp\left[-H(z, z_n)-\frac{\pi}{2} H(z_n,z_n)\right] \theta(z+z_n)/\theta(z)
		\\=\exp\left[-\eta(z_n)\left(z+\frac{z_n}{2}\right)\right]\sigma(z+z_n)/\sigma(z). 
	\end{multline*}
	Hence its pull back by $n: E \rightarrow E$ is  $f^{\mathcal{M}}_{[0], P_n}$. 
	Therefore if we use the trivialization of $\cO_E([0])$ defined by $\sigma(z)$, then
	$f_{[0], P_n} ^{\mathcal{M}}(z)\cdot n^*\sigma(z)$ is the translation 
	of the section $\sigma(z)$ by Mumford's isomorphism. 
	Since the Taylor coefficients of  $\sigma(z)$ at the origin are in $F$, we have the assertion.  
\end{proof}

\subsection{The $\fp$-adic integrality  of reduced theta functions 
} \label{subsection:3-3}
%

We use the same notation as in \S \ref{subsection:3-1}. In addition, 
in what follows, let $\bbC_p$ be the completion of the algebraic closure $\ol\bbQ_p$ of $\bbQ_p$.
We denote by $\cO_{\bbC_p}$ and $\fm_{\bbC_p}$ the ring of integers and the maximal ideal of $\bbC_p$.  
We fix once and for all an embedding
$
	i_p : \ol K \hookrightarrow \bbC_p. 
$
Let $W= W(\ol\bbF_p)$ be the ring of Witt vectors with
coefficients in $\ol\bbF_p$. 

The purpose of this section is to prove the $\fp$-integrality property of the Taylor expansion of  
the reduced theta function $\vartheta_s(z)$ with respect to a formal group parameter of $A$.
For this result, it is necessary to assume some ordinarity condition on $p$. 
Let $\fp$ be a prime ideal of the CM field $K$ over $p$.  
We assume that 
\begin{equation}\label{ordinarity}
	\prod_{i=1}^g \phi_i(\fp) \quad \text{ is prime to}\quad \prod_{i=1}^g \ol \phi_i( \fp).
\end{equation}

Let $\fP$ be a prime of $F$ over $\fp$ such that the completion $F_{\fP}$ of $F$ at $\fP$ is the completion of $F$ in 
$\bbC_p$ with respect to the inclusion $i_p$, and let 
$R$ be the ring of integers $\cO_{F_{\fP}}$ of $F_{\fP}$. 
We assume that the CM abelian variety $A/F$ has a proper 
 model $\mathcal{A}$ over $\cO_F$ which is smooth over $\fP$ and 
the invariant differentials $\omega_1, \dots, \omega_g$ give 
a global basis of $\Omega_{\mathcal{A}/R}^1$. 
We sometimes denote the pull-back  $\mathcal{A} \otimes_{\cO_F} R$ of $\mathcal{A}$ 
 on $R$ also by $\mathcal{A}$. 
 Let $t_1, \dots, t_g$ be a local parameter of $\mathcal{A}/R$ at the origin and 
we consider the formal completion $\wh{\mathcal{A}}/R$
of $\mathcal{A}/R$ at the origin 
with this parameter. 
We put $\lambda(t):=(\lambda_1(t), \dots, \lambda_g(t))$ where 
$\lambda_i(t):=\lambda_i(t_1, \dots, t_g)$ is the formal logarithm 
of $\wh{\mathcal{A}}/R$ corresponding to 
the differential form $\omega_i$. 
(Namely, $\lambda_i(t)$ is the power series of $t$ such that 
$d\lambda_i(t)=\wh\omega_i(t)$ and $\lambda_i(0)=0$. This power series exists 
since $d\omega_i=0$).
Let $\vartheta_s(v)$ be the reduced theta function associated to 
a section of $\Phi$-admissible symmetric line bundle $\mathscr{L}$ 
as in Proposition \ref{proposition; algebraicity at the origin}. 
Since the coefficients of the Taylor expansion of $\vartheta_s(v)$ are 
algebraic, we can consider the formal composition of $\vartheta_s(z)$ with 
$v=\lambda(t)$. 
We denote it by $\wh \vartheta_s(t)$, which is a priori in $F_{\fP}[[t_1, \dots, t_g]]$. 
We show that it is actually in  $R[[t_1, \dots, t_g]]$.

\begin{lemma}\label{useful}
Let $D$ be an arithmetic divisor of $\mathcal{A}$ defined over $R$.
Suppose that $D$ is principal and the divisor of the poles of $D$ 
does not intersect $(\fP, t_1, \dots, t_g)$.  
Then the formal completion along  $(\fP, t_1, \dots, t_g)$ of a rational function $f$ with divisor 
$D$ is an element in $R[[t_1,\dots, t_g]]$. 
\end{lemma}
\begin{proof}
By assumption,   for a sufficiently small open neighborhood $U$ of the 
ideal  $(\fP, t_1, \dots, t_g)$, 
the rational function $f$ defines a morphism $U \rightarrow \mathbb{A}^1_R$.  
Since the completion of $U$ along $(\fP, t_1, \dots, t_g)$ is 
$R[[t_1, \dots, t_g]]$,  the completion of $f$ is an element in 
$R[[t_1, \dots, t_g]]$. 
\end{proof}

\begin{proposition}\label{proposition: integrality of theta} 
Let $\fp$ be a prime ideal of $F$ that satisfies the ordinarity condition 
(\ref{ordinarity}). 
Let  $\mathscr{L}=\mathscr{L}(H, \alpha)$ be a $\Phi$-admissible line bundle on 
$\mathcal{A}$, and assume that $\alpha^N=1$ for some integer $N$ prime to $p$. 
Let $s$ be a meromorphic section of $\mathscr{L}$ 
with a (arithmetic) Cartier  divisor $D$. 
Let $f$ be a rational function which defines 
the Cartier divisor $D$ on $\mathcal{A}$ in a neighborhood of the origin. 
If $(\vartheta_s/f)(0) \in R^\times$, 
then we have
		$$
			\wh f(t)^{-1} \cdot \wh\vartheta_s(t) \in R[[t_1, \dots, t_g]]^\times.
		$$
\end{proposition}

\begin{proof}
By considering $f(v)^{-1} \cdot \vartheta(v)$ instead of $\vartheta(v)$, 
we may assume that $f=1$ and the  divisor  $D$ 
does not  intersect the point $(\fP, t_1, \dots, t_g)$. 
The proof is similar to that of Proposition \ref{proposition; algebraicity at the origin}. 
In this case, by the ordinarity condition, we take $a \in N\cO_K$ such that 
$\phi_i(a) \in \fp$ but  $\ol\phi_i(a) \notin \fp$ for our fixed $\fp$. 
With the same notation as in the proof of Proposition \ref{proposition; algebraicity at the origin}, 
the function $F_a(z)$ has no zeros and poles at the origin except those coming from the
arithmetic divisors of $\mathrm{Spec} R$. Therefore 
by Lemma \ref{useful} 
the power series $\wh F_a(t)$ and the power series $\wh F_a(t)^{-1}$
are elements of $R[[t_1, \dots, t_g]]\otimes \Q$. 
Since $\wh F_a(0)=1$, we have $\wh F_a(t)$ is in $R[[t_1, \dots, t_g]]^\times $. 
If necessary, by changing the local parameter $t_1, \dots, t_g$, 
we may assume that $z_i =\lambda_i(t) \equiv t_i \mod (t_1, \dots, t_g)^2$. 
Then since $\iota(a)^*\omega_i=\phi_i(a)\omega_i$, 
we have $[a]t_i \equiv \phi_i(a) t_i \mod (t_1, \dots, t_g)^2$. 
Then as in the proof of  Proposition $\ref{proposition; algebraicity at the origin}$, 
if we write $\wh\vartheta_s(t)=\sum c_n t^n$ using the multi-index $n$, 
we have 
	\begin{equation*}
		\left(\prod_{i=1}^g \phi_i(a)^{n_i}\right) c_{n} =
		\left( \prod_{i=1}^g \ol\phi_i(a)^{n_i}\right) c_{n} + u(\{c_k\}),
	\end{equation*}
	where $u(\{z_k\})$ is a polynomial over $R$ of several variables $\{z_k\}$ 
          where $k$ runs through 
	multi-indices whose total degree is less than $n$. 
	Then by the assumption on $a$, the element 
	$\prod_{i=1}^g \phi_i(a)^{n_i}$ is not a $\fp$-adic unit but 
	$ \prod_{i=1}^g \ol\phi_i(a)^{n_i}$ is a $\fp$-adic unit if $n \not=0$. 
Since $\wh\vartheta_s(0) \in R^\times$, we have $c_0 \in R^\times$. Hence by induction,
 we have $c_n \in R$. 
\end{proof}

\begin{proposition}\label{proposition: integrality of translation} 
We use the same notations and the assumption in Proposition \ref{proposition: integrality of theta}. 
Here we also assume that $\mathscr{L}$ is symmetric. 
Let $v_0$ be an element of $V$ such that $N v_0 \in \Lambda$ for $N$ prime to $p$.  
Let $g_{v_0}(v)$ be a rational function which defines the 
Cartier divisor of $\tau_{v_0}^*D$ in a neighborhood of the origin. 
Then we have 
$$
	\wh g_{v_0}(t)^{-1} \cdot \wh \vartheta_{v_0}^{\cM}(t) \in 
	\cO_{F_\fP}(A[4N^2])[[t_1, \dots, t_g]]^\times. 
$$
\end{proposition}
\begin{proof}
Let $R'$ be the ring $\cO_{F_\fP}(A[4N^2])$, and for $n$ such that $2N|n$, let 
$P_n$ be an element corresponding to the point $v_0/n$ through
the uniformization $\pi: V/\Lambda \rightarrow A(\C)$. 
Then the rational function
	$$
	n^*g_{v_0}^{-1}
		\cdot f_{D,P_n}^{\mathcal{M}} \cdot n^*f
		$$
	on $\mathcal{A}/R$ does not have a zero nor a pole  around the origin.  
	Therefore by Lemma \ref{useful}, its
	Taylor expansion with respect to the parameter $t$ is in $R'[[t]]^\times$.  
	Hence by Proposition \ref{proposition: integrality of theta}, the power series 
	\begin{align*}
	& n^*\wh g_{v_0}(t)^{-1}   \cdot n^*\wh\vartheta_{v_0}^\cM(t) \\
	& =n^*\wh g_{v_0}(t)^{-1}
		\cdot \wh f_{D,P_n}^{\mathcal{M}}(t) \cdot n^*\wh f (t)
	 \cdot n^*\left(\wh f(t)^{-1}\cdot\wh\vartheta(t)\right)
	\end{align*}
	is in $R'[[t]]^\times$. 
	Since the multiplication by $n$  is \'etale over $R'$, it induces an 
	isomorphism on $R'[[t]]$. Hence we have the assertion.
\end{proof}

Now we assume 
that $K$ is an imaginary  quadratic  field. 
 We assume that $p \geq 5$ and 
we fix a Weierstrass model 
$$
	\cE : y^2 = 4 x^3 - g_2 x - g_3
$$
of $E$ over $\cO_F$ which is good for a prime $\fP$ over $\fp$.
We let $\wh\cE$ be the formal group associated to $\cE$ with respect to the parameter $t = -2x/y$. 
As before, we denote by $\lambda(t)$ the formal logarithm of $\wh\cE$, 
which is a power series giving a homomorphism of formal groups
$\wh\cE \rightarrow \wh\bbG_a,  z = \lambda(t)$, and normalized so that $\lambda'(0)=1$.

\begin{corollary}\label{proposition: integrality of robert theta} 
	Let $p \geq 5$ be a prime that splits as $(p) = \fp \ol\fp$ in $\cO_K$.
	Let $z_0$, $w_0 \in n^{-1}\Gamma$ for an integer $n$ prime to $p$. 
	We put  $\delta(z)=1$ if $z \in \Gamma$ and  $\delta(z)=0$ otherwise. 
	Let $\theta(z)$ be the theta function of Example \ref{robert}, and 
	let $\Theta(z, w)$ be the Kronecker theta function. 
	Then we have 
	\begin{enumerate}
 		\item  $t^{-\delta(z_0)} \wh\theta_{z_0}^{\cM}(t) \in \cO_{F_{\fP}}[[t]]^\times.$ \\
	 	\item 
		$
			s^{\delta(z_0)} t^{\delta(w_0)}   \wh\Theta_{z_0,w_0}(s,t) \in 
			\cO_{F_{\fP}}(E[4n^2])[[s,t]].
		$
	\end{enumerate} 
\end{corollary}

The integrality i) for  $z_0=0$ has already been obtained by Bernardi, Goldstein and Stephens 
(\cite{BGS} Proposition III.1) in connection with the $p$-adic height pairing.   
See also Perrin-Riou (\cite{Per} Chapitre III, \S 1.2, Lemma 2).  
Their proofs are based on the fact
the elliptic function $\theta(\alpha z)/\theta(z)^{\mathrm{deg}\,\alpha}$ has 
integral coefficients if $\alpha$ is an \'etale isogeny (at $\fP$) of odd degree. 
Mazur and Tate generalized their method to non-CM elliptic curve to 
construct their $p$-adic $\sigma$-function.  

\begin{proof}
	First we remark that the parameter $t$ defines the Cartier divisor $[0]$ in a neighborhood of 
	the origin not only for $E/F$ but also for the abelian scheme $\mathcal{E}/R$ as an arithmetic 
	divisor.  We also have $(\theta/t)(0)=1$. 
	The divisor of   $\theta_{z_0}^{\cM}$ is $\tau^*_{z_0}[0]=[-z_0]$.  
	Applying Proposition \ref{proposition: integrality of translation} by taking $g_{v_0}$ to be $t^{\delta(z_0)}$, 
	we obtain i).  For ii), since 
	$\Theta_{z_0,w_0}(z,w)$ and $\Theta_{z_0,w_0}^{\cM}(z,w)$ differ only 
	by a $n$-th root of unity, it suffices to prove the statement for $\Theta_{z_0,w_0}^{\cM}(z,w)$. 
	Since 
	$$
		\Theta_{z_0,w_0}^{\cM}(z,w)={\theta_{z_0+w_0}^{\cM}(z+w)}
		{\theta_{z_0}^{\cM}(z)^{-1}\theta_{w_0}^{\cM}(w)^{-1}},
	$$ 
	case ii) follows from i). 
\end{proof}

Together with Theorem \ref{theorem: generating function theorem}, we have 
\begin{corollary}\label{corollary: for measure} 
	$$
		\wh\Theta_{z_0, w_0}(s,t) - \pair{w_0, z_0} \delta(z_0) s^{-1} - \delta(w_0) t^{-1}
		\in \cO_{F_\fP}(E[4n^2])[[s,t]].
	$$
\end{corollary}

\subsection{The $p$-adic translation by $\fp$-power torsion points}\label{subsection:3-4}

We keep the notations of \S \ref{subsection:3-1} and \S \ref{subsection:3-3}. 
Since the power series $\wh\vartheta_{v_0}^{\cM}(t)$ has integral coefficients
(Proposition \ref{proposition: integrality of translation}), 
we can consider the composed power series
\begin{equation}\label{equation: MTP}
\wh \vartheta_{v_0}^{\cM}(t \oplus t_{\fp^n})
\end{equation} 
of  $\wh\vartheta_{v_0}^{\cM}(t)$ with the power series $t \oplus t_{\fp^n}$, where 
$\oplus$ is the formal addition and $ t_{\fp^n}$ is a $\fp^n$-torsion point of the formal group 
$\wh{\mathcal{A}}$.
In this subsection,  we explicitly determine this power series.

Let $\pi_p$ is the morphism 
\begin{equation}\label{equation: map pip}
	\pi_p : \wh{\mathcal{A}}(\fm_{\bbC_p})_{\tor} \rightarrow A(\bbC_p)_{\tor} = A(\ol\bbQ)_{\tor}
	= A(\bbC)_{\tor} \xrightarrow{\cong} (\Lambda \otimes \bbQ) / \Lambda
\end{equation}
obtained by the complex uniformization 
$\pi: \C^g/\Lambda \cong A(\C)$. 
The difficulty in calculating (\ref{equation: MTP}) stems from the fact that 
$\vartheta_{v_0}(v)$ 
is \textit{not periodic} on $V=\C^g$.  For rational functions on $A$, we have the 
following result.

\begin{lemma}\label{lemma: parallel for rational} 
Let $t_{\fp^n}$ be a $\fp^n$-power torsion point, and let 
$v_n$ be an element in $\Lambda\otimes \Q$ 
that represents the image of $t_{\fp^n}$ 
by $\pi_p$ of (\ref{equation: map pip}). 
	Let $f(v)$ be a rational function on $A/\ol F$ and 
	we put $f_{v_n}(v)=f (v+v_n)$. Then we have
	$$
		\wh f(t \oplus t_{\fp^n}) 
		 = \wh{f}_{v_n}(t). 
	$$
\end{lemma}
\begin{proof}
Let $P_n$ be a point in $A(\C_p)$ corresponding to $t_{\fp^n}$. 
Since the additive laws of $A$ and $\wh A$ are compatible, 
we have $\wh{\tau^*_{P_n}f}(t)=\wh f(t \oplus t_{\fp^n})$. 
Since $\tau^*_{P_n}f(v)=f_{v_n}(v)$, we have the desired result. 
\end{proof}

The above lemma seems to be trivial but we need some care since the equation $v = \lambda(t)$ 
has meaning only as an equality of formal power series. 
For example,  $v_n$ is not equal to $\lambda(t_{\fp^n})$ since the latter is always zero. 
The left hand side is calculated as an infinite 
sum in the $p$-adic field and the right hand side 
is as an infinite sum in the complex number field. 

For the theta function $ \vartheta_{v_0}^{\cM}(v)$, the result is as follows.

\begin{proposition}[Calculation of $p$-adic Translation]\label{proposition: parallel for theta}
Let $\fa$ be an integral ideal of $K$ prime to $p$ and let $v_0$ be an  $\fa$-torsion point of 
$\C^g/\Lambda$, namely, $v_0 \in \C^g$ and $\iota(\fa) v_0 \in \Lambda$.  
Let $t_{\fp^n}$ be a $\fp^n$-torsion point of the formal group $\wh{\mathcal{A}}$, 
and let $v_n$ be an element in $\Lambda\otimes \Q$ 
that represents the images of $t_{\fp^n}$ 
by $\pi_p$. 
	Let $\smallalpha$ be an element of  $\cO_K$ such that 
	\begin{equation}\label{equation; epsilon p}
	\smallalpha \equiv 1 \mod{\phi_i(\fp^n)}, 
	\qquad \smallalpha \equiv 0 \mod{2 \ol\phi_i(\fp^n)}.
	\end{equation}
	  Then we have
	$$
		\wh \vartheta^{\cM}_{\smallalpha v_0} (t \oplus t_{\fp^n} )
		=  \pair{ \smallalpha v_n, \smallalpha v_0/2}_\mathscr{L}\;
		\wh\vartheta^{\cM}_{\smallalpha v_0+\smallalpha v_n }( t). 
			$$
	Moreover, if  $\smallalpha \equiv 1 
	\mod{\phi_i(\mathfrak{a})}$ for all $i$, we have 
	\begin{equation}
		\wh\vartheta^{\cM}_{ v_0}(t \oplus t_{\fp^n}) = 
		\pair{\epsilon v_n, (2\epsilon-1) v_0/2}_\mathscr{L}\;
	 \wh \vartheta^{\cM}_{ v_0 + \smallalpha v_n}(t).
	\end{equation}
\end{proposition}

The proof of this proposition will be given below.
The power series $\wh \vartheta^{\cM}_{v_0+v_n }(t)$ would depends on the choice of 
$v_n$.  The exponential factors and $\smallalpha$ effectively allows us to choose a ``direction", 
giving a power series independent of the choice.
The integrality of $\wh \vartheta^{\cM}_{\smallalpha v_0}$ and this proposition says that if $\fp$-power torsion points 
$t_{\fp^m}$ and $t_{\fp^n}$ are $p$-adically close, 
 there exist  $p$-adic congruences 
between the Taylor  coefficients of  
$\wh\vartheta^{\cM}_{\smallalpha v_0+\smallalpha v_m }$ and 
$\wh\vartheta^{\cM}_{\smallalpha v_0+\smallalpha v_n }$, even though these numbers are 
a priori complex numbers. 

\begin{lemma}
	 For any $\epsilon \in \cO_K$, the function
	\begin{equation*}
		f(v) : = {\vartheta^{\cM}_{\smallalpha v_0}(\smallalpha v)}/
		{\vartheta^{\cM}_{\ol\smallalpha v_0}(\ol\smallalpha v)}	\end{equation*}
	is meromorphic, periodic with respect to the lattice $\Lambda$. 
	Here we denote $\iota(\epsilon)v$ by $\epsilon v$ for simplicity. 
	If $\epsilon$ is as in (\ref{equation; epsilon p}) and $v_n$ is 
	 as in Proposition \ref{proposition: parallel for theta}, this function satisfies
	\begin{equation}\label{equation: AFT4}
		f_{v_n}(v) = \pair{ \smallalpha v_n, \smallalpha v_0/2}_\mathscr{L}\,
	{\vartheta^{\cM}_{\smallalpha v_0 + \smallalpha v_n}(\smallalpha v)}/{
		\vartheta^{\cM}_{\ol\smallalpha v_0}(\ol\smallalpha v)}. 
	\end{equation}
\end{lemma}

\begin{proof}
Since $\mathscr{L}$ is  $\Phi$-admissible,   we have 
       $e^{\cM}_{\epsilon u}( \epsilon v)=e^{\cM}_{\ol\epsilon u}( \ol\epsilon v)$ 
and $ \pair{ \smallalpha u, \smallalpha v}_\mathscr{L} 
=\pair{ \ol\smallalpha u, \ol\smallalpha v}_\mathscr{L}$ for $u, v \in V$. 
We also have $\ol \epsilon v_n \in \Lambda$. 
       Then these formulas follow from Proposition 
	\ref{proposition; transformations } ii) and iii). 
	\end{proof}

\begin{proof}[Proof of Proposition \ref{proposition: parallel for theta}]
We apply Lemma \ref{lemma: parallel for rational} to 
 $f$  in the previous lemma. Then 
		$\wh f(t \oplus t_{\fp^n})=
		 \wh{f}_{v_n}(t).$
	By our choice of $\smallalpha$, we have	$[\smallalpha] t_{\fp^n} = t_{\fp^n}$ and 
 $[\ol\smallalpha] t_{\fp^n} = 0$.  This implies
		$$
		\wh\vartheta^{\cM}_{\smallalpha v_0}([\smallalpha]
		(t \oplus t_{\fp^n})) = \wh\vartheta^{\cM}_{\smallalpha v_0}([\smallalpha] t\oplus t_{\fp^n}), 
		\quad 
			\wh\vartheta^{\cM}_{\ol\smallalpha v_0}([\ol\smallalpha]
		(t \oplus t_{\fp^n})) = \wh\vartheta^{\cM}_{\ol\smallalpha v_0}([\ol \smallalpha] t). 
	$$
	As a consequence, \eqref{equation: AFT4}  gives 
	\begin{equation*}
		\wh \vartheta^{\cM}_{\smallalpha v_0} ([\smallalpha] t \oplus t_{\fp^n})
		=  \pair{ \smallalpha v_n, \smallalpha v_0/2}_\mathscr{L}\;
		\wh \vartheta^{\cM}_{\smallalpha v_0 +\smallalpha v_n}([\smallalpha ] t). 		\end{equation*}
	 Our assertion follows by substituting the inverse power series of  $[\smallalpha] t$ into $t$ of the above equality.  
	If $\smallalpha \equiv 1 \pmod{\phi_i(\mathfrak a)}$ for all $i$, we have 
	\begin{equation*}
	\vartheta^{\cM}_{\smallalpha v_0 + \smallalpha v_n}
	(v)=\alpha_\mathscr{L}((\epsilon-1)v_0)\pair{ v_0+\smallalpha v_n, (\smallalpha-1) v_0/2 }_\mathscr{L}
	\vartheta^{\cM}_{ v_0 + \smallalpha v_n}
	(v)
	\end{equation*}
	and 
	\begin{equation*}
	 \vartheta^{\cM}_{\smallalpha v_0 }
	(v)=\alpha_\mathscr{L}((\epsilon-1)v_0)\pair{ v_0, (\smallalpha-1) v_0/2 }_\mathscr{L}
	\vartheta^{\cM}_{ v_0}
	(v).
	\end{equation*}
	The last formula follows from these equations. 
\end{proof}

\begin{corollary}
Let $\fa$ be as in  Proposition \ref{proposition: parallel for theta}. 
Let $v_0, w_0$ be  $\fa$-torsion points of 
$\C^g/\Lambda$.  
Let $s_{\fp^n}$  and  $t_{\fp^n}$ be $\fp^n$-torsion points of the formal group $\wh{\mathcal{A}}$, 
and let $v_n$ and $w_n$ be elements in $\Lambda\otimes \Q$ 
that respectively represents the images of $s_{\fp^n}$  and $t_{\fp^n}$ by $\pi_p$. 
Let $\smallalpha$ be an element of  $\cO_K$ such that 
	$$\smallalpha \equiv 1 \mod{\phi_i(\fp^n)}, 
	\qquad \smallalpha \equiv 0 \mod{ \ol\phi_i(\fp^n)}.$$
Then we have
	$$
		\wh \varTheta_{\smallalpha v_0, \smallalpha w_0} (s \oplus s_{\fp^n}, t \oplus t_{\fp^n} )
		=  \pair{ \smallalpha v_n, \smallalpha w_0}_\mathscr{L}\;
		\wh \varTheta_{\smallalpha v_0+\smallalpha v_n, 
		\smallalpha w_0+\smallalpha w_n  }( s, t). 
			$$
	Moreover,  if   $\smallalpha \equiv 1 
	\mod{\phi_i(\mathfrak{a})}$ for all $i$, we have 
	$$
		\wh \varTheta_{ v_0,  w_0} (s \oplus s_{\fp^n}, t \oplus t_{\fp^n} )
		=  \pair{ \smallalpha v_n, \smallalpha w_0}_\mathscr{L}\;
		 \pair{ \smallalpha w_n, (\smallalpha-1) v_0}_\mathscr{L}\;
		\wh \varTheta_{ v_0+\smallalpha v_n, 
		 w_0+\smallalpha w_n  }( s, t). 
			$$
\end{corollary}
\begin{proof}
Let $\epsilon' \in \cO_K$ be such that $\smallalpha' \equiv 1 \mod{\phi_i(\fp^n)}$ and 
	$\smallalpha' \equiv 0 \mod{2\ol\phi_i(\fp^n)}$. 
	Then the corollary for $\epsilon'$ instead of $\epsilon$ 
	follows directly from Proposition  \ref{proposition: parallel for theta} and 
	Proposition \ref{proposition; transformations }. 	
	The fact that the formulas are invariant after changing  $\epsilon'$ to $\epsilon$ follows again from 
	Proposition \ref{proposition; transformations }.
\end{proof}

\subsection{Relation with Norman's $p$-adic theta function.}\label{subsection:3-5}
%

There are several algebraic and $p$-adic methods to construct power series theta function. 
First, algebraic properties of such power series associated to theta functions 
with algebraic divisors were
systematically studied by Barsotti \cite{Bar}.  Expanding on this theory, modulo $p$ and $p$-adic
properties of such theta functions were studied by Cristante \cite{Cri1} \cite{Cri2},  Candilera-Cristante \cite{CC}, and independently by Norman \cite{Nor} based on the technic of 
Mumford's theory of algebraic theta functions.  
Mazur-Tate also constructed $p$-adic power series 
theta function associated to the divisor $[0]$ of elliptic curves. 

In this subsection, we digress to relate our theta function with the theory of $p$-adic theta functions by Norman.
First we review Norman's theory of $p$-adic theta functions (\cite{Nor}, \S 4). 
Let $R$ be a complete discrete valuation ring of the maximal ideal $\mathfrak m$ 
with residue field $k$  of characteristic $p>0$. 
Let $A$ be an abelian scheme over $R$ of relative dimension $g$ and 
we assume that its special fiber $A \times k$ is {\it ordinary}.
Then there exist abelian schemes $A^{(m)}$ over $R$ and 
isogenies 
$$F_m : A \rightarrow A^{(m)}, 
\quad V_m: A^{(m)} \rightarrow A$$
satisfying
\begin{enumerate}
	\item $F_m \circ V_m$  and $V_m \circ F_m$ are  multiplications by 
		$p^m$ on  $A^{(m)}$ or $A$. 
	\item The degree of $F_m$ and $V_m$ are equal. 
	\item  $V_m$ and $F_m^\vee$ are \'etale. 
\end{enumerate}
We put  $R_n=R/\mathfrak{m}^n$, and let $R_{m,n}$ be the ring 
$R_n[T_1, \dots, T_g]/(T_1, \dots, T_g)^m$. 
Let $t_1, \dots, t_g$ be a local parameter of $A/R$ at the origin and 
let $\Delta(m,n)$ be 
 a diagonal mapping at the origin 
$$\Delta(m,n): \; \mathrm{Spec}\,R_{m,n} \rightarrow 
A \times R_{m,n}, \quad t_i \mapsto T_i.$$  
Since $V_k$ is \'etale, we have a unique lifting $\Delta(m,n)_k$  of $\Delta(m,n)$  
compatible with $m,n$ 
such that $\Delta(1,1)_{k}$ is the zero section of $A\times R_{1,1}$. 
\begin{equation*}
\xymatrix{
 &&& 
\;A^{(k)} \times R_{m,n} \ar[d]^{V_k}  \\
\mathrm{Spec}\,R_{m,n}\;
\ar[rrru]^{\Delta(m,n)_{k}} \;\ar[rrr]_{\Delta(m,n)} &&& \;A \times R_{m,n}. 
}
\end{equation*} 
For simplicity, we sometimes denote  $\Delta(m, n)_k$ again  by $\Delta(m,n)$.   
Then for a line bundle $\mathscr{L}$ on $A$, we consider 
an infinitesimal translation of $\mathscr{L}$ by $\Delta(m,n)$.  

\begin{lemma}
Let $G$ be a finite flat subgroup scheme of $(A^{(k)} \times R_{m,n})[p^N]$ over $R_{m,n}$.  
We assume that $G$ is a connected.  
Let $\Delta$ be a valued point of $G$. 
Then the line bundle 
$ \tau^*_{\Delta}\left(V_N^*\mathscr{L}\right)
 \otimes \left(V_N^*\mathscr{L}\right)^{-1}$ is trivial. 
\end{lemma}
\begin{proof}
The theorem of the square shows that the line bundle  
$$F^*_N(\tau^*_{\Delta}\left(V_N^*\mathscr{L}\right)
 \otimes \left(V_N^*\mathscr{L}\right)^{-1})
\cong \tau^*_{p^N\Delta}\mathscr{L}
 \otimes \mathscr{L}^{-1}
$$
is trivial. Hence 
the correspondence 
$$\Delta \;\longmapsto\; \tau^*_{\Delta}\left(V_N^*\mathscr{L}\right)
 \otimes \left(V_N^*\mathscr{L}\right)^{-1}$$
defines a morphism $G$ to the finite flat group scheme 
$\mathrm{Ker} F_N^\vee$ of $(A^{(N)})^\vee$. 
However, $G$ is connected and $\mathrm{Ker} F_N^\vee$ is \'etale, 
this morphism should be trivial. 
\end{proof}

For sufficiently large $N$  depending on only $m$ and $n$, 
we have $p^N\Delta(m,n)=0$.  Moreover, 
 since $p$ is not equal to $2$, there exists a lifting $\widetilde{\Delta}(m,n): \; \mathrm{Spec}\,R_{m,n} \rightarrow 
A \times R_{m,n}$ such that $2\widetilde{\Delta}(m,n)=\Delta(m,n)$ compatible with 
$m, n$ and 
$\widetilde{\Delta}(1,1)$ is the
zero section.  

We apply the above lemma to the subgroup scheme $G$ of $A \times R_{m,n}$ 
generated by $\widetilde{\Delta}(m,n)$ to construct Mumford's (infinitesimal)
 translation.   
Namely, as in \S \ref{subsection:3-2},  for a symmetric line bundle $\mathscr{L}$, 
by using $V_N$ instead of using the multiplication map by $n$,
we have Mumford's isomorphism 
\begin{equation}\label{equation; mumford iso norman}
V^*_N\tau^\mathcal{M}_{\Delta(m,n)}\;:\; V_N^*(\tau^*_{\Delta(m,n)}\mathscr{L})
\longrightarrow V_N^*\mathscr{L}.
\end{equation} 
We fix a trivialization $\varphi: [0]^*\mathscr{L} \cong \mathbb{A}^1_{R}$. 
Then this induces a trivialization 
$$ [0]^*V^*_N\mathscr{L} = [0]^*\mathscr{L} \cong \mathbb{A}^1_{R},$$ and 
for a section $s$ of $\mathscr{L}$ which is holomorphic around the origin,  the section 
$[0]^*V_N^*\tau^*_{\Delta(m,n)}s$ defines a morphism 
\begin{multline*}
\mathrm{Spec}\,R_{m,n} \xrightarrow{[0]^*V_N^*\tau^*_{\Delta(m,n)}s}
[0]^*V_N^*(\tau^*_{\Delta(m,n)}\mathscr{L}) \\
\xrightarrow{\quad [0]^* V^*_N\tau^\mathcal{M}_{\Delta(m,n)} \quad} [0]^*V_N^*\mathscr{L} 
	\xrightarrow{\qquad\varphi\qquad} \mathbb{A}^1_{R_{m,n}}. 
\end{multline*}
Hence this determines an element of $R_{m,n}$ compatible with respect to $m$ and $n$. 
Taking the limit by $m,n$, 
we obtain an element $\vartheta_s^{\mathcal{N}}$ of $R[[T_1, \dots, T_g]]$. This is 
Norman's $p$-adic theta function associated to $s$.

Now we compare Norman's $p$-adic theta function and 
our reduced theta function $\wh\vartheta_s(t)$. 
We again use the notations and assumptions of sections \ref{subsection:3-1} and \ref{subsection:3-3}. In particular, $A$ has CM. 

\begin{proposition}
Let $D$ be a Cartier divisor of $\mathcal{A}$. 
We assume that $\mathscr{L}:=\mathscr{L}(D)$ is  $\Phi$-admissible and 
 $(-1)^*D=D$ for simplicity.   
Let $s=s_D$ be a section of $\mathscr{L}$ corresponding to $D$. 
We assume that $s$ is holomorphic at the origin, namely, 
the divisor of the poles  does not intersect $(\fP, t_1, \dots, t_g)$. 
Let $U$ be a open neighborhood of the origin of $\mathcal{A}$ which 
trivializes the line bundle $\mathscr{L}$ and suppose that 
$s |_U=g |U$ for a rational function $g$. 
We consider the isomorphism  $\mathscr{L}|_{U}=\cO_{U} g \cong \cO_U$
by the multiplication by  $g^{-1}$ and  
the trivialization $\varphi: [0]^*\mathscr{L} \cong \mathbb{A}^1_{R}$
induced by this  isomorphism. 
Then we have 
$$\vartheta_s^{\mathcal{N}}(T)=\wh \vartheta_s(T)/(\wh g^{-1} \wh \vartheta_s)(0).$$ 
\end{proposition}

\begin{proof}
Let $h$ be the class number of $K$. 
If a natural number $N$ is a multiple of $h$, then 
$p^N$ splits in $\cO_K$ as $p^N=\alpha_N \alpha_N^c$ for $\alpha_N \in \cO_K$ 
which is prime to $\phi_i(\fp)$ ($i=1, \dots, g$) and 
$c$ is a non-trivial element of $\mathrm{Gal}(K/K_0)$. 
Then the kernels of the morphisms $[\alpha_N]$ and $V_N$ are equal (to the 
unique \'etale subgroup scheme of $\mathrm{Ker} \;p^N$ of degree $p^N$) and 
we may assume that $A=A^{(N)}$ and $[\alpha_N]=V_N$.   

We put $D'=2D$, $\mathscr{L}'=\mathscr{L}^{\otimes 2}$ and $s'=s^{\otimes 2}$. 
Then Norman's $p$-adic theta function associated to the section 
$s'$ with trivialization $\varphi^{\otimes 2}$ is given by 
\begin{equation}\label{equation; norman vs ours}
\lim_{m,n \rightarrow \infty}\lim_{N \rightarrow \infty}\;   
f_{D',\Delta(m,n)_N}^{\mathcal{M}}(t,T) 
\cdot \wh g^{-2}([\alpha_N]t)\;
 \vert_{t=0} 
\end{equation}
where $f_{D',\Delta(m,n)_N}^{\mathcal{M}}(t,T)$ is the rational function  
defined by the isomorphism of Mumford in \eqref{equation; mumford iso norman}
($t$ is the parameter of $A$ and $T$ is that of $R_{m,n}$.)   
 The function 
$f_{D',\Delta(m,n)_N}^{\mathcal{M}}$ is explicitly given by 
$$
	 \tau^*_{\widetilde{\Delta}(m,n)_N} \left( \rho_{\widetilde{\Delta}(m,n)_N} (t) 
	\cdot \rho_{\widetilde{\Delta}(m,n)_N} ([-1]t)^{-1} \right),
$$
where $\rho_{\widetilde{\Delta}(m,n)_N}$ is  any rational function  whose divisor is  
$$\tau^*_{\widetilde{\Delta}(m,n)_N}\left(\alpha_N^*D'\right)
- \left(\alpha_N^*D'\right).$$
(Since $(-1)^*D=D$, we took $\rho_{-1}=1$.)
For sufficiently large $N$, we have $\alpha_N^c \Delta(m,n)_N=0$,
and we can take $\rho_{\widetilde{\Delta}(m,n)_N}$ as 
$$
	\tau^*_{\widetilde{\Delta}(m,n)_N}\wh{f}_{\alpha_N}(t) \cdot \wh{f}_{\alpha_N}(t)^{-1},
$$
where ${f}_{\alpha_N}(z)$ is a rational function $a \vartheta_s(\alpha_N z)^2/
\vartheta_s(\alpha^c_N z)^2$ for some constant $a$. 
(The function $\sqrt{f_\alpha}$ might not be rational and 
this is why we consider $D'$ instead of $D$).
Hence we have 
$$
	 \tau^*_{\widetilde{\Delta}(m,n)_N} \left( \rho_{\widetilde{\Delta}(m,n)_N} (t) 
	\cdot \rho_{\widetilde{\Delta}(m,n)_N} ([-1]t)^{-1} \right) = \tau^*_{{\Delta}(m,n)_N}
	\wh f_{\alpha_N}(t)\cdot \wh f_{\alpha_N}(t)^{-1}.
$$
Note that $\wh f_{\alpha}([-1]t)=\wh f_{\alpha}(t)$. 
Therefore   (\ref{equation; norman vs ours}) is equal to
$$\wh \vartheta_s(T)^2/(\wh g^{-1} \wh \vartheta_s)(0)^2. $$
By definition, the number
$(\wh g^{-1}\vartheta_s^{\mathcal{N}})(T)_{\vert T=0}$ must be $1$ and 
$\vartheta_{s'}^{\mathcal{N}}(T)= \vartheta_s^{\mathcal{N}}(T)^2$. 
Hence we have 
$\vartheta_s^{\mathcal{N}}(T)=\wh \vartheta_s(T)/(\wh g^{-1} \wh \vartheta_s)(0)$. 
\end{proof}

Suppose that $A$ is an elliptic curve and $D$ is the divisor $[0]$. 
We consider the section $s=s_D$ which is defined by
$t=-2x/y$ in a neighborhood of the origin. Then the proposition above says that 
Norman's theta function is the pull back by $z=\lambda(t)$ of 
$\theta(z)$.

We assumed $(-1)^*D=D$ in the proposition just for simplicity and it suffices to assume that 
$\mathscr{L}$ is symmetric.

\section{$p$-adic measures and  $p$-adic $L$-functions}\label{section:4}

In \S \ref{subsection:4-1}, we will construct the $p$-adic measure $\mu_{z_0, w_0}$
on $\Z_p \times \Z_p$ and determine the restriction to $\Z^\times_p \times \Z^\times_p$. 
Then  in \S \ref{subsection:5-1} and  \S \ref{subsection:5-2}, we use this measure to construct 
various $p$-adic measures constructed by various authors which were used to construct the
$p$-adic $L$-function interpolating special values of algebraic Hecke characters.
Such $p$-adic measure was first constructed by Manin-Vishik \cite{MV}, and subsequently 
Katz \cite{Ka2} gave an alternative construction using algebraic differentials on modular forms.
From the Galois group side, expanding on the construction of Coates and Wiles \cite{CW1} 
\cite{CW2} of the one variable case,
Yager \cite{Yag1} gave an alternative construction of this $p$-adic measure from a system of elliptic units,
when the class number of $K$ is equal to $1$.
This method has been generalized to the case of higher class numbers 
by various authors (\cite{dS}, \cite{Til}).
Colmez-Schneps \cite{CS} gave another construction of this $p$-adic measure 
by direct analysis of the Eisenstein-Kronecker-Lerch series. 

In \S \ref{subsection:5-1},  we will give a construction using our $\mu_{z_0, 0}$
of the $p$-adic measure defined by Yager \cite{Yag1} used to define his $p$-adic $L$-function. 
In \S \ref{subsection:5-2}, we will give a construction of Katz's $p$-adic $L$-function. 
By using  coefficients of our theta functions, we also give 
an explicit power series giving the Mazur-Swinnerton-Dyer measure. 

\subsection{Construction of the $p$-adic measure}\label{subsection:4-1}
%
%

In this section, we will construct and study the $p$-adic measure $\mu_{z_0,w_0}$ on $\bbZ_p^\times 
\times \bbZ_p^\times$ which interpolates the Eisenstein-Kronecker numbers at $(z_0, w_0)$.  
We keep the notation of \S \ref{subsection:3-3} for CM elliptic curves.  
In particular,  $p \geq 5$ is a prime 
that splits as $(p) = \fp \ol\fp$ in $\cO_K$ and 
the elliptic curve $E$ comes from a fixed integral Weierstrass model $\mathcal{E}$ which 
is good at $p$. 
We first briefly recall the relation 
between $p$-adic measures on $\mathbb{Z}_p \times \mathbb{Z}_p$ and power series $f(S, T) \in W[[S, T]]$. 
See for example \cite{Yag1} \S 6 for details.

\begin{proposition}\label{proposition: measure basic}
	Let $f(S,T) \in W[[S,T]]$ be a two-variable formal power series with coefficients in $W$.   
	Then there exists a 
	unique measure $\mu_f(x,y)$ on $\mathbb{Z}_p \times \mathbb{Z}_p$ with values in $W$ such that
	$$\int_{\mathbb{Z}_p \times \mathbb{Z}_p} (1+S)^x (1+ T)^y d \mu_f(x,y) = f(S,T).$$
	In particular, for any integers $m$, $n \geq 0$, we have the relation
	$$
		\int_{\mathbb{Z}_p \times \mathbb{Z}_p} x^m y^n d\mu_f(x,y) =
		\bigl.  \partial_{S, \log}^m \partial_{T, \log}^n f(S,T) \bigr|_{(S,T)=(0,0)},
	$$
	where $\partial_{X, \log}$ for any variable $X$ is the differential operator $(1+X) \partial_X$.
	This correspondence gives a bijection between the ring of $p$-adic measures on 
	$\mathbb{Z}_p \times \mathbb{Z}_p$ and power series in $W[[S,T]]$.
\end{proposition}

It is known that there exists an isomorphism $\eta_p$ of formal groups
$$
	\eta_p: \widehat{E} \xrightarrow{\cong} \widehat{\mathbb{G}}_m
$$
over $W$. 
We fix an isomorphism $\eta_p$ and then $\eta_p$ is given by a power series 
$$\eta_p(t) = \Exp{\Omega_\fp^{-1}\lambda(t)}-1,$$ 
where  $\Omega_\fp$ is a suitable $p$-adic period $\Omega_\fp \in W^{\times}$.   
We denote by $\iota(T)$ the inverse of $\eta_p(t)$.  In other words, $\iota(T)$ is such that 
$\iota(\eta_p(t)) = t$ and $\eta_p(\iota(T))=T$.   Since $\Omega_\fp \in W^{\times}$, 
the coefficients of $\iota(T)$ are also in $W$.    For any power series $g(s,t) \in W[[s, t]]$ 
on $\wh E \times \wh E$,
we may associate a measure $\mu_g$ on $\bbZ_p \times \bbZ_p$ by taking the measure
associated to the power series $g^\iota(S,T) := g(s, t)|_{s = \iota(S), t=\iota(T)} \in W[[S,T]]$ 
on $\wh{\G}_m \times \wh{\G}_m$.  

In order to define our measure $\mu_{z_0, w_0}$, we first let
\begin{equation}\label{equation: artificial removal}
	\wh\Theta^*_{z_0, w_0}(s,t ; \Gamma):= \wh\Theta_{z_0,w_0}(s,t; \Gamma) - \pair{w_0, z_0} \delta(z_0) s^{-1}
		 - \delta(w_0) t^{-1}.
\end{equation}
Then since $\Gamma$ is the period lattice of $(E, \omega_E)$ 
which comes from a good integral Weierstrass model,  
 Corollary \ref{corollary: for measure} shows
$$
	 \wh\Theta^*_{z_0, w_0}(s,t ; \Gamma) \in W[[s,t]].
$$

\begin{definition}\label{definition: definition of measure} 
        Let $\Gamma$ be the period lattice in $\C$  of $(E, \omega_E)$ 
       which comes from 
     a good integral Weierstrass model at $p$. 
	Let $z_0$, $w_0 \in \Gamma \otimes \Q$.  
	We define $\mu_{z_0, w_0}=\mu_{z_0, w_0}(x,y ; \Gamma)$ to be the measure on $\bbZ_p \times \bbZ_p$ corresponding to 
	the power series
	$$
		\wh\Theta^{* \iota}_{z_0, w_0}(S, T; \Gamma): = 
		\wh\Theta^*_{z_0, w_0}(s,t ; \Gamma)\bigl.\bigr|_{s = \iota(S), t=\iota(T)}. 
	$$
\end{definition}

If there is no fear of confusion, we omit $\Gamma$ in the notation. 
When $z_0$, $w_0 \not\in\Gamma$, then we have $\wh\Theta^*_{z_0, w_0}(s,t) = \wh\Theta_{z_0, w_0}(s,t)$, and 
we have 
$$
		\partial_{S, \log}^{b-1} \partial_{T, \log}^a \wh\Theta^{\iota}_{z_0, w_0}(S,T) \bigl. \bigr|_{S=T=0} =
		\Omega^{a+b-1}_\fp \partial_{z}^{b-1} \partial_{w}^a \Theta_{z_0, w_0}(z, w) \bigl. \bigr|_{z=w=0}.
	$$
Hence our $p$-adic measure in this case satisfies the following interpolation property.

\begin{proposition}\label{proposition: M1}
	Let $z_0$, $w_0 \in (\Gamma \otimes \Q) \setminus \Gamma$ be such that the orders of
	$z_0$ and $w_0$ modulo $\Gamma$ is prime to $p$.
	Then for integers $a \geq 0$, $b>0$, we have
 	$$
		\frac{1}{\Omega_\fp^{a+b-1}} \int_{\bbZ_p \times \bbZ_p} x^{b-1} y^a d \mu_{z_0, w_0}(x,y) = (-1)^{a+b-1}
		(b-1)! \frac{e^*_{a,b}(z_0, w_0)}{A^a}.
	$$
\end{proposition}

In the definition of $\wh\Theta_{z_0, w_0}^*(s,t)$, we are subtracting factors such as $s^{-1}$ and 
$t^{-1}$ instead of $z^{-1}=\lambda(s)^{-1}$ and $w^{-1}=\lambda(t)^{-1}$.  Hence the interpolation property satisfied
by $\mu_{z_0, w_0}$ for $z_0 \in \Gamma$ or $w_0 \in \Gamma$ is not immediately clear. We will
instead calculate the interpolation property of the measure 
which is obtained  by the restriction to $\bbZ_p^\times
\times \bbZ_p^\times$.   


\begin{lemma}\label{lemma: integral formula}
	We have
	\begin{multline}\label{equation: integral formula}
		\int_{\bbZ_p^\times \times \bbZ^\times_p} (1+S)^x (1+T)^y d \mu_{z_0, w_0}(x,y) \\
		= \wh\Theta^\iota_{z_0,w_0}(S, T) - \frac{1}{p} \sum_{S_1} %
		\wh\Theta^\iota_{z_0, w_0}(S \oplus S_1, T)\\
		- \frac{1}{p} \sum_{T_1} \wh\Theta^\iota_{z_0, w_0}(S,  T \oplus T_1)
		+ \frac{1}{p^2} \sum_{S_1, T_1} \wh\Theta^\iota_{z_0, w_0}(S \oplus S_1,  T \oplus T_1),
	\end{multline}
	where the sum is taken over $p$ torsion points $S_1$ and $T_1$ of $\wh\bbG_m$.
\end{lemma}

Note that the functions appearing on the right hand side is $\wh\Theta^\iota$, not $\wh\Theta^{* \iota}$. 

\begin{proof}
For a measure $\mu_f$ on $\bbZ_p$ associated to the
	one-variable power series $f(T) \in W[[T]]$, 
	 the restriction of $\mu_f$
	to $\bbZ_p^{\times}$  is associated to the power series
	$$
	 f(T) - \frac{1}{p} \sum_{T_1 \in \wh\bbG_m[p]} f(T \oplus T_1) \in W[[T]],
	$$
	where the sum is taken over all $p$ torsion points $T_1$ of $\wh\bbG_m$  
	(see for example Lang \cite{Lan1} Chapter 4 \S 2, \textbf{Meas 4}.)
	By applying this to the two variable case, we have the assertion, except for the
	terms appearing in the difference between $\wh\Theta$ and $\wh\Theta^*$.  
	They are differ by the poles coming from $\iota(S)^{-1}$ and  $\iota(T)^{-1}$, but 
	direct
	calculation shows that the term involving $\iota(S)^{-1}$
	add up as
	$$
	 \left( \frac{1}{\iota(S)} -
	\frac{1}{p} \sum_{S_1} \frac{1}{\iota(S \oplus S_1)} - \frac{1}{p} \sum_{T_1} \frac{1}{\iota(S)} 
	+ \frac{1}{p^2} \sum_{S_1, T_1} \frac{1}{\iota(S \oplus S_1)} \right) = 0,$$
	and similarly for terms involving $\iota(T)^{-1}$. 
\end{proof}

In the case of the Kubota-Leopoldt $p$-adic $L$-function associated to the special values of
the Riemann zeta function, it is necessary to take an auxiliary integer $c$ to remove the pole
of the generating function.  
Similarly in the case of CM elliptic curves, one had to take an auxiliary choice of an ideal $\fa$ to 
remove the poles  and afterwards, one had to eliminate again this auxiliary ideal. 
The above proposition shows that by simultaneously restricting the
two-variable $p$-adic measure to $\bbZ_p^\times \times \bbZ_p^\times$, the poles 
automatically disappear.  This enable us to avoid auxiliary ideals and justifies the artificial``removal of poles" in 
\eqref{equation: artificial removal} used to define the measure.

\begin{proposition}\label{proposition: M2}
	Let $z_0$, $w_0 \in \Gamma_\bbQ$ be such that $z_0 \mathfrak{a}$ and $ w_0 \mathfrak{a} $ are in $\Gamma$ and $N\!\mathfrak a$  is prime to $p$.   Then 
	\begin{multline*}
			\int_{\bbZ_p^\times \times \bbZ_p^\times}   (1+S)^x (1+T)^y d \mu_{\smallalpha z_0, \smallalpha w_0}(x,y)\\
			= \biggl[\biggr.\Theta_{\smallalpha z_0, \smallalpha w_0}(z, w; \Gamma)
			 - \Theta_{p \smallalpha z_0, \epsilon w_0}\left(p z, w ;\ol{\fp}\,\Gamma\right)			 - \Theta_{\smallalpha z_0, p\smallalpha w_0} \left(z, p w; \ol{\fp}\,\Gamma\right)\\%
			+\Theta_{p \smallalpha z_0, p \smallalpha w_0}\left(p z, 
			p w ; \ol{\fp}^2 \Gamma\right)
			\biggl. \biggr]\biggl. \biggr|_{z = \lambda^\iota(S), w=\lambda^\iota(T)},
	\end{multline*}
	where $\smallalpha \in \cO_K$ is such that $\smallalpha \equiv 1 \pmod{ \fp}$ and
	$\smallalpha \equiv 0 \pmod{\ol\fp}$. If in addition $\smallalpha \equiv 1 \pmod{ \mathfrak a}$, we have 
\begin{multline*}
			\int_{\bbZ_p^\times \times \bbZ_p^\times}   (1+S)^x (1+T)^y d \mu_{ z_0, w_0}(x,y)\\
			= \biggl[\biggr.\Theta_{ z_0,  w_0}(z, w; \Gamma)
			 - \Theta_{p  z_0, \epsilon w_0}\left(p z, w ;\ol{\fp}\,\Gamma\right)			 -\pair{(\epsilon-1)z_0, w_0}_\Gamma \Theta_{\smallalpha z_0, p w_0} \left(z, p w; \ol{\fp}\,\Gamma\right)\\%
			+\pair{(\epsilon-1)z_0, w_0}_\Gamma\Theta_{p \smallalpha z_0, p \smallalpha w_0}\left(p z, 
			p w ; \ol{\fp}^2 \Gamma\right)
			\biggl. \biggr]\biggl. \biggr|_{z = \lambda^\iota(S), w=\lambda^\iota(T)},
	\end{multline*}
		\end{proposition}

\begin{proof}
	By definition, we have
	\begin{equation*}
		\wh\Theta_{\smallalpha z_0, \smallalpha w_0}(s, t) = \Theta_{\smallalpha z_0,\smallalpha w_0}(z,w)%
			 \bigl.\bigr|_{z=\lambda(s),w=\lambda(t)}.
	\end{equation*}
	By Proposition \ref{proposition; distribution} and 
	Proposition \ref{proposition: parallel for theta}, we have 

\begin{align*}
			\sum_{s_1} \wh\Theta_{\smallalpha z_0,\smallalpha w_0}(s\oplus s_1, t) %
			&= p\sum_{z_1 \in \fp^{-1}\Gamma/\Gamma} 
			\pair{\smallalpha z_1, \smallalpha w_0}
			\Theta_{\smallalpha z_0 + \smallalpha z_1, \smallalpha w_0}(z,w)|_{z = \lambda(s), w=\lambda(t)}\\
			&= p\Theta_{p \smallalpha z_0, \epsilon w_0}\left(p z, w ;\ol{\fp}\,\Gamma\right)%
			 \bigl.\bigr|_{z = \lambda(s), w=\lambda(t)}\nonumber. 
	\end{align*}
	Similarly, we have 
	\begin{align*}
			\sum_{t_1}\wh\Theta_{\smallalpha z_0,\smallalpha w_0}  \left(s,t\oplus  t_1 \right)
			=p\Theta_{\smallalpha z_0, p\smallalpha w_0} \left(z, p w; \ol{\fp}\,\Gamma\right)
			 \bigl.\bigr|_{z = \lambda(s), w=\lambda(t)}\nonumber,
	\end{align*}
	\begin{align*}
		\sum_{s_1, t_1} \wh\Theta_{\smallalpha z_0,\smallalpha w_0} 
		 (s \oplus s_1, t \oplus t_1) 
			=p^2
			\Theta_{p \smallalpha z_0, p \smallalpha w_0}\left(p z, 
			p w ; \ol{\fp}^2 \Gamma\right)
			 \bigl.\bigr|_{z = \lambda(s), w=\lambda(t)}\nonumber.
	\end{align*}
	The first assertion now follows from the previous lemma. The last 
	assertion follows from the first assertion and Proposition 
	\ref{proposition; transformations }. 
\end{proof}

%
\subsection{Relation to the $p$-adic measure of Yager}\label{subsection:5-1}
%
%
%

Yager constructed a certain two variable $p$-adic measure that 
interpolates the special values of the Hecke $L$-function of a 
Gr\"ossencharakter $\varphi$ of an imaginary  quadratic  field $K$ whose  class number  is equal to $1$.   

Let $\varphi$ be a Hecke character of infinity type $(1,0)$ 
of an imaginary  quadratic  field $K$ of 
the class number $1$. Let $\ff$ be its conductor.  
Let $\Omega$ be a complex number such that $\Gamma=\Omega\ff$ 
is a period lattice of a Weierstrass integral model of $E$ over $\cO_K$. 
For a fractional ideal $\fa$, we put $\Gamma_\fa=\varphi(\fa)\fa^{-1}\Gamma$. 

\begin{definition}
Let $\alpha_\fa$ be an element of $\fa^{-1} \cap P_K(\ff)$. 
	We define the $p$-adic measure $\mu_\varphi$ by
	\begin{align*}
		\mu_\varphi(x,y) :=\sum_{\fa \in I_K(\ff)/P_K(\ff)}
			\mu_{\varphi(\alpha_\fa\fa)\Omega, 0}(x,N(\ff)y ; \Gamma_{\fa}). 
	\end{align*}
	\end{definition}
	
	Note that since we assumed the class number of $K$ is $1$, the value of 
	$\varphi$ is in $K$ and $\Gamma_{\fa}=\Gamma$. 
	Hence we can use the same integral models and $p$-adic periods for 
	all elliptic curves $\C/\Gamma_{\fa}$. 
	This is the reason that we assume the class number of $K$ is equal to $1$.
	We define the measure $\mu_\varphi$ as 
	\begin{align*}
		\mu_\varphi(x,y) =\sum_{\alpha \in (\cO_K/\ff)^\times}
			\mu_{\varphi(\alpha)\Omega, 0}(x,N(\ff)y ; \Gamma). 
	\end{align*}
		The function $\Theta_{z_0, 0}(z,w; \Gamma)$ is periodic in $z_0$ 
		with respect to the lattice $\Gamma$. Hence the definition
	of $\mu_\varphi$ does not depend on the choice of the representative 
	$\fa$ and  $\alpha_\fa$.

\begin{theorem}\label{theorem: Yager}
Let $a,b$ be  integers such that $b>a \geq 0$.  
	Then we have
	\begin{multline*}
		\frac{1}{\Omega_\fp^{a+b}} \int_{\bbZ_p^\times \times \bbZ_p^\times}
		x^{b-1} y^a d \mu_{\varphi} (x,y) \\
		= (-1)^{a+b-1} (b-1)!
		\left[ \frac{2\pi}{\sqrt{d_K}}\right]^a
		\left(1 - \frac{\varphi(\fp)^{a+b}}{N\!\fp^{a+1}} \right) \left( 1 - \frac{\ol\varphi(\ol\fp)^{a+b}}{N\!\ol\fp^{b}}\right)
		\frac{L_\ff(\ol\varphi^{a+b}, b)}{\Omega^{a+b}}.
	\end{multline*}
\end{theorem}

\begin{proof}
	Let $\smallalpha \in \cO_K$ such that $\smallalpha \equiv 1 \pmod{\ff\fp}$ and
	$\smallalpha \equiv 0 \pmod{\ol\fp}$. 
	Then by Proposition \ref{proposition: M2}, we have 
	\begin{multline*}
			\int_{\bbZ_p^\times \times \bbZ_p^\times}   (1+S)^x (1+T)^y d \mu_{\varphi(\alpha_\fa \fa)\Omega, 0}(x,y; \Gamma_{\fa})\\
			= \biggl[\biggr.\Theta_{\varphi(\alpha_\fa\fa)\Omega,  0}(z, w;\Gamma_{\fa})
			 - \Theta_{p\varphi(\alpha_\fa \fa)\Omega, 0}
			\left(p z, w ;\ol{\fp}\,\Gamma_{\fa}\right)	
			- \Theta_{\smallalpha \varphi(\alpha_\fa \fa)\Omega, 0} \left(z, p w; \ol{\fp}\,\Gamma_{\fa}\right)\\%
			+\Theta_{p \smallalpha\varphi(\alpha_\fa\fa)\Omega, 0}\left(p z, 
			p w ; \ol{\fp}^2 \Gamma_{\fa}\right)
			\biggl. \biggr]\biggl. \biggr|_{z = \lambda^\iota(S), w=\lambda^\iota(T)},
	\end{multline*}	
	We have 
	\begin{multline*}
			 \Theta_{p\varphi(\alpha_\fa\fa)\Omega, 0}
			\left(p z, w ;\ol{\fp}\,\Gamma_{\fa}\right)	=
		{p}^{-1}{\varphi(\fp)}
			 \Theta_{\varphi(\alpha_\fa\fa\fp)\Omega, 0}
			\left(\varphi(\fp) z, {p}^{-1}{\varphi(\fp)}w ;\,\Gamma_{\fa\fp}\right). 
							\end{multline*}	
Hence by Proposition \ref{proposition: EK and L}  i) and 
Theorem \ref{theorem: generating function theorem}, we have
\begin{multline*}
\sum_{\fa \in I_K(\ff)/P_K(\ff)} \partial_z^{b-1}\partial_w^{a}
 \Theta_{p\varphi(\alpha_\fa\fa)\Omega, 0} 
			\left(p z, w ;\ol{\fp}\,\Gamma_{\fa}\right)	\vert_{z=w=0}\\=
			(-1)^{a+b-1} \frac{(b-1)!}{A(\Gamma)^a} \,  \frac{\varphi(\fp)^{a+b}}{N\!\fp^{a+1}}\,
			\sum_{\fa \in I_K(\ff)/P_K(\ff)} 
			 e^*_{a,b}(\varphi(\alpha_\fa\fa\fp)\Omega, 0
			 ; \Gamma_{\alpha\fa\fp})\\
			=(-1)^{a+b-1} |\mu_K| \frac{(b-1)!}{A(\ff)^a} \,  \frac{\varphi(\fp)^{a+b}}{N\!\fp^{a+1}}\,
			 \frac{L(\ol\varphi^{a+b}, b)}{\Omega^{a+b}}. 
\end{multline*}
Similarly, we have 
\begin{multline*}
\sum_{\fa \in I_K(\ff)/P_K(\ff)} \partial_z^{b-1}\partial_w^{a} \Theta_{\epsilon \varphi(\alpha_\fa\fa)\Omega, 0} 
			\left(z, p w ;\ol{\fp}\,\Gamma_{\fa}\right)	\vert_{z=w=0}\\			=(-1)^{a+b-1} |\mu_K| \frac{(b-1)!}{A(\ff)^a} \,  \frac{\ol\varphi(\ol\fp)^{a+b}}{N\!\fp^{b}}\,
			 \frac{L(\ol\varphi^{a+b}, b)}{\Omega^{a+b}}
\end{multline*}	
and 
\begin{multline*}
\sum_{\fa \in I_K(\ff)/P_K(\ff)} \partial_z^{b-1}\partial_w^{a} \Theta_{p\epsilon \varphi(\alpha_\fa\fa)\Omega, 0} 
			\left(pz, p w ;\ol{\fp}^2\,\Gamma_{\fa}\right)	\vert_{z=w=0}\\			=(-1)^{a+b-1} |\mu_K| \frac{(b-1)!}{A(\ff)^a} \,  
			\frac{\varphi(\fp)^{a+b}\ol\varphi(\ol\fp)^{a+b}}{N\!\fp^{a+b+1}}\,
			 \frac{L(\ol\varphi^{a+b}, b)}{\Omega^{a+b}}. 
\end{multline*}	
	(Note that $\varphi(\smallalpha)=\smallalpha$). 
	Since 
	$\partial_{S, \log} = \Omega_\fp \partial_z$ and $\partial_{T, \log} = \Omega_\fp \partial_w$,  
	our assertion now follows from the above.
\end{proof}

The above interpolation property shows that the measure $\mu_\varphi$ defined above on 
$\bbZ_p^\times \times \bbZ_p^\times$ is in fact the $p$-adic measure used by Yager  
to define his $p$-adic $L$-function (See \cite{Yag1} \S 9).

\subsection{Relation to the $p$-adic measure of Katz}\label{subsection:5-2}
%
%
%

First we introduce some notations and definitions following 
\cite{GS}.  See also \cite{dS}, Chapter 2. 
Let $\ff$ be an integral ideal of $\cO_K$ such that 
$w_\ff=1$ where $w_\ff$ is the number of units congruent to $1$ modulo $\ff$. 
Let $p$ be a rational prime which splits in $\cO_K$ as $p=\fp \ol\fp$ and 
$(p, N\!\ff)=1$. 
Let $F$ be the ray class field $K(\ff)$. 
Then there exists a Hecke character $\varphi_0$ of $K$ whose conductor is equal to $\ff$ and 
an elliptic curve 
over $E/F$ whose Gr\"o{\ss}encharakter $\psi_{E/F}$ is 
written as $\psi_{E/F}=\varphi_0 \circ N_{F/K}$. 
 Since $F(E_{\mathrm{tor}})/K$ is abelian, we have 
$\psi_{E^\sigma}= \psi_{E}$ for $\sigma \in \mathrm{Gal}(F/K)$. 
Hence considering its conjugate if necessary we may assume that 
the period lattice for an invariant differential for $E/F$ is given by 
$\Gamma:=\ff\Omega $ for some complex number $\Omega$. 
Moreover, by changing $\Omega$ by a constant multiple, 
we may also assume that the Weierstrass equation for $(E, \omega_E)$  
gives a proper smooth model $\mathcal{E}/\cO_{K_\fp}$ at $\fp$. 
Let $B$ be the Weil restriction of $E$ over $F$ to $K$ 
and $\wt \varphi$  the Serre-Tate character $\wt \varphi : K_{\mathbb{A}}^\times 
\rightarrow \mathrm{End}_K B \otimes \Q$. 
Since $\mathrm{End}_K B \cong \prod_{\sigma \in \mathrm{Gal}(F/K)} 
\mathrm{Hom}(E, E^\sigma)$, for an integral ideal $\fa$ of $K$ prime to $\ff$ 
the Serre-Tate character $\wt \varphi_0$  induces  an isogeny 
$\wt{\varphi}_0(\fa): E \rightarrow E^\sigma$. 
We define $\Lambda(\fa) \in F$ to be the element such that 
$\wt{\varphi}_0(\fa)^*(\omega_E^{\sigma_\fa})=\Lambda(\fa) \omega_E$ 
where $\sigma_\fa=(\fa, F/K)$. 
Then $\Lambda$ satisfies the cocycle condition 
$\Lambda(\fa\mathfrak{b})
=\Lambda(\fa)^{\sigma_\mathfrak{b}} \Lambda(\mathfrak{b})
$ and 
this enable us to extend the definition of $\Lambda$ to 
any fractional ideal prime to $\ff$ so that this cocycle condition 
remains valid. 
Then the period lattice $\Gamma^{\sigma_\fa}$
corresponding to $(E^{\sigma_\fa}, \omega_E^{\sigma_\fa})$ is 
given by $\Gamma^{\sigma_\fa}=\Lambda(\fa)\fa^{-1}\Gamma$. 
As before, let $W$ be the Witt ring $W(\ol{\mathbb{F}_p})$ and 
we fix a prime $\fP$ of $F$ over $\fp$ and an embedding 
$F \subset F_\fP \hookrightarrow W(\ol{\mathbb{F}_p})\otimes \Q$. 
We let $\sigma_\fa=(\fa, F(E[\ol\fp^\infty])/K)$ for $\fa$ such that$(\fa, \ff \ol\fp)=1$. 
Then the action of $\sigma_\fa$ is extended to $W$ continuously. 
There exists a $p$-adic period $\Omega_p \in W^\times$ such that 
$\Omega_p^{\sigma_\fa}=\Lambda(\fa)N\!\fa^{-1} \Omega_p$ 
if $(\fa, \ff \ol\fp)=1$, and this gives an isomorphism 
$\eta_{\wh E} : \wh E \cong \wh{\mathbb{G}}_m $ over $W$ by 
$\eta_{\wh E}(t)=\exp(
\Omega_p^{-1}\lambda_{\wh E}(t))$. 
The $p$-adic period is unique up to an element of $\Z_p^\times$ 
and we take it so that we have 
$$\zeta_m:=\pair{1, \epsilon}_{\fp^m}=\eta_{\wh E^\sigma}(t) \vert_{t=t_{\fp^m}}$$
where $\epsilon$ is an element of $\cO_K$ such that 
$\epsilon \equiv 1 \mod \fp^m \ff$ and  $\epsilon \equiv 0 \mod \ol\fp^m$, 
the element $\sigma$ is $\sigma_{\fp^{-m}}$ and 
$t_{\fp^m}$ is a $\fp^m$-torsion point of $\wh E^\sigma$ obtained as 
the image of 
$$\Lambda(\fp^{-m})\nu_m\Omega \; \in \; 
\Gamma^\sigma\otimes \Q$$ 
where $\nu_m$ be an element of $\ff$ such that 
$\nu_m \equiv 1 \mod \fp^m$ 
through the isomorphism
$(\Gamma^\sigma\otimes\Q)/\Gamma^\sigma \cong
E^\sigma(\ol\Q)_{\mathrm{tor}}$ defined by the correspondence
$$z \mapsto (\wp(z; \Gamma^\sigma), \wp'(z; \Gamma^\sigma)). $$

Let $\fa$ be a fractional ideal of $K$ prime to $\ff \ol\fp$, and  
let $z_0, w_0 \in \C$ be torsion points of $\C/\fa^{-1}\Gamma$ whose orders are prime to $p$. 
In \S \ref{subsection:4-1}, we defined the Kronecker theta 
measure $\mu_{z_0, w_0}(x,y; \Gamma)$ for a $\Gamma$ 
which comes from a good Weierstrass integral model at $p$. 
Now we define a $p$-adic measure $\mu_{z_0, w_0}(x,y; \fa^{-1}\Gamma)$ even if 
$\fa^{-1}\Gamma$ does not come from a good model. 

We consider the formal power series compositions of the Laurent expansion at 
the origin of the 
Kronecker theta function 
$$\Theta^\iota_{\Lambda(\fa)z_0, \Lambda(\fa)w_0}
(z, w; \Gamma^{\sigma_\fa}) $$ and 
the power series 
$$z=\Omega_p^{\sigma_\fa}\log(1+S), \quad w=\Omega_p^{\sigma_\fa}\log(1+T).$$ 
We denote this power series by 
$\wh{\Theta}^\iota_{\Lambda(\fa)z_0,\Lambda(\fa)w_0}(S, T; \Gamma^{\sigma_\fa})$. 
Since $\Gamma^{\sigma_\fa}$ is the period lattice of $(E^{\sigma_\fa}, \omega_E^{\sigma_\fa})$ 
whose Weierstrass model gives a proper smooth model 
$\mathcal{E}^{\sigma_\fa}/\cO_{K_\fp}$ at $\fp$, and 
$\Omega_p^{\sigma_\fa}$ is a $p$-adic period for $\mathcal{E}^{\sigma_\fa}/\cO_{K_\fp}$, 
as we have seen in \S \ref{subsection:4-1}, we have 
$$\wh{\Theta}^\iota_{\Lambda(\fa)z_0,\Lambda(\fa)w_0}(S, T; \Gamma^{\sigma_\fa}) \in W((S,T)).$$ 
We define the measure $\mu_{z_0, w_0}(x,y ;\fa^{-1}\Gamma)$ as the measure corresponding to 
the power series 
$$\Lambda(\fa)\wh{\Theta}^{* \iota}_{\Lambda(\fa)z_0,\Lambda(\fa)w_0}(S, T; \Gamma^{\sigma_\fa}) \in W[[S,T]]\otimes \Q.$$ 
 If $\fa$ is prime to $\fp$, then the above power series is in $W[[S,T]]$. 
We remark that 
$\fa$ does not have to be prime to $\fp$ but prime to $\ff \ol\fp$. 
Suppose that $\fa$ is prime to $\ff p$. 
Then we consider  
another measure  $\wt\mu_{z_0, w_0}(x,y ;\fa^{-1}\Gamma)$ corresponding to 
the power series 
$$\Lambda(\fa)\wh{\Theta}^{* \iota}_{\Lambda(\fa)z_0,\Lambda(\fa)w_0}([N(\fa)]S, [N(\ff)]T; \Gamma^{\sigma_\fa}) \in W[[S,T]]\otimes \Q.$$

\begin{definition}
For a fractional ideal $\fa$ prime to $\ff p$ and elements $z_0, w_0 \in \C$ which 
are  torsion points of $\C/\fa\Gamma$ whose orders are prime to $p$. 
We define the Katz measure 
$\mu_{z_0, w_0}^{\mathscr{K}}(\alpha ; \fa\Gamma)$ 
 on $(\cO_K \otimes \Z_p)^\times$ 
by 
\begin{align*}
\int_{(\cO_K \otimes \Z_p)^\times} f(\alpha) &\;d\mu_{z_0, w_0}^{\mathscr{K}}(\alpha ; \fa\Gamma) \\
&:=
{\Omega_p}\int_{\Z_p^{\times} \times \Z_p^{\times}} f(x,y^{-1}) x^{-1}d\wt\mu_{z_0, w_0}(x,y ; \fa\Gamma) 
\end{align*}
for a $W$-valued continuous function $f$ on $(\cO_K \otimes \Z_p)^\times$ which 
is regarded as the function $f(x,y)$ on 
$\Z_p \times \Z_p$ through the isomorphism 
\begin{equation}\begin{split}\label{equation: identification}
  \cO_K \otimes \bbZ_p  \cong \cO_{K_{\fp}} \times 
	\cO_{K_{\ol\fp}} &\xrightarrow{\cong}
	\cO_{K_{\fp}} \times 
	\cO_{K_{\fp}} \cong \bbZ_p \times \bbZ_p \\
	(\alpha_1, \alpha_2) & \mapsto (\alpha_1, \ol\alpha_2). 
\end{split}\end{equation}
\end{definition}

Then
\begin{align}\label{scalar change}
\int_{(\cO_K \otimes \Z_p)^\times} f(a\alpha) \;d\mu_{az_0, aw_0}^{\mathscr{K}}(\alpha ; a\fa\Gamma) =
\int_{(\cO_K \otimes \Z_p)^\times} f(\alpha) \;d\wt\mu^{\mathscr{K}}_{z_0, w_0}(\alpha; \fa \Gamma) 
\end{align}
for an element of $c \in K^\times$ prime to $\ff p$. 
Hence the Katz measure essentially depends only on the 
ideal class of $\fa$. 

As before, 
let $I_K(\fa)$ be the ideal group of $K$ consists of fractional ideals prime to $\fa$ 
and
let $P_K(\fa)$ be the set of principal ideals $(\alpha)$ such that 
$\alpha \equiv 1 \mod \fa$. If there is no fear of confusion,  
we also denote the set $\{\alpha \in K | \alpha \equiv 1 \mod \fa\}$ by $P_K(\fa)$. 

\begin{definition}
Let $\rho$ be a  function on $I_K(\ff)$ to $\C_p$ such that 
for a fixed $\fa \in I_K(\ff)$ the function 
$\alpha \mapsto \rho(\fa\alpha)$ on $P_K(\ff)\cap I_K(p)$ can be extended to 
a continuous function on $(\cO_K \otimes \Z_p)^\times$.
This extension is unique if it exists, and we denote it by $\rho^{(p)}(\fa\alpha)$.
For a fractional ideal $\fa$ of $K$ prime to $\ff p$ 
we define the partial $p$-adic $L$-function $\mathcal{L}_{p,\ff}(\rho, \fa)$ by 
$$\mathcal{L}_{p, \ff}(\rho, \fa):=
\int_{(\cO_K \otimes \Z_p)^\times} \rho^{(p)}(\fa\alpha) \;d\mu_{\alpha_0\Omega, 0}^{\mathscr{K}}(\alpha ; \fa\Gamma)
$$
and  the $p$-adic $L$-function $\mathcal{L}_{p, \ff}(\rho)$ by 
$$\mathcal{L}_{p, \ff}(\rho):=\sum_{\fa \in I_K(\ff)/P_K(\ff)}
\mathcal{L}_{p, \ff}(\rho, \fa)
$$
where $\alpha_0$ is any element of $\fa\cap P_K(\ff)$. 
The partial $p$-adic $L$ function does not depend on the choice of the representative of 
$\alpha_0$ and depends only on the ideal class of $\fa$ in $I_K(\ff)/P_K(\ff)$.  
\end{definition}

Let $\varphi$ be a Hecke character  of $K$ of infinity type $(a,b)$ whose conductor 
divides $\ff p^{\infty}$. 
Then there exists a finite character 
$\chi_{p}: (\cO_K/p^m)^\times \rightarrow \ol\Q^\times$  such that 
$\varphi(\alpha)=\chi_{p}(\alpha)\alpha^a \ol\alpha^b$ on $P_K(\ff)$, and 
 the function 
$\alpha \mapsto \varphi(\fa\alpha)$ 
on $P_K(\ff)\cap I_K(p)$ can be extended 
continuously  to $(\cO_K \otimes \Z_p)^\times$.  Therefore  
the value $\mathcal{L}_{p,\ff}(\varphi)$ is defined. 
Now we calculate this value.  
The  character $\chi_{p}$ is of the form $\chi_{\fp} \chi_{\ol \fp}$,
where $\chi_{\fp}$ (resp. $\chi_{\ol \fp}$) is a character of $\cO_K^\times$ whose conductor divides 
$\fp^m$ (resp. ${\ol \fp}^m$). 
Let  $\chi_1$ and $\chi_2$ be Dirichlet characters 
on $(\Z/p^m\Z)^\times$ such that 
$\chi_\fp(\alpha)=\chi_1(\alpha)$ and $\chi_{\ol\fp}(\alpha)=\chi_2^{-1}(\ol\alpha)$,
where we identify $(\cO_K/\fp^m)^\times$ with $(\Z/p^m\Z)^\times$ by 
the natural projection. 
We put $\wt\chi_1=\chi_1 \circ N_{K/\Q}$. 
Then
the partial $p$-adic $L$-function $\mathcal{L}_{p,\ff}(\varphi, \fa^{-1})$ is given by 
\begin{align*}
&\Omega_p\varphi(\fa^{-1})\int_{ \Z_p^\times \times \Z_p^\times} \chi_1(x)\chi_2(y) x^{k-1}y^l \;d\wt\mu_{\alpha_0\Omega, 0}(x,y ; \fa^{-1}\Gamma)\\
&= \Omega_p \varphi(\fa^{-1})
N(\fa)^{k-1} N(\ff)^{l}\\
&\qquad \qquad \qquad \times \int_{ \Z_p^\times \times \Z_p^\times} \chi_1(N(\fa)  x)\chi_2(N(\ff)y) x^{k-1}y^l \;d\mu_{\alpha_0\Omega, 0}(x,y ; \fa^{-1}\Gamma)\\
&= \Omega_p
 (\varphi^{-1}\wt\chi_1)(\fa)  N(\fa)^{k-1}N(\ff)^{l} \\
&\qquad \qquad \qquad \qquad \times \int_{ \Z_p^\times \times \Z_p^\times} \chi_1(x)\chi_2(N(\ff)y) x^{k-1}y^l \;d\mu_{\alpha_0\Omega, 0}(x,y ; \fa^{-1}\Gamma). 
\end{align*}
The above formulas are valid for a fractional ideal $\fa$ prime to $\ff p$, but 
we can also give a natural meaning to the last formula for a fractional ideal prime only to $\ff \ol\fp$.  
In fact, the measure $\mu_{\alpha_0\Omega, 0}(x,y ; \fa^{-1}\Gamma)$ is  defined for 
a fractional ideal $\fa$ prime to $\ff \ol\fp$ and we can define 
$\varphi\wt\chi_1^{-1}(\fp)$ as 
$$\varphi\wt\chi_1^{-1}(\fp):=\varphi\wt\chi_1^{-1}(\fa_0)\chi_1^{-1}(\ol c)\chi_2^{-1}(\ol c)c^k \ol c^{-l}$$ where 
$\fp=c\fa_0$ for some $c \in P_K(\ff)$ and $(\fa_0, p)=1$. 
(Note that $\ol c$ is prime to $\fp$ and the right hand side in the definition of 
$\varphi\wt\chi_1^{-1}(\fp)$ is well-defined. 
Note also that  if $\alpha \in P_K(\ff)$ is prime to $p$, 
we have $\varphi\wt\chi_1^{-1}(\alpha)
=\chi_1^{-1}(\ol \alpha)\chi_2^{-1}(\ol \alpha)\alpha^k \ol \alpha^{-l}$).
We extend the definition of $\varphi\wt\chi_1^{-1}(\fp)$ 
to any fractional ideal $\fa$ prime to $\ff \ol\fp$ multiplicatively.  

Hence we can extend the definition of the partial $p$-adic $L$-function $\mathcal{L}_{p,\ff}(\varphi, \fa^{-1})$ to any fractional ideal $\fa$ prime to $\ff \ol \fp$, and this depends only on the ideal class of $\fa$ in $I_K(\ff)/P_K(\ff)$. 
Therefore we calculate 
$$\mathcal{L}_{p,\ff}(\varphi)
=\sum_{\fa \in I_K(\ff)/P_K(\ff)}
\mathcal{L}_{p, \ff}(\varphi, \fa^{-1})
=\sum_{\fa \in I_K(\ff)/P_K(\ff)}
\mathcal{L}_{p, \ff}(\varphi, \fa^{-1}\fp^m).
$$
where $\fa$ runs through a representative of ${\fa \in I_K(\ff)/P_K(\ff)}$ 
prime to $\ff p$. 

Let $\rho_1$ and $\rho_2$ be functions on $\Z/p^m\Z$ with values in $\C_p$. 
We identify $\cO_K/\fp^m$ with $ \Z/p^m\Z$ naturally. 
Let $\mathcal{F}_{\rho_1}^{(m)}(\alpha)$ be  the Fourier coefficient of $\rho_1$ 
with respect to the character $x \mapsto \zeta_m^{\alpha_\fp x}$, where $\alpha_\fp$ 
is the image of $\alpha$ by the natural identification $\cO_K/\fp^m \cong \Z/p^m\Z$.
In other words,
$$\mathcal{F}^{(m)}_{\rho_1}(\alpha)=
p^{-m} \sum_{x \in \Z/p^m\Z}  \rho_1(x) \zeta_m^{-\alpha_\fp x}=
p^{-m} \sum_{x \in \Z/p^m\Z}  \rho_1(x) 
\pair{1, -x \epsilon \alpha}_{\fp^m}.$$

\begin{proposition}\label{proposition; katzinterpolation}
Let $\fa$ be an integral ideal of $K$ prime to $p$ and 
let $\alpha_0$ be an element of $(\fa^{-1}\fp^m)\cap P_K(\ff)$. 
For functions $\rho_1$ and $\rho_2$ on $\Z/p^m\Z$, we have 
\begin{multline*}
	\frac{N(\fa\fp^{-m})^{k-1}N(\ff)^{l}}{\Omega_p^{k+l-1} }
	\int_{\Z_p^2} \rho_1(N(\fa) x) \rho_2(N(\ff)y) 
		x^{k-1} y^l  d \mu_{\epsilon\alpha_0\Omega,0}(x,y; \fa^{-1}\fp^{m}\Gamma)
		\\=
		\left[ \frac{2 \pi}{\sqrt{d_K}} \right]^l
		\frac{(-1)^{k+l-1} (k-1)! }{\Omega^{k+l}}
		\left. \left( \sum_{\alpha \in \fa^{-1}\cap P_K(\ff)} 
		\frac{\mathcal{F}^{(m)}_{\rho_1}(a\alpha) {\rho}_2(a \ol\alpha) \ol{\alpha}^l}{\alpha^{k} |\alpha|^{2s}}
		\right) \right|_{s=0}
	\end{multline*}
	where $a \in \Z$ is an element such that $N\fa | a $ and $a \equiv 1 \mod p^m$.  
\end{proposition}

\begin{proof}
Since both sides are bilinear for variables $\rho_1$ and $\rho_2$,  
it suffices to prove the theorem when $\rho_1$ and $\rho_2$ are given by 
$\rho_1(x)=\pair{1,  c \epsilon x}_{\fp^m}$ and 
$\rho_2(y)=\pair{1, d \epsilon y}_{\fp^m}$ for some 
$c, d \in \Z $. Then $\mathcal{F}_{\rho_1}(\alpha)$ is equal to $1$ if $\alpha \equiv c  
\mod \fp^m$ and 
$0$ otherwise. 
Therefore for $c'=N(\fa)c$ and $d'=N(\ff)d$, we calculate 
\begin{equation}\label{integralcalculation}
\int_{ \Z_p \times \Z_p} \zeta_m^{c'x+d'y} x^{k-1}y^l \;d\mu_{\epsilon\alpha_0\Omega, 0}
(x,y ; \fa^{-1}\fp^m\Gamma).
\end{equation}
We consider a commutative diagram 
\begin{equation*}
\begin{CD}
\C/\Gamma^{\sigma} @>\pi_\sigma>>  E^{\sigma}(\C)\\
@V\times \Lambda(\fa)^{\sigma}VV @VV\wt \varphi(\fa)^{\sigma} V \\
\C/\Gamma^{\tau} @>\pi_\tau>>  E^{\tau}(\C)\\
\end{CD}
\end{equation*}
where $\sigma=\sigma_{\fp^{-m}}$ and $\tau=\sigma_{\fp^{-m}\fa}$.
We  have 
$\eta_{\wh E^\tau}(t_{\fp^m}^{\sigma_\fa})=\zeta_m^{\sigma_\fa}=\zeta_m^{N\!\fa}$. 
Hence for $\delta \in \Z$ such that $\delta N(\fa) \equiv 1 \mod p^m$, the image of 
 $\delta\Lambda(\fa\fp^{-m}) \nu_m\Omega$  by 
$\eta_{\wh E^\tau} \circ \pi_\tau$ goes to $\zeta_m=S_m+1$. 
Hence noting that 
\begin{multline*}
	\left. 
		 \partial_{S, \log}^{k-1}\partial_{T, \log}^l
		\wh\Theta^{* \iota}_{\Lambda(\fa\fp^{-m})\epsilon\alpha_0\Omega,0}
		(S\oplus[c']S_m, T\oplus [d']S_m; \Gamma^\tau) \right |_{S=T=0}\\
		= \left.  (\Omega_p^\tau)^{k+l-1}
		 \partial_{z}^{k-1}\partial_{w}^l
		\Theta^{*}_{\Lambda(\fa\fp^{-m})\epsilon\Omega (\alpha_0+
	 \delta c' \nu_m), \Lambda(\fa\fp^{-m})\epsilon \delta d' \nu_m \Omega}
		(z, w; \Gamma^\tau) \right |_{z=w=0},
\end{multline*}
(\ref{integralcalculation}) is equal to 
\begin{multline*}
		\left. 
		 \Lambda(\fa\fp^{-m})
		 \partial_{S, \log}^{k-1}\partial_{T, \log}^l
		\wh\Theta^{* \iota}_{\Lambda(\fa\fp^{-m})\epsilon\alpha_0\Omega,0}(S\oplus[c']S_m, T\oplus [d']S_m; \Gamma^\tau) \right |_{S=T=0}\\
		={C}{A(\cO_K\Omega)^{-l}}
		e^*_{l,k}(\epsilon\alpha_0\Omega+ \epsilon c'\delta \nu_m\Omega, \epsilon d' \delta \nu_m\Omega; 
		\fa^{-1}\fp^m\Gamma)\\
		=
		{C}\Omega^{-k-l}\left[\frac{2\pi}{\sqrt{d_K}}\right]^l
		\left.
		\sum_{\alpha \in \fa^{-1}\cap P_K(\ff), \alpha \equiv c \;\mathrm{mod}\, \fp^m}
		\frac{\ol\alpha^l}{\alpha^k|\alpha|^{2s}}\pair{a\alpha, \epsilon d}_{\fp^m} \right|_{s=0}
	\end{multline*}
where $C=(-1)^{k+l-1}(k-1)!N(\fa^{-1}\fp^m)^{k-1}N(\ff)^{-l} \Omega_p^{k+l-1}$ and 
we may take $a=\delta N(\fa)$.
\end{proof}

\begin{theorem}\label{theorem: Katz}
Let $\varphi$ be a Hecke character of $K$ of infinity type $(k,-l)$ ($k>l \geq 0$) whose 
conductor divides $\ff p^\infty$. 
Let $\chi_p$  be  a (unique) character of $(\cO_K/p^m)^\times$ such that 
$\varphi((\alpha))=\chi_p(\alpha) \alpha^k \ol\alpha^{-l}$ for $\alpha  \in P_K(\ff)$ prime to $p$. 
We write $\chi_p$ as $\chi_p=\chi_\fp \chi_{\ol\fp}$ by 
a character $\chi_\fp$ on  $(\cO_K/\fp^m)^\times$ 
 and 
a character $\chi_{\ol\fp}$ on  $(\cO_K/\ol\fp^m)^\times$.
We put $\wt\chi_\fp=\chi_\fp \circ N_{K/\Q}$ on $I_K(p)$.  Then 
$\varphi\wt\chi_\fp^{-1}$ can naturally be extend to a character on $I_K(\ol \fp)$ as follows. 
For an ideal $\mathfrak{c} \in I_K(\ol\fp)$, we 
take a fractional ideal $\fa$ of $K$ prime to $p$ such that $\mathfrak{c}=\fa c$ for some $c \in P_K(\ff)$
and put 
$$\varphi\wt\chi_\fp^{-1}(\mathfrak{c}):=  \varphi\wt\chi_\fp^{-1}(\fa)  
\chi_{\ol\fp}(c) \chi_\fp^{-1}(\ol c)c^k \ol c^{-l}.$$
The element $\ol c$ is prime to $\fp$ and the right hand side  
is well-defined.
This value does not depend on the choices of $\fa$ and $c$. 
Let $\epsilon$ be an element of $\cO_K$ such that 
$\epsilon \equiv 1 \mod \fp^m$ and  $\epsilon \equiv 0 \mod \ol\fp^m$. 
Then the $p$-adic $L$-function has the following interpolation property. 
\begin{align*}
{\Omega_p^{-k-l}}\mathcal{L}_{p,\ff}(\varphi)={(-1)^{k+l-1} (k-1)! }\left[ \frac{2 \pi}{\sqrt{d_K}} \right]^l
                 \left(1-\frac{\varphi(\fp)}{p}\right)
		\frac{ \tau_{\varphi}(\chi_\fp)L_{\ff p}(\varphi^{-1}, 0)}{
		\Omega^{k+l}}
			\end{align*}
			where the Gauss sum $\tau_{\varphi}(\chi_\fp)$ is defined by 
			$$\tau_{\varphi}(\chi_\fp):=
			\frac{\varphi  \wt\chi_{\fp}^{-1}(\fp^n)}{p^{n}} 
			\sum_{a \in (\Z/p^n\Z)^\times} \chi_\fp(a) \pair{1, -a\epsilon}_{\fp^n}
			$$
			if $\chi_\fp \not=1$ where $\fp^n$ is the conductor of $\chi_\fp$, and 
			we put $\tau_{\varphi}(\chi_\fp)=1$ if $\chi_\fp=1$. 
\end{theorem}
\begin{proof}
We may replace $\alpha_0$ by $\epsilon\alpha_0$ in the definition of 
the partial $p$-adic $L$-function. 
Let $\mathbf{1}_{\Z_p^\times}$ be the characteristic function of $\Z_p^\times$, 
namely, $\mathbf{1}_{\Z_p^\times}(x)=1$ if $x \in \Z_p^\times$ and $0$ otherwise. 
We apply the previous proposition for $\rho_1=\chi_1 \cdot \mathbf{1}_{\Z_p^\times}$ and 
$\rho_2=\chi_2 \cdot \mathbf{1}_{\Z_p^\times}$.  We use the same notations 
as in the previous proposition. 
First we assume that $\chi_1$ is not trivial. Then $\chi_1 \cdot \mathbf{1}_{\Z_p^\times}=\chi_1$. 
If $v_\fp(\alpha)\not=m-n$ for $\alpha \in \fa^{-1}$, 
then $\pair{1, x\epsilon a \alpha }_{\fp^m}$  is not a $p^n$-th power 
root of unity and the Gauss sum $\mathcal{F}_{\chi_1}(\alpha)$ is equal to $0$.  
Suppose that $v_\fp(\alpha)=m-n$. 
Let $\fa_0$ be a fractional ideal of $K$ prime to $p$ such that 
 $\fp^{m-n}=c\fa_0$ for some $c \in P_K(\ff)$ and we put $\beta=\alpha c^{-1}$. 
Then $\beta$ is a $\fp$-adic unit and $\beta \in \fa_0\fa^{-1}$, 
and we have  
\begin{align*}
	\mathcal{F}^{(m)}_{\chi_1}(a\alpha)&=p^{-m} \sum_{x \in (\Z/p^m)^\times}  \chi_1( x) 
	\pair{1, -x \epsilon a\alpha}_{\fp^m} \\
	&=p^{-m}   \sum_{x \in (\Z/p^m)^\times}  \chi_1(\ol{c}N\!\fa_0 x) 
		\pair{1, -\ol{c} N\!\fa_0 x \epsilon c  a\beta}_{\fp^m}\\
	&=\frac{ \chi_1(\ol{c}N\!\fa_0 )}{p^n}\sum_{x \in (\Z/p^n)^\times}  \chi_1( x) 
		\pair{1, - x \epsilon a \beta}_{\fp^n} =\chi_1(\ol{c}N\!\fa_0 )\tau(\chi_\fp)\chi_1^{-1}(\beta).
\end{align*} 
Note that we have
\begin{multline*}
	 \chi_\fp(\ol{c}N\!\fa_0 )\varphi^{-1}(\fa)  
		(\varphi^{-1}\wt\chi_\fp)(\fp^{-m}) 
		 \left. \left( \sum_{\substack{\alpha \in \fa^{-1}\cap P_K(\ff)\cap I_K(\ol \fp),\\ v_{\fp}(\alpha)=m-n}} 
		\frac{\chi_\fp^{-1}(\beta) {\chi}_{\ol \fp}^{-1}(\alpha) \ol{\alpha}^l}{\alpha^{k} |\alpha|^{2s}}
		\right) \right|_{s=0} \\
		= \varphi\wt\chi_\fp^{-1}(\fp^n)
		\left. \left( \sum_{\beta \in \fa_0\fa^{-1}\cap P_K(\ff)\cap I_K(p)} 
		\frac{\varphi^{-1}(\fa\fa_0^{-1} \beta)}{|\beta|^{2s}}
		\right) \right|_{s=0}. 
\end{multline*}
Hence by Proposition \ref{proposition; katzinterpolation}, we have
\begin{align*}
	& \mathcal{L}_{p,\ff}(\varphi, \fa^{-1}\fp^m) \\
	&= \Omega_p
 (\varphi^{-1}\wt\chi_1)(\fa\fp^{-m})  N(\fa\fp^{-m})^{k-1}N(\ff)^{l} \\
& \qquad \qquad \qquad \times \int_{ \Z_p^\times \times \Z_p^\times} \chi_1(x)\chi_2(N(\ff)y) x^{k-1}y^l \;d\mu_{\alpha_0\Omega, 0}(x,y ; \fa^{-1}\Gamma) \\	
	&={(-1)^{k+l-1} (k-1)! }
	\left[ \frac{2 \pi}{\sqrt{d_K}} \right]^l		
		 \tau(\chi_\fp)\\
                  &  \qquad \qquad \qquad \qquad  \times \varphi\wt\chi_\fp^{-1}(\fp^n)
		\left. \left( \sum_{\beta \in \fa_0\fa^{-1}\cap P_K(\ff)\cap I_K(p)} 
		\frac{\varphi^{-1}(\fa\fa_0^{-1} \beta)}{|\beta|^{2s}}
		\right) \right|_{s=0}. 
	\end{align*}
	Hence we have 
	$$\sum_{\fa \in I_K(\ff)/P_K(\ff)} 
   \frac{\mathcal{L}_{p,\ff}(\varphi, \fa^{-1}\fp^m)}{\Omega_p^{k+l}}
	=(-1)^{k+l-1} (k-1)! 
	\left[ \frac{2 \pi}{\sqrt{d_K}} \right]^l		
		 \frac{\tau_\varphi(\chi_\fp)L_{\ff p}(\varphi^{-1}, 0)}{\Omega^{k+l}}.$$
	Suppose that $\chi_1=1$. 
	Then  $\mathcal{F}^{(m)}_{\chi_1 \cdot \mathbf{1}_{\Z_p^\times}}(a\alpha)$ for 
	$\alpha \in \fa^{-1}$ is 
	equal to $(p-1)/p$ if $v_{\fp}(\alpha)=m$,  
	to $-1/p$ if $v_{\fp}(\alpha)=m-1$ and to $0$ otherwise. 
	It is also straightforward to see that $ \varphi\wt\chi_\fp^{-1}(\fp)= \varphi(\fp)$. 
	Hence Proposition \ref{proposition; katzinterpolation} shows that for
	$\mathfrak{c} :=  \fa^{-1}\cap P_K(\ff) \cap I_K(\ol \fp)$, we have
\begin{align*}
&\mathcal{L}_{p,\ff} (\varphi, \fa^{-1}\fp^m)/{(-1)^{k+l-1} (k-1)! }\Omega^{-k-l}
\left[ \frac{2 \pi}{\sqrt{d_K}} \right]^l		 \\
                  & =\left.\frac{\varphi(\fp^m)}{p}
		 \left( (p-1)\sum_{\substack{\alpha \in \mathfrak{c},\\ v_\fp(\alpha)=m} }
		\frac{\varphi^{-1}(\fa \alpha)}{|\alpha|^{2s}}-
				 \sum_{\substack{\alpha \in \mathfrak{c},\\ v_\fp(\alpha)=m-1} }
		\frac{\varphi^{-1}(\fa \alpha)}{|\alpha|^{2s}}
		\right) \right|_{s=0}\\
		 & = \left.\left(
		\sum_{\substack{\alpha \in  \mathfrak{c}, \\v_\fp(\alpha)=m} }
		\frac{\varphi^{-1}(\fp^{-m} \fa \alpha)}{|\alpha|^{2s}} \right)
		-\frac{\varphi(\fp)}{p}
		\left(
				 \sum_{\substack{\alpha \in \mathfrak{c},\\ v_\fp(\alpha)=m-1}} 
		\frac{\varphi^{-1}(\fp^{-m+1} \fa \alpha)}{|\alpha|^{2s}}
		 \right) \right|_{s=0} \\
		 & =\left(1-\frac{\varphi(\fp)}{p} \right) \left.\left(
		\sum_{\alpha \in  \mathfrak{c} \cap I_K(\fp)} 
		\frac{\varphi^{-1}(\fa \alpha)}{|\alpha|^{2s}} \right)
		 \right|_{s=0}. 
	\end{align*}
	Hence we have the desired formula. 
	\end{proof}

Now we consider a modular elliptic curve $E$ over $\Q$, and we
compare our measure to the Mazur and Swinnerton-Dyer measure. 
Our approach gives an explicit  expression of the Mazur and Swinnerton-Dyer measure 
in terms of coefficients of our theta function.

\begin{corollary}[Katz]
We fix a N\'eron differential $\omega_E$.  
Let $\Omega_E^+$ be a real N\'eron period of $E$ and 
let $\Omega$ be a complex number such that $\Gamma=\ff \Omega$.  
Let $u$ be a $\fp$-adic unit in $\cO_K$ such that 
$\Omega^+_E=u\Omega$. 
Let $\psi$ be the Gr\"ossencharakter over $K$ associated to $E$ and 
 let $\chi: (\bbZ/p^m \bbZ)^\times \rightarrow \mu_{p^\infty}$ be a Dirichlet character 
of conductor $p^m$. 
We regard $\chi$ as a character on $\cO_K$ by the norm map.
Then 
we have 
$$\Omega_p^{-1}\mathcal{L}_{p, \ff}(\psi\chi)=u \mathcal{L}_{p, \ff}^{MS}(\chi)$$
where $\mathcal{L}_{p, \ff}^{MS}$ is the Mazur and Swinnerton-Dyer $p$-adic $L$-function 
 for the system of roof of unity $(\zeta_{p^m})_m$ such that 
$\zeta_{p^m}=\pair{1, -\alpha^m\epsilon}_{\fp^m}$ where 
$\epsilon$ is any element such that 
$\epsilon\equiv 1 \mod \fp^m$ and $\epsilon\equiv 0 \mod \ol\fp^m$.  
\end{corollary}
\begin{proof}
Let $L(E/\bbQ, s)$ be the Hasse-Weil $L$-function of $E$ over $\Q$, and 
we let $L(E/\bbQ, \chi, s)$ be its $\chi$-twist. 
By a theorem of Deuring, we have
\begin{equation}\label{equation: Deuring}
	L(E/\bbQ, \ol\chi, s) =L(\ol{\psi\chi}, s)= L((\psi\chi)^{-1}, s-1). 
\end{equation}
For a generator $(\zeta_{p^m})_m$ of $\Z_p(1)$,  
Mazur and Swinnerton-Dyer 
(see also \cite{MTT}) 
defined a measure $\nu$ on $\bbZ_p^\times$ such that
\begin{align*}
	\mathcal{L}_{p, \ff}^{MS}(\chi):=\int_{\bbZ_p^\times}& \chi(u) d \nu(u) = 
	\frac{p^m}{\alpha^m}
	\frac{ L(E/\bbQ, \ol\chi, 1) }
	{\tau(\ol\chi)  \Omega_E^+},%
\end{align*}
where 
$\tau(\chi) : = \sum_{a \in (\bbZ/p^m\Z)^\times} \chi(a) \zeta_{p^m}^a$ 
and $\alpha$ is a unit root of the Euler factor $X^2-a_pX+p$. 
In our case, $\alpha$ is equal to $\ol{\psi}(\fp)$. 

By Theorem \ref{theorem: Katz} and (\ref{equation: Deuring}),  we have 
$$\frac{\mathcal{L}_{p, \ff}(\psi\chi)}{\Omega_p}
=\frac{\tau(\chi) \psi(\fp)^m}{p^m}
	\frac{ L(E/\bbQ, \ol\chi, 1) }
	{  \Omega}=
	\frac{p^m}{\ol\psi(\fp)^m}
	\frac{ L(E/\bbQ, \ol\chi, 1) }
	{\tau(\ol\chi)  \Omega},
$$ 
where $\zeta_{p^m}=\pair{1, \epsilon}_{\fp^m}$. 
(Note that $\psi\chi\chi^{-1}_{\fp}(\fp)=\psi(\fp)$). 
Our assertion now follows from this calculation.
\end{proof}

It is also simple to give the measure of Mazur and Swinnerton-Dyer 
using the coefficients of our theta functions. 

\begin{proposition}
Suppose that 
$$\psi(\alpha) = \varepsilon(\alpha) \alpha, \qquad (\alpha, \ff) = 1$$
for some finite character $\varepsilon: (\cO_K/\ff)^\times \rightarrow \ol\bbQ^\times$. \\
i) For  $\alpha \in \cO_K\setminus\ff$,  we let  
$$
\Theta_{ \alpha \Omega,  0}(z, N(\ff)w; \ff \Omega)-\frac{1}{w}=
\sum_{k,l} c_{k,\ell}(\alpha) \frac{z^kw^{\ell}}{k! \ell!}. 
$$
We put 
$$ c_{k, \ell}:=  \Omega_p^{k+\ell}\sum_{\alpha \in  (\cO_K/\ff)^\times} \varepsilon(\alpha) c_{k,\ell}(\alpha).$$
Then 
the $p$-adic limit $c_k:=\lim_{\ell \rightarrow -k}  c_{k, \ell}$ 
 exists in $W$ where $\ell$ runs through positive integers divisible by $(p-1)$ which goes to $-k$ 
with the $p$-adic topology. \\
ii) We put $f(z)=\sum_{n=0}^\infty c_n z^n/n!$ and 
$f_{\ell}(z)=\sum_{n=0}^\infty c_{n, n \ell} z^n/n!$. 
Then 
$$
\wh f_{\ell}(T):=f_{\ell}(\Omega_p \log (1+T))
\in W[[T]].$$
In particular,  $\wh f(T):=f(\Omega_p \log (1+T)) \in W[[T]]$. \\
iii) Let $\mu$ be the measure corresponding to the power series 
$\wh f(T)-\wh f([p]T)$. 
Then we have 
$$\Omega_p^{-1}{\mathcal{L}_{p, \ff}(\psi\chi)}=
\int_{\Z_p^\times} \chi(x) d\mu(x). 
$$
Namely, $\mu$ 
is the measure of Mazur and Swinnerton-Dyer. 
$($Since $c_{k, \ell}(\alpha)$ is essentially $e^*_{\ell, k+1}(\alpha, 0)$, 
the number $c_{k}$ is a $p$-adic version of the special value $L(E, k+1)$ 
and $f(z)$ is the generating function for these values. $)$
\end{proposition}
\begin{proof}
i) follows from  
$$\sum_{a \in  (\cO_K/\ff)^\times} \varepsilon(\alpha) c_{k, \ell}(\alpha)=\Omega_p^{k+\ell}\sum_{a \in  (\cO_K/\ff)^\times} \varepsilon(\alpha) \int_{\Z_p^\times \times \Z_p^\times} x^k y^{\ell} 
d\wt \mu_{\alpha\Omega,0}(x,y; \ff\Omega).$$
ii) follows from  
$$\Omega_p^{k+\ell}\sum_{a \in  (\cO_K/\ff)^\times} \varepsilon(\alpha) \int_{\Z_p^\times \times \Z_p^\times} 
\binom{x y^{\ell}}{n} d\wt\mu_{\alpha\Omega,0}(x,y; \ff\Omega) \in W.$$
Now we prove iii). 
Since the class number of $K$ is equal to $1$ and 
$w_\ff=1$ by the existence of $\psi$, we have $I(\ff)/P(\ff) \cong (\cO_K/\ff)^\times$.
Hence we can choose the representative of $(\cO_K/\ff)^\times$ as 
the representative of $I(\ff)/P(\ff)$. 
Then by (\ref{scalar change}) we have 
\begin{align*}
\mathcal{L}_{p, \ff}(\psi\chi)&=
\sum_{a \in  (\cO_K/\ff)^\times}\varepsilon(a)  \int_{(\cO_K \otimes \Z_p)^\times} 
\chi(\alpha) \alpha 
 \;d\mu_{a\Omega, 0}^{\mathscr{K}}(\alpha ; \ff \Omega) \\
&={\Omega_p}\sum_{a \in  (\cO_K/\ff)^\times} \varepsilon(a) 
\int_{\Z_p^{\times} \times \Z_p^{\times}} 
 \chi(xy^{-1}) d\wt \mu_{a\Omega, 0}(x,y ; \ff \Omega). 
\end{align*}
By Proposition \ref{proposition: M2}, for $\pi=\psi(\fp)$ we have 
	\begin{multline*}
                \sum_{\alpha \in (\cO_K/\ff)^\times} \varepsilon(\alpha)  
			\int_{\bbZ_p^\times \times \bbZ_p^\times}   x^k y^{\ell} 
			d\wt \mu_{ \alpha \Omega, 0}(x,y)\\
			=\Omega_p^{k+\ell}  
			 \left(
			1- \frac{\pi^k }{\ol \pi^{\ell} \varepsilon(\pi)} 
			 \right)
			 \left(
			1
			 -\frac{\varepsilon(\ol\pi )\pi^{\ell}}{\ol \pi^{k}} 
			 \right)
			 \sum_{\alpha \in (\cO_K/\ff)^\times} \varepsilon(\alpha) 
			c_{k, \ell}(\alpha).
	\end{multline*}
Hence by taking $p$-adic  limit $\ell \rightarrow -k$ for  a positive $\ell$ such that $(p-1)|\ell$, 
	\begin{equation*}
                \sum_{\alpha \in (\cO_K/\ff)^\times} \varepsilon(a)  
			\int_{\bbZ_p^\times \times \bbZ_p^\times}   x^k y^{-k} 
			d\wt \mu_{ \alpha \Omega, 0}(x,y)
			=
			 (
			1- {p^k } 
			 )
			c_{k}
	\end{equation*}
	(Note that $\pi^\ell \rightarrow 0$ and $\varepsilon(\pi)=1$.) 
	Therefore if we write as 
	$\chi(x)=\sum_{k=0}^\infty a_k \binom {x}{k} $ with $a_k \rightarrow 0$, we have 
	\begin{align*}
	\Omega_p^{-1} \mathcal{L}_{p, \ff}(\psi\chi)
&=\sum_{k=0}^\infty a_k \sum_{a \in  (\cO_K/\ff)^\times} \varepsilon(a) 
\int_{\Z_p^{\times} \times \Z_p^{\times}} 
 \binom{xy^{-1}}{k} d\wt \mu_{a\Omega, 0}(x,y ; \ff \Omega) \\
&=\sum_{k=0}^\infty a_k  
\int_{\Z_p } 
 \binom{x}{k} d \mu(x)= \int_{\Z_p} \chi(x) d\mu(x). 
	\end{align*}
\end{proof}

\section{Explicit Calculation for the Supersingular Case}\label{SupersingularCase}
%
%

In this section, we consider the case when $p$ is a supersingular prime.  
In the ordinary case, the $p$-adic measure which defines  the  two variable $p$-adic $L$-function 
is  constructed from the power series  
$$
	\wh{\Theta}_{z_0, w_0}(s, t) := \left.\Theta_{z_0, w_0}(z,w)\right|_{z = \lambda(s), w = \lambda(t)}
$$
and the existence of  the $p$-adic $L$-function is equivalent to 
the integrality of this power series. 
This power series is also well-defined in the supersingular case, and
contains information about the $p$-adic properties of  Eisenstein-Kronecker numbers.
Since Mumford's translation works over integral basis, 
we do not lose generality if we only consider the case $(z_0, w_0)=(0,0)$, 
namely, the power series $\wh{\Theta}(s, t)$. 
By Kronecker's formula, we have $\Theta(z,w)=\theta(z+w)/\theta(z)\theta(w)$.
Hence it would be important to investigate the $p$-adic properties of $\wh\theta(t)$.
For example, $\Theta(z,z)=\theta(2z)/\theta(z)^2$ is equal to 
$\theta(z)^2$ up to an elliptic function.  Since $p$ is odd, its 
$p$-adic  properties are essentially equivalent to that of $\theta(z)$.

\begin{figure}\caption{CM elliptic curves over $\Q$.}\label{fig: one}
$$
\begin{array}{|c||c|c||c|}
\hline
\mathrm{End}(E_{\mathbb{C}}) & g_2 & g_3 & e_2^{*} \\
\hline
\hline
\Z[\frac{1+\sqrt{-3}}{2}] & 0 & u & 0 \\
\Z[\sqrt{-3}] & 15u^2 & 11u^3 & \frac{u}2 \\
\Z[\frac{1+3\sqrt{-3}}{2}] & 120u^2 & 253u^3 & 2u \\
\hline
\Z[\sqrt{-1}] & u & 0 & 0 \\
\Z[2\sqrt{-1}] & 44u^2 & 56u^3 & u \\
\hline
\Z[\frac{1+\sqrt{-7}}{2}] & 35u^2 & 49u^3 & \frac{u}2 \\
\Z[\sqrt{-7}] & 5\cdot 7 \cdot 17u^2 
& 3\cdot 7^2 \cdot 19u^3 & \frac{9u}2 \\
\hline
\Z[\sqrt{-2}] & 30u^2 & 28u^3 & \frac{u}2 \\
\hline
\Z[\frac{1+\sqrt{-11}}{2}] & 8\cdot 3
\cdot 11u^2 & 7\cdot 11^2 u^3 & 2u \\
\hline
\Z[\frac{1+\sqrt{-19}}{2}] & 8\cdot 19u^2 & 19^2u^3 & 2u \\
\hline
\Z[\frac{1+\sqrt{-43}}{2}] & 16\cdot 5\cdot 43u^2 & 
3\cdot 7\cdot 43^2u^3 & 12u \\
\hline
\Z[\frac{1+\sqrt{-67}}{2}] & 8\cdot 5\cdot 11\cdot 67u^2 & 
7\cdot 31\cdot 67^2u^3 & 38u \\
\hline
\Z[\frac{1+\sqrt{-163}}{2}] & 
16\cdot 5 \cdot 23\cdot 29\cdot 163u^2 & 
7\cdot 11\cdot 19\cdot 127\cdot 163^2 u^3 & 724u \\
\hline
\end{array}
$$
\end{figure}

Using Mathematica, we explicitly calculated the coefficients of the power series $\wh\Theta(s,t)$
for a supersingular $p$.  The table of Figure \ref{fig: one} classifying CM elliptic curves 
$$
	E: y^2=4x^3-g_2x-g_3
$$ 
defined over $\Q$,   wans kindly given to us by Seidai Yasuda and was helpful for our calculations.  
The $u$ in the table is any non-zero rational number.  The table also
contains the values for $e_2^*$, and was used to find examples for various $e_2^*$.
Yasuda also notified us the expansion of the formal logarithm $\lambda(t)$ of the above elliptic curve
for the parameter $t= - 2x/y$.  It is given by 
$$
	\lambda(t)=\sum_{m, n=0}^\infty \, 
	\frac{(2m+3n)!}{(m+2n)! m!n!}
	 \left(-\frac{g_2}{4}\right)^m \left(-\frac {g_3}{4}\right)^n \frac{ t^{4m+6n+1}}{4m+6n+1}.
$$
We will give a proof of this formula in \cite{BK}, Appendix. 
We calculated the power series $\wh\Theta(s,t)$ using these formulas.
The calculation and graphs in this section were created using Mathematica.

\begin{figure}
	\includegraphics[scale=1.0]{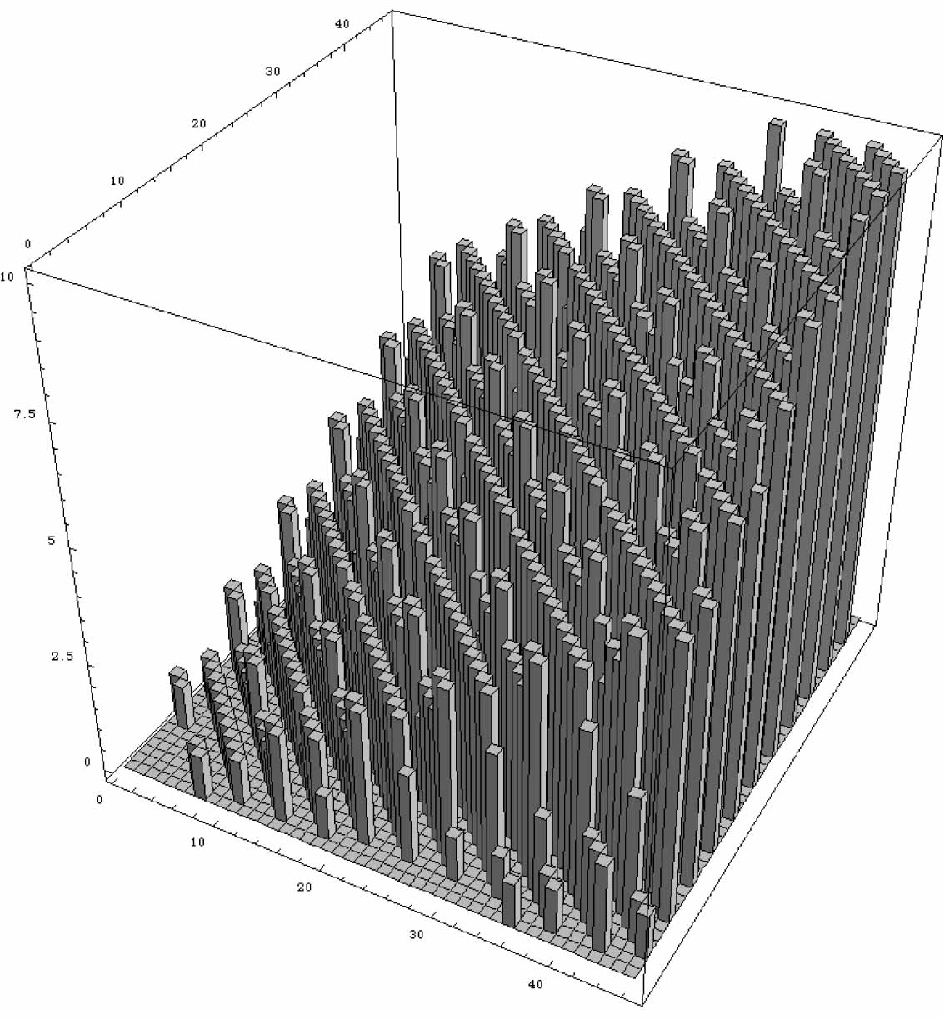}
	\caption{Exponent of $p=7$ in denominator of $\wh c_{m,n}$.}
	\label{figure: two}
\end{figure}

For example, consider the elliptic curve  $E: y^2=4 x^3 - 4x$.  Then this curve has good reduction for $p \geq 5$.
Since $E$ has complex multiplication in $\bbZ[i]$,  it has good ordinary reduction for 
$p \equiv 1 \pmod{4}$ and good supersingular reduction for $p\equiv 3 \pmod{4}$.  
Consider the power series
$$
	\wh{\Theta}(s, t) = \frac{1}{s} + \frac{1}{t} + \sum_{m, n \geq 0} \wh c_{m,n} s^m t^n
$$
and prime $p=7$ where $E$ has good supersingular reduction.
Figure \ref{figure: two} is a graph of the exponent of $p=7$ in the denominator of $\wh c_{m,n}$.
We observe that the growth of the exponent of $p=7$ in the denominator of $\wh c_{m, n}$
increases greatly in the diagonal $m=n$ direction.  

\begin{figure}
	\begin{psfrags}
		\psfrag{A1}{$\wh c_{m,0}$}
		\psfrag{A2}{$\wh c_{m,4}$}
   		\psfrag{A3}{$\wh c_{m,8}$}
    		\psfrag{A4}{$\wh c_{m,m}$}
		\psfrag{A5}{$pn/(p^2-1)$}
		\includegraphics[scale=1.0]{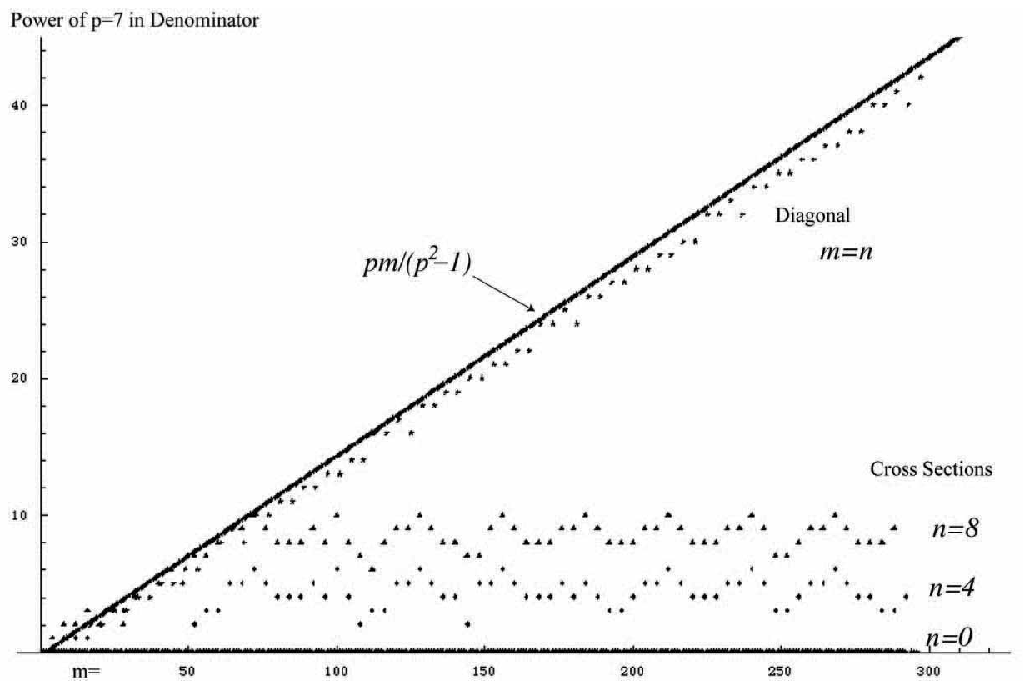}
	\end{psfrags}
	\caption{Cross Sections of Figure \ref{figure: two}}
	\label{figure: three}
\end{figure}

Figure \ref{figure: three} represents the various cross sections of Figure \ref{figure: two}.
Note that for a constant $n$, the increase of the exponent of $p=7$ in the denominator of $\wh c_{m,n}$ 
levels off to a certain pace, but exponent in the denominator of $\wh c_{m,m}$ rises at a consistent rate.
The slope $p/(p^2-1)$ is  the valuation of the $p$-adic period  of $E$. Although we have explained the example 
of a specific elliptic curve for a specific supersingular prime, we have been able to check that this phenomenon holds 
for any CM elliptic curve over $\bbQ$ and supersingular prime $p$ that we have tried to calculate.

This explains why it was possible to construct $p$-adic $L$-functions for supersingular 
elliptic curves in one variable (\cite{Box1}, \cite{Box2}, \cite{ST}, \cite{Fou} and \cite{Yam}), but not in two-variables.
If we consider the two-variable power series $\wh\Theta(s,t)$, then the above calculation would
lead us to expect the radius of convergence is equal to $p^{-p/(p^2-1)} < 1$.  Since the general 
theory of $p$-adic Fourier transforms constructed by Schneider-Teitelbaum \cite{ST} calls for a
power series with radius of convergence $1$, we are not able to apply the general theory in order
to construct the $p$-adic measure in this case.  

In a subsequent paper \cite{BK}, we proceed with the investigation of $\wh\Theta(s,t)$ for a supersingular prime
$p \geq 5$.  We show that the $p$-adic radius of convergence of $\widehat{\theta}(t)$ is precisely equal to 
$p^{- p/(p^2-1)}$, which is the valuation of the $p$-adic period.  This implies that the radius of convergence of
$\wh\Theta(s,t)$ is $p^{- p/(p^2-1)}$.  We then prove a refined version of the $p$-adic Fourier transform of Schneider-Teitelbaum.  
Using this theory, we were able to derive precise divisibility results about 
the critical values of the two-variable Hecke $L$-function.  This generalizes one-variable 
congruence and divisibility results by Katz \cite{Ka6} and  Chellali \cite{Ch}, and improves upon the 
two-variable divisibility results by Fujiwara \cite{Fuj}.  See \cite{BK} for details.  

It seems that the radius of convergence of $p$-adic theta functions 
are closely related to the valuation of the $p$-adic periods. 
M. Kurihara pointed out to the second author that 
the $\mu$-invariant of the cyclotomic $p$-adic $L$-function 
at supersingular primes is also conjectured to be equal to  $p/(p^2-1)$ 
(the valuation of the $p$-adic period in the supersingular case). 
(See \cite{Ku}, Remark 0.2 or \cite{Ko}, Corollary 10.10).
On the other hand,  under mild assumptions, 
the $\mu$-invariant in the ordinary case is conjectured to be $0$, 
which is also the valuation of the $p$-adic period in this case. 
We do not know the reason of these coincidences.



\end{document}